\documentclass{article}%
\usepackage{amssymb}
\usepackage{amsmath}
\usepackage{makeidx}
\usepackage{cite}
\usepackage{amsfonts}
\usepackage{graphicx}%
\begin{document}

\author{Hartmut Wachter\thanks{E-Mail: Hartmut.Wachter@gmx.de}\\An der Schafscheuer 56\\D-91781 Wei\ss enburg, Federal Republic of Germany}
\title{Tutorial on one-dimensional $q$-Fourier transforms}
\maketitle
\date{}

\begin{abstract}
This paper is an introductory text to the theory of $q$-de\-formed Fourier
transforms, as first discussed by Rogov and Olshanetsky. We derive the
well-known results in detail, present them in a format that suits our needs,
and include some new findings on specific aspects of the theory.

\end{abstract}
\tableofcontents

\section{Introduction}

The first three chapters review some basics of the so-called $q$-cal\-cu\-lus.
These include the Jackson derivatives and the Jackson integrals as $q$-analogs
of ordinary derivatives or integrals. The eigenfunctions of the Jackson
derivative lead us to $q$-analogs of the exponential function. Using the
$q$-ex\-po\-nen\-tial, we introduce $q$-analogs of the trigonometric functions
and derive $q$-ver\-sions of the Taylor expansion. In the subsequent chapters,
we present a $q$-ver\-sion of the Hilbert space of square-integrable functions
and study the spaces of test functions and $q$-dis\-tri\-bu\-tions. Finally,
we discuss $q$-analogs of the classical Fourier transforms of position or
momentum space.

Our approach uses the considerations from Ref.~\cite{olshanetsky1998q.alg},
adapting them to our specific requirements. We have also added some new
results relevant to our further investigations. Remember, there are other
$q$-ver\-sions of Fourier transforms apart from the ones discussed in this
paper \cite{Askey1993,1996q.alg.....8008K}. The $q$-de\-formed Fourier
transforms used in this paper are the best option for our needs because they
closely resemble the classical Fourier transform.

\section{Exponential functions\label{KapqAnaExpTriFkt}}

The \textbf{Jackson derivative} of a real function $f$ is defined as follows
\cite{Jackson:1910yd}:%
\begin{equation}
D_{q}f(x)=\frac{f(x)-f(qx)}{x-qx}. \label{DefJacAbl}%
\end{equation}
In the limit $q\rightarrow1$, the Jackson derivative becomes the ordinary
derivative of a real function \cite{Klimyk:1997eb}:%
\begin{equation}
\lim_{q\hspace{0.01in}\rightarrow1}D_{q}f(x)=\frac{\text{d}f(x)}{\text{d}x}.
\label{KlaGreJacWdh}%
\end{equation}

The definition in Eq.~(\ref{DefJacAbl}) implies the following rule for the
Jackson derivative of a power function ($\alpha\in\mathbb{R}$)%
\begin{equation}
D_{q}\hspace{0.01in}x^{\alpha}=[[\alpha]]_{q}\hspace{0.01in}x^{\alpha-1},
\end{equation}
The \textbf{antisymmetric }$q$\textbf{-num\-bers} are defined as follows:%
\begin{equation}
\lbrack\lbrack\alpha]]_{q}=\frac{1-q^{\alpha}}{1-q}. \label{DefAntSymQZah2}%
\end{equation}
For $q\rightarrow1$, the antisymmetric $q$-num\-ber becomes the corresponding
real number:%
\begin{equation}
\lim_{q\hspace{0.01in}\rightarrow1}[[\alpha]]_{q}=\alpha. \label{GreAntQZah}%
\end{equation}
The $q$\textbf{-fac\-to\-ri\-al} is introduced in analogy to the undeformed
case ($n\in\mathbb{N}$):%
\begin{equation}
\lbrack\lbrack\hspace{0.01in}n]]_{q}!=[[1]]_{q}\hspace{0.01in}[[2]]_{q}%
\ldots\lbrack\lbrack\hspace{0.01in}n-1]]_{q}\hspace{0.01in}[[\hspace
{0.01in}n]]_{q},\qquad\lbrack\lbrack0]]_{q}!=1. \label{DefQFakJacBas}%
\end{equation}
Accordingly, we define the $q$\textbf{-bi\-no\-mi\-al coefficients} as follows
($k\in\mathbb{N}$):%
\begin{equation}%
\genfrac{[}{]}{0pt}{}{n}{k}%
_{q}=\frac{[[\hspace{0.01in}n]]_{q}!}{[[\hspace{0.01in}n-k]]_{q}%
!\hspace{0.01in}[[k]]_{q}!}. \label{qBinKoeBas}%
\end{equation}

By mathematical induction, we can proof the following formula for the Jackson
derivative \cite{Klimyk:1997eb}:%
\begin{equation}
D_{q}^{n}f(x)=(-1)^{n}\frac{q^{-n(n\hspace{0.01in}-1)/2}}{(1-q)^{n}x^{n}}%
\sum_{m\hspace{0.01in}=\hspace{0.01in}0}^{n}%
\genfrac{[}{]}{0pt}{}{n}{m}%
_{q}(-1)^{m}q^{m(m\hspace{0.01in}-1)/2}f(q^{n-m}x).\label{MehAnwJacAblFkt}%
\end{equation}
Moreover, the Jackson derivative satisfies a \textbf{product rule}:%
\begin{align}
D_{q}\left(  f(x)\hspace{0.01in}g(x)\right)   &  =f(qx)\hspace{0.01in}%
D_{q}\hspace{0.01in}g(x)+g(x)\hspace{0.01in}D_{q}f(x)\nonumber\\
&  =f(x)\hspace{0.01in}D_{q}\hspace{0.01in}g(x)+g(qx)\hspace{0.01in}%
D_{q}f(x).\label{ProRegJacAblBas}%
\end{align}

We introduce the $q$\textbf{-ex\-po\-nen\-tial} as eigenfunction of the
Jackson derivative \cite{Majid:1993ud}, i.~e.%
\begin{equation}
D_{q}\exp_{q}(x)=\exp_{q}(x) \label{Bed1qExp}%
\end{equation}
with the normalization condition%
\begin{equation}
\exp_{q}(0)=1. \label{Bed2qExp}%
\end{equation}
To get a power series expansion of the $q$-ex\-po\-nen\-tial, we make the
following ansatz:%
\begin{equation}
\exp_{q}(x)=\sum_{k\hspace{0.01in}=\hspace{0.01in}0}^{\infty}\hspace
{0.01in}(C_{q})_{k}\,x^{k}.
\end{equation}
Due to Eqs.~(\ref{Bed1qExp}) and (\ref{Bed2qExp}), the coefficients
$(C_{q})_{k}$ satisfy the following recurrence relations:%
\begin{equation}
(C_{q})_{k\hspace{0.01in}+1}=\frac{1}{[[\hspace{0.01in}k+1]]_{q}}%
\hspace{0.02in}(C_{q})_{k},\qquad(C_{q})_{0}=1.
\end{equation}
By iteration, we find:%
\begin{equation}
(C_{q})_{k}=\frac{1}{[[k]]_{q}!}.
\end{equation}
Thus, the power series expansion for the $q$-ex\-po\-nen\-tial takes the
following form:%
\begin{equation}
\exp_{q}(x)=\sum_{k\hspace{0.01in}=\hspace{0.01in}0}^{\infty}\hspace
{0.01in}\frac{1}{[[k]]_{q}!}\,x^{k}. \label{AusQExpEin}%
\end{equation}
For $q\rightarrow1$, the $q$-ex\-po\-nen\-tial becomes the ordinary
exponential function [cf. Eq.~(\ref{GreAntQZah})]:%
\begin{equation}
\lim_{q\rightarrow1}\exp_{q}(x)=\operatorname*{e}\nolimits^{x}.
\label{KlaGreQExp}%
\end{equation}

The $q$-ex\-po\-nen\-tial in Eq.~(\ref{AusQExpEin}) is related to the
following functions ($0<q<1$) \cite{Gasper:1990,Klimyk:1997eb}:%
\begin{align}
\operatorname{E}_{q}(z)  &  =\sum_{k\hspace{0.01in}=\hspace{0.01in}0}^{\infty
}\frac{q^{k(k\hspace{0.01in}-1)/2}z^{k}}{(q\hspace{0.01in};q)_{k}}%
=(-z\hspace{0.01in};q)_{\infty},\nonumber\\
\operatorname{e}_{q}(z)  &  =\sum_{k\hspace{0.01in}=\hspace{0.01in}0}^{\infty
}\frac{z^{k}}{(q\hspace{0.01in};q)_{k}}=\frac{1}{(z\hspace{0.01in};q)_{\infty
}}. \label{qExpLit}%
\end{align}
In the formulae above, we have used the $q$\textbf{-Poch\-ham\-mer symbol}
($a\in\mathbb{C}$, $n\in\mathbb{N}_{0}$):%
\begin{equation}
(a\hspace{0.01in};q)_{n}=\left\{
\begin{tabular}
[c]{cl}%
$1$ & f\"{u}r\quad$n=0,$\\
$(1-a)(1-a\hspace{0.01in}q)\ldots(1-a\hspace{0.01in}q^{n-1})$ & f\"{u}r\quad
$n\geq1.$%
\end{tabular}
\ \ \ \ \right.  \label{qPocHamSym}%
\end{equation}
Accordingly, we have:%
\begin{equation}
(a\hspace{0.01in};q)_{\infty}=\lim_{n\rightarrow\hspace{0.01in}\infty
}\,(a\hspace{0.01in};q)_{n}=\prod\limits_{j=\hspace{0.01in}0}^{\infty
}(1-qa^{j}). \label{PocSymUneDef}%
\end{equation}

The power series of $\operatorname{e}_{q}(z)$ converges for $\left\vert
z\right\vert <1$ and that of $\operatorname{E}_{q}(z)$ for $z\in\mathbb{C}$.
Due to the identities%
\begin{equation}
\lbrack\lbrack\hspace{0.01in}n]]_{q}!=\frac{(q\hspace{0.01in};q)_{n}%
}{(1-q)^{n}} \label{ZusqFakPocSym}%
\end{equation}
and \cite{Gasper:1990}%
\begin{equation}
(a\hspace{0.01in};q)_{n}=(a^{-1};q^{-1})_{n}\hspace{0.01in}(-a)^{n}%
q^{n(n\hspace{0.01in}-1)/2},\quad(n\in\mathbb{N}_{0}), \label{ZusPocInv}%
\end{equation}
it holds:%
\begin{equation}
\operatorname{e}_{q}((1-q)z)=\exp_{q}(z),\qquad\operatorname{E}_{q}%
((1-q)z)=\exp_{q^{-1}}(z). \label{ZusExpQ1}%
\end{equation}
Comparing the two identities in Eq.~(\ref{ZusExpQ1}), we find:%
\begin{equation}
\operatorname{e}_{q}((1-q)z)=\operatorname{E}_{q^{-1}}((1-q^{-1})z).
\end{equation}
Due to Eq.~(\ref{ZusPocInv}) and the two power series expansions in
Eq.~(\ref{qExpLit}), we have:%
\begin{equation}
\operatorname{e}_{q^{-1}}(z)=\operatorname{E}_{q}(-qz).
\end{equation}
The representations of $\operatorname{E}_{q}(z)$ and $\operatorname{e}_{q}(z)$
as products imply the following identity [cf. Eq.~(\ref{qExpLit})]:%
\begin{equation}
\operatorname{e}_{q}(z)\operatorname{E}_{q}(-z)=1.
\end{equation}

As we can see from its representation in Eq.~(\ref{qExpLit}), the
$q$-ex\-po\-nen\-tial function $\operatorname{e}_{q}(z)$ is a meromorphic
function with poles at $z=q^{-k}$, $k\in\mathbb{N}_{0}$. Hence, there is a
partial fraction decomposition for $\operatorname{e}_{q}(z)$
\cite{olshanetsky1998q.alg}, i.~e.%
\begin{equation}
\operatorname{e}_{q}(z)=\sum_{k\hspace{0.01in}=\hspace{0.01in}0}^{\infty}%
\frac{c_{k,n}}{1-zq^{k}} \label{ParExpqDef}%
\end{equation}
with%
\begin{equation}
c_{k,n}=\lim_{z\rightarrow\hspace{0.01in}q^{-k}}(1-zq^{k})\operatorname{e}%
_{q}(z)=\frac{(-1)^{k}\,q^{k(k\hspace{0.01in}+1)/2}}{(q\hspace{0.02in}%
;q)_{\infty}\,(q\hspace{0.02in};q)_{k}}.
\end{equation}
The series in Eq.~(\ref{ParExpqDef}) converges for $0<q<1$ and $z\neq q^{-2k}$.

We use the $q$-ex\-po\-nen\-tial in Eq.~(\ref{AusQExpEin}) to calculate
$q$\textbf{-trans\-la\-tions}
\cite{Chryssomalakos:1993zm,Majid:1993ud,rogov2000q.alg}. To this end, we
introduce the operator:%
\begin{equation}
\exp_{q}(a|D_{q})=\sum_{k\hspace{0.01in}=\hspace{0.01in}0}^{\infty}\frac
{a^{k}}{[[k]]_{q}!}\,D_{q}^{k}.\label{qTraOpeDef}%
\end{equation}
Applying this operator to a function $f(x)$, we get the following
$q$\textbf{-ver\-sion of Taylor's formula}:%
\begin{equation}
f(a\,\bar{\oplus}\,x)=\exp_{q}(a|D_{q})\triangleright f(x)=\sum_{k\hspace
{0.01in}=\hspace{0.01in}0}^{\infty}\frac{a^{k}}{[[k]]_{q}!}\,D_{q}^{k}%
\hspace{0.01in}f(x).\label{qTraVerGerUnKon}%
\end{equation}
Similiarly, it holds:%
\begin{equation}
f(x\,\bar{\oplus}\,a)=\sum_{k\hspace{0.01in}=\hspace{0.01in}0}^{\infty}%
\frac{1}{[[k]]_{q}!}\hspace{0.01in}(D_{q}^{k}\hspace{0.01in}f(x))\,a^{k}%
.\label{qTraVerGerKonUnk2}%
\end{equation}
Applying the operator in Eq.~(\ref{qTraOpeDef}) to a power function with a
natural number as the exponent, we obtain a $q$\textbf{-ver\-sion of the
binomial theorem}:%
\begin{align}
(a\,\bar{\oplus}\,x)^{n} &  =\exp_{q}(a|D_{q})\triangleright x^{n}%
=\sum_{k\hspace{0.01in}=\hspace{0.01in}0}^{\infty}\frac{a^{k}}{[[k]]_{q}%
!}\,D_{q}^{k}\triangleright x^{n}\nonumber\\
&  =\sum_{k\hspace{0.01in}=\hspace{0.01in}0}^{\infty}\frac{a^{k}}{[[k]]_{q}%
!}\,[[\hspace{0.01in}n]]_{q}\cdots\lbrack\lbrack\hspace{0.01in}n-k+1]]_{q}%
\,x^{n-k}=\sum_{k\hspace{0.01in}=\hspace{0.01in}0}^{n}%
\genfrac{[}{]}{0pt}{}{n}{k}%
_{q}a^{n-k}x^{k}.\label{qBinForGan}%
\end{align}
The $q$-binomial theorem implies the identities%
\begin{align}
(a\,\bar{\oplus}\,0)^{n} &  =\left.  (a\oplus x)\right\vert _{x\hspace
{0.01in}=\hspace{0.01in}0}=a^{n}\nonumber\\
(0\,\bar{\oplus}\,a)^{n} &  =\left.  (x\,\bar{\oplus}\,a)\right\vert
_{x\hspace{0.01in}=\hspace{0.01in}0}=a^{n},
\end{align}
which can be generalized as follows:%
\begin{align}
f(a\,\bar{\oplus}\,0) &  =\left.  f(a\,\bar{\oplus}\,x)\right\vert
_{x\hspace{0.01in}=\hspace{0.01in}0}=f(a),\nonumber\\
f(0\,\bar{\oplus}\,a) &  =\left.  f(x\,\bar{\oplus}\,y)\right\vert
_{x\hspace{0.01in}=\hspace{0.01in}0}=f(a).\label{NeuEleQTra}%
\end{align}

There is a second operator for calculating $q$-trans\-la\-tions:%
\begin{align}
\exp_{q^{-1}}(-a|D_{q})  & =\sum_{k\hspace{0.01in}=\hspace{0.01in}0}^{\infty
}\frac{(-a)^{k}}{[[k]]_{q^{-1}}!}\,D_{q}^{k}\nonumber\\
& =\sum_{k\hspace{0.01in}=\hspace{0.01in}0}^{\infty}\frac{q^{k(k\hspace
{0.01in}-1)/2}}{[[k]]_{q}!}\,(-a)^{k}D_{q}^{k}.\label{AusQExpEinInv}%
\end{align}
Applying this operator to a function $f(x)$, we get:%
\begin{align}
f((\bar{\ominus}\,a)\,\bar{\oplus}\,x)  & =\exp_{q^{-1}}(-a|D_{q}%
)\triangleright f(x)\nonumber\\
& =\sum_{k\hspace{0.01in}=\hspace{0.01in}0}^{\infty}\frac{q^{k(k\hspace
{0.01in}-1)/2}}{[[k]]_{q}!}\,(-a)^{k}D_{q}^{k}\hspace{0.01in}%
f(x).\label{QTayForTyp2}%
\end{align}
Accordingly, we have:%
\begin{equation}
f(x\,\bar{\oplus}\,(\bar{\ominus}\,a))=\sum_{k\hspace{0.01in}=\hspace
{0.01in}0}^{\infty}\frac{q^{k(k\hspace{0.01in}-1)/2}}{[[k]]_{q}!}\,(D_{q}%
^{k}\hspace{0.01in}f(x))(-a)^{k}.\label{QTayForTyp3}%
\end{equation}

The operators in Eqs.~(\ref{qTraOpeDef}) and (\ref{AusQExpEinInv}) are
inverses of each other. We show this by the following calculation:%
\begin{align}
&  \exp_{q^{-1}}(-a|D_{q})\exp_{q}(a|D_{q})\triangleright f(x)=\exp
_{q}(a|D_{q})\exp_{q^{-1}}(-a|D_{q})\triangleright f(x)\nonumber\\
&  \qquad=\sum_{k,j\hspace{0.01in}=\hspace{0.01in}0}^{\infty}\frac
{(-1)^{k}\hspace{0.01in}q^{k(k\hspace{0.01in}-1)/2}}{[[k]]_{q}!\,[[\hspace
{0.01in}j]]_{q}!}\,a^{k\hspace{0.01in}+j}\hspace{0.01in}D_{q}^{k\hspace
{0.01in}+j}\hspace{0.01in}f(x)\nonumber\\
&  \qquad=\sum_{k,j\hspace{0.01in}=\hspace{0.01in}0}^{\infty}\frac
{(-1)^{k}\hspace{0.01in}q^{k(k\hspace{0.01in}-1)/2}}{[[k+j]]_{q}!}%
\frac{[[\hspace{0.01in}k+j]]_{q}!}{[[k]]_{q}!\,[[\hspace{0.01in}j]]_{q}%
!}\,a^{k\hspace{0.01in}+j}\hspace{0.01in}D_{q}^{k\hspace{0.01in}+j}%
\hspace{0.01in}f(x)\nonumber\\
&  \qquad=\sum_{n\hspace{0.01in}=\hspace{0.01in}0}^{\infty}\frac{1}%
{[[\hspace{0.01in}n]]_{q}!}\sum_{k\hspace{0.01in}=\hspace{0.01in}0}%
^{n}(-1)^{k}\hspace{0.01in}q^{k(k\hspace{0.01in}-1)/2}\frac{[[\hspace
{0.01in}n]]_{q}!}{[[k]]_{q}!\,[[\hspace{0.01in}n-k]]_{q}!}\,a^{n}%
\hspace{0.01in}D_{q}^{n}\hspace{0.01in}f(x)\nonumber\\
&  \qquad=\sum_{n\hspace{0.01in}=\hspace{0.01in}0}^{\infty}\frac{1}%
{[[\hspace{0.01in}n]]_{q}!}\sum_{k\hspace{0.01in}=\hspace{0.01in}0}%
^{n}(-1)^{k}\hspace{0.01in}q^{k(k\hspace{0.01in}-1)/2}%
\genfrac{[}{]}{0pt}{}{n}{k}%
_{q}\,a^{n}\hspace{0.01in}D_{q}^{n}\hspace{0.01in}f(x)=f(x). \label{NacInvTra}%
\end{align}
In the last step of the above calculation, we used the following identity
\cite{Klimyk:1997eb}:%
\begin{equation}
\sum_{k\hspace{0.01in}=\hspace{0.01in}0}^{n}\hspace{0.01in}(-1)^{k}%
\hspace{0.01in}q^{k(k\hspace{0.01in}-1)/2}%
\genfrac{[}{]}{0pt}{}{n}{k}%
_{q}=(1\hspace{0.01in};q)_{n}=\left\{
\begin{tabular}
[c]{cl}%
$1$ & if\quad$n=0,$\\
$0$ & if\quad$n\geq1.$%
\end{tabular}
\ \ \right.  \label{SumBinNul}%
\end{equation}

If we shift the $q$-ex\-po\-nen\-tial by a $q$-trans\-la\-tion, the following
\textbf{addition theorem} holds \cite{Chryssomalakos:1993zm}:%
\begin{equation}
\exp_{q}(x\,\bar{\oplus}\,a)=\exp_{q}(x)\exp_{q}(a).\label{AddTheqExp1Dim}%
\end{equation}
We prove this rule as follows [cf. Eq.~(\ref{qBinForGan})]:%
\begin{align}
&  \exp_{q}(x)\exp_{q}(a)=\sum_{k\hspace{0.01in}=\hspace{0.01in}0}^{\infty
}\sum_{l\hspace{0.01in}=\hspace{0.01in}0}^{\infty}\frac{1}{[[k]]_{q}!}\frac
{1}{[[\hspace{0.01in}l]]_{q}!}\,x^{k}a^{l}\nonumber\\
&  \qquad=\sum_{k\hspace{0.01in}=\hspace{0.01in}0}^{\infty}\sum_{l\hspace
{0.01in}=\hspace{0.01in}0}^{\infty}\frac{1}{[[k+l]]_{q}!}\frac{[[k+l]]_{q}%
!}{[[k]]_{q}!\hspace{0.01in}[[\hspace{0.01in}l]]_{q}!}\,x^{k}a^{l}\nonumber\\
&  \qquad=\sum_{n\hspace{0.01in}=\hspace{0.01in}0}^{\infty}\frac{1}%
{[[\hspace{0.01in}n]]_{q}!}\sum_{l\hspace{0.01in}=\hspace{0.01in}0}^{n}%
\genfrac{[}{]}{0pt}{}{n}{l}%
_{q}x^{n-l}a^{l}\nonumber\\
&  \qquad=\sum_{n\hspace{0.01in}=\hspace{0.01in}0}^{\infty}\frac{1}%
{[[\hspace{0.01in}n]]_{q}!}\,(x\,\bar{\oplus}\,a)^{n}=\exp_{q}(x\,\bar{\oplus
}\,a).
\end{align}

Using Eq.~(\ref{QTayForTyp2}), we can calculate so-called $q$%
\textbf{-in\-ver\-sions}:\footnote{In Chap.$~$\ref{KapEinDimQFouTraN}, we use
the following notation for the $q$-in\-ver\-sion: $\underline{\bar{S}%
}_{\hspace{0.01in}x\hspace{0.01in}}\hspace{-0.03in}(f(x))=f(\bar{\ominus
}\,x).$}%
\begin{equation}
f((\bar{\ominus}\,x))=\sum_{k\hspace{0.01in}=\hspace{0.01in}0}^{\infty}%
\frac{q^{k(k\hspace{0.01in}-1)/2}}{[[k]]_{q}!}\,(-x)^{k}(D_{q}^{k}%
\hspace{0.01in}f)(0).
\end{equation}
For powers of $x$, we find:%
\begin{equation}
(\bar{\ominus}\,x)^{n}=q^{n(n\hspace{0.01in}-1)/2}\,(-x)^{n}. \label{KonAntXn}%
\end{equation}
Taking Eqs.~(\ref{Bed1qExp}) and (\ref{Bed2qExp}) into account, we also get
[cf. Eq.~(\ref{AusQExpEin})]:%
\begin{align}
\exp_{q}(\bar{\ominus}\,x)  &  =\sum_{k\hspace{0.01in}=\hspace{0.01in}%
0}^{\infty}\frac{q^{k(k\hspace{0.01in}-1)/2}}{[[k]]_{q}!}\,(-x)^{k}\exp
_{q}(0)\nonumber\\
&  =\sum_{k\hspace{0.01in}=\hspace{0.01in}0}^{\infty}\frac{1}{[[k]]_{q^{-1}}%
!}\,(-x)^{k}=\exp_{q^{-1}}(-x). \label{ExpForInvQExp}%
\end{align}
By similar reasoning as in Eq.~(\ref{NacInvTra}), we can prove the following
identities:%
\[
\exp_{q}(\bar{\ominus}\,x)\exp_{q}(x)=\exp_{q}(x)\exp_{q}(\bar{\ominus
}\,x)=1.
\]

\section{Trigonometric functions\label{KapQTrigFkt}}

We use the $q$-ex\-po\-nen\-tial function $\operatorname{e}_{q}(z)$ to
introduce $q$\textbf{-ver\-sions of trigonometric functions}
\cite{Klimyk:1997eb}, i.~e.%
\begin{align}
\widetilde{\sin}_{q}(z) &  =\frac{1}{2\text{i}}\left(  \operatorname{e}%
_{q}(\text{i}z)-\operatorname{e}_{q}(-\text{i}z)\right)  =\sum_{n\hspace
{0.01in}=\hspace{0.01in}0}^{\infty}\frac{(-1)^{n}\hspace{0.01in}%
z^{2n\hspace{0.01in}+1}}{(q\hspace{0.02in};q)_{2n\hspace{0.01in}+1}%
},\nonumber\\
\widetilde{\cos}_{q}(z) &  =\frac{1}{2}\left(  \operatorname{e}_{q}%
(\text{i}z)+\operatorname{e}_{q}(-\text{i}z)\right)  =\sum_{n\hspace
{0.01in}=\hspace{0.01in}0}^{\infty}\frac{(-1)^{n}\hspace{0.01in}z^{2n}%
}{(q\hspace{0.02in};q)_{2n}},\label{TilDefqDefTriFkt}%
\end{align}
and%
\begin{align}
\widetilde{\operatorname*{Sin}}_{q}(z) &  =\frac{1}{2\text{i}}\left(
\operatorname{E}_{q}(\text{i}z)-\operatorname{E}_{q}(-\text{i}z)\right)
,\nonumber\\
\widetilde{\operatorname*{Cos}}_{q}(z) &  =\frac{1}{2}\left(  \operatorname{E}%
_{q}(\text{i}z)+\operatorname{E}_{q}(-\text{i}z)\right)
.\label{TilDefqDefTriFkt2}%
\end{align}

Using the partial fraction decomposition in Eq.~(\ref{ParExpqDef}) of the
previous chapter, we get the following expressions for the $q$%
-trigono\-met\-ric functions in Eq.~(\ref{TilDefqDefTriFkt}):%
\begin{align}
\widetilde{\sin}_{q}(z)  &  =\frac{z}{(q\hspace{0.02in};q)_{\infty}}%
\sum_{k\hspace{0.01in}=\hspace{0.01in}0}^{\infty}\frac{(-1)^{k}\hspace
{0.02in}q^{k(k+3)/2}}{(q\hspace{0.02in};q)_{k}\hspace{0.02in}(1+z^{2}q^{2k}%
)},\nonumber\\
\widetilde{\cos}_{q}(z)  &  =\frac{1}{(q\hspace{0.02in};q)_{\infty}}%
\sum_{k\hspace{0.01in}=\hspace{0.01in}0}^{\infty}\frac{(-1)^{k}\hspace
{0.02in}q^{k(k\hspace{0.01in}+1)/2}}{(q\hspace{0.02in};q)_{k}\hspace
{0.02in}(1+z^{2}q^{2k})}. \label{ParZerTriFktqDef}%
\end{align}
For $x\in\mathbb{R}$, we can estimate both series in
Eq.~(\ref{ParZerTriFktqDef}) by the following inequalities:%
\begin{align}
\left\vert \hspace{0.01in}\sum_{k\hspace{0.01in}=\hspace{0.01in}0}^{\infty
}\frac{(-1)^{k}\hspace{0.02in}q^{k(k+3)/2}}{(q\hspace{0.02in};q)_{k}%
\hspace{0.02in}(1+x^{2}q^{2k})}\right\vert  &  \leq\sum_{k\hspace
{0.01in}=\hspace{0.01in}0}^{\infty}\frac{q^{k(k\hspace{0.01in}-1)/2}%
}{(q\hspace{0.02in};q)_{k}\hspace{0.02in}(1+x^{2})}=\frac{\operatorname{E}%
_{q}(1)}{1+x^{2}},\nonumber\\
\left\vert \hspace{0.01in}\sum_{k\hspace{0.01in}=\hspace{0.01in}0}^{\infty
}\frac{(-1)^{k}\hspace{0.02in}q^{k(k\hspace{0.01in}+1)/2}}{(q\hspace
{0.02in};q)_{k}\hspace{0.02in}(1+x^{2}q^{2k})}\right\vert  &  \leq
\sum_{k\hspace{0.01in}=\hspace{0.01in}0}^{\infty}\frac{q^{k(k\hspace
{0.01in}-1)/2}q^{-k}}{(q\hspace{0.02in};q)_{k}\hspace{0.02in}(1+x^{2})}%
=\frac{\operatorname{E}_{q}(q^{-1})}{1+x^{2}}. \label{AbsSumParZerTriFkt}%
\end{align}
For $0<q<1$ and $k\in\mathbb{N}_{0}$, the inequalities above result from%
\begin{equation}
\frac{1}{1+x^{2}q^{2k}}\leq\frac{q^{-2k}}{1+x^{2}}.
\end{equation}
Using the results in Eq.~(\ref{AbsSumParZerTriFkt}), we finally get [also cf.
Eq.~(\ref{qExpLit}) of the previous chapter]:%
\begin{equation}
\big|\hspace{0.01in}\widetilde{\sin}_{q}(x)\big|\leq\frac{(-1\hspace
{0.02in};q)_{\infty}\left\vert \hspace{0.01in}x\right\vert }{(q\hspace
{0.02in};q)_{\infty}\hspace{0.02in}(1+x^{2})},\qquad\big|\hspace
{0.01in}\widetilde{\cos}_{q}(x)\big|\leq\frac{(-q\hspace{0.02in};q)_{\infty}%
}{(q\hspace{0.02in};q)_{\infty}\hspace{0.02in}(1+x^{2})}.
\label{AbsqTriSinCos}%
\end{equation}
There are also inequalities for the $q$-trigono\-met\-ric functions in
Eq.~(\ref{TilDefqDefTriFkt2}) (cf. Ref.~\cite{olshanetsky1998q.alg}):%
\begin{equation}
\Big|\hspace{0.01in}\widetilde{\operatorname*{Sin}}_{q}(x)\Big|\leq\left\vert
\hspace{0.01in}x\right\vert ,\qquad\Big|\hspace{0.01in}\widetilde
{\operatorname*{Cos}}_{q}(x)\Big|\leq1. \label{AbsqTriSinCos2}%
\end{equation}

The $q$-trigono\-met\-ric functions in Eqs.~(\ref{TilDefqDefTriFkt}) and
(\ref{TilDefqDefTriFkt2}) do not become the ordinary trigonometric functions
if we take $q\rightarrow1$. For this reason, we introduce the new
$q$-trigono\-me\-tric functions [cf. Eq.~(\ref{ZusExpQ1}) of the previous
chapter]%
\begin{align}
\sin_{q}(x)  &  =\widetilde{\sin}_{q}((1-q)\hspace{0.01in}x)=\frac
{1}{2\text{i}}\left(  \exp_{q}(\text{i}\hspace{0.01in}x)-\exp_{q}%
(-\text{i}\hspace{0.01in}x)\right)  ,\nonumber\\
\cos_{q}(x)  &  =\widetilde{\cos}_{q}((1-q)\hspace{0.01in}x)=\frac{1}%
{2}\left(  \exp_{q}(\text{i}\hspace{0.01in}x)+\exp_{q}(-\text{i}%
\hspace{0.01in}x)\right)  , \label{DefqTriSinCosExp1}%
\end{align}
and%
\begin{align}
\operatorname*{Sin}\nolimits_{q}(x)  &  =\widetilde{\operatorname*{Sin}}%
_{q}((1-q)\hspace{0.01in}x)=\frac{1}{2\text{i}}\left(  \exp_{q^{-1}}%
(\text{i}\hspace{0.01in}x)-\exp_{q^{-1}}(-\text{i}\hspace{0.01in}x)\right)
,\nonumber\\
\operatorname*{Cos}\nolimits_{q}(x)  &  =\widetilde{\operatorname*{Cos}}%
_{q}((1-q)\hspace{0.01in}x)=\frac{1}{2}\left(  \exp_{q^{-1}}(\text{i}%
\hspace{0.01in}x)+\exp_{q^{-1}}(-\text{i}\hspace{0.01in}x)\right)  .
\label{DefqTriSinCosExp2}%
\end{align}
Due to Eq.~(\ref{KlaGreQExp}) of the previous chapter, we regain the ordinary
trigonometric functions if $q\rightarrow1$:%
\begin{equation}
\lim_{q\rightarrow1^{-}}\sin_{q}(x)=\sin(x),\qquad\lim_{q\rightarrow1^{-}}%
\cos_{q}(x)=\cos(x).
\end{equation}
From the expressions in Eqs.~(\ref{DefqTriSinCosExp1}) and
(\ref{DefqTriSinCosExp2}), we get the following identities:%
\begin{equation}
\operatorname*{Sin}\nolimits_{q}(x)=\sin_{q^{-1}}(x),\qquad\operatorname*{Cos}%
\nolimits_{q}(x)=\cos_{q^{-1}}(x).
\end{equation}
Furthermore, there is a $q$-ver\-sion of the identity $\sin^{2}(x)+\cos
^{2}(x)=1$ (also cf. Ref.~\cite{Kac:2002eb}):%
\begin{equation}
\sin_{q}(x)\sin_{q^{-1}}(x)+\cos_{q}(x)\cos_{q^{-1}}(x)=1.
\end{equation}

By applying the Jackson derivatives to the $q$-trigono\-met\-ric functions in
Eqs.~(\ref{DefqTriSinCosExp1}) and (\ref{DefqTriSinCosExp2}), we find (also
cf. Ref.~\cite{Kac:2002eb}):%
\begin{align}
D_{q}\sin_{q}(x)  &  =\cos_{q}(x), & D_{q}\operatorname*{Sin}\nolimits_{q}(x)
&  =\operatorname*{Cos}\nolimits_{q}(qx),\nonumber\\
D_{q}\cos_{q}(x)  &  =-\sin_{q}(x), & D_{q}\operatorname*{Cos}\nolimits_{q}%
(x)  &  =-\operatorname*{Sin}\nolimits_{q}(qx).
\end{align}
Accordingly, $\sin_{q}(\omega x)$ and $\cos_{q}(\omega x)$ are linearly
independent solutions of the difference equation (also cf.
Ref.~\cite{Klimyk:1997eb})%
\begin{equation}
D_{q}^{2}f(x)+\omega^{2}f(x)=0,
\end{equation}
and the $q$-trigono\-met\-ric functions $\operatorname*{Sin}\nolimits_{q}%
(\omega x)$ and $\operatorname*{Cos}\nolimits_{q}(\omega x)$ satisfy the
difference equation%
\begin{equation}
D_{q}^{2}f(x)+\omega^{2}f(q^{2}x)=0.
\end{equation}

\section{Integrals\label{KapQIntTrig}}

For $z>0$ and $0<q<1$, the one-di\-men\-sion\-al \textbf{Jackson integral} is
defined as follows \cite{Jackson:1908}:%
\begin{align}
\int_{z}^{\hspace{0.01in}z.\infty}\text{d}_{q}x\hspace{0.01in}f(x)  &
=(1-q)\hspace{0.01in}z\sum_{j\hspace{0.01in}=1}^{\infty}q^{-j}f(q^{-j}%
z),\nonumber\\
\int_{0}^{\hspace{0.01in}z}\text{d}_{q}x\hspace{0.01in}f(x)  &  =(1-q)\hspace
{0.01in}z\sum_{j\hspace{0.01in}=\hspace{0.01in}0}^{\infty}q^{\hspace{0.01in}%
j}f(q^{\hspace{0.01in}j}z). \label{QInt0klQkl1}%
\end{align}
Accordingly, we have:%
\begin{align}
\int_{0}^{\hspace{0.01in}z.\infty}\text{d}_{q}x\hspace{0.01in}f(x)  &
=\int_{0}^{\hspace{0.01in}z}\text{d}_{q}x\hspace{0.01in}f(x)+\int_{z}%
^{\hspace{0.01in}z.\infty}\text{d}_{q}x\hspace{0.01in}f(x)\nonumber\\
&  =(1-q)\hspace{0.01in}z\sum_{j\hspace{0.01in}=-\infty}^{\infty}%
q^{-j}f(q^{-j}z). \label{ImpJacInt}%
\end{align}

So far, we have only considered $q$-in\-te\-grals with a non-neg\-a\-tive
integration range. We can define $q$-in\-te\-grals with a negative integration
range as well ($z<0$ and $0<q<1$):%
\begin{align}
\int_{-\infty}^{\hspace{0.01in}z}\text{d}_{q}x\hspace{0.01in}f(x)  &
=(q-1)\,z\sum_{j\hspace{0.01in}=1}^{\infty}q^{-j}f(q^{-j}z),\nonumber\\
\int_{z}^{\hspace{0.01in}0}\text{d}_{q}x\hspace{0.01in}f(x)  &  =(q-1)\,z\sum
_{j\hspace{0.01in}=\hspace{0.01in}0}^{\infty}q^{\hspace{0.01in}j}%
f(q^{\hspace{0.01in}j}z). \label{QInt0klQkl1Neg}%
\end{align}
Accordingly, we have:%
\begin{equation}
\int_{z.\infty}^{\hspace{0.01in}0}\text{d}_{q}x\hspace{0.01in}%
f(x)=(q-1)\,z\sum_{j\hspace{0.01in}=\hspace{0.01in}-\infty}^{\infty}%
q^{\hspace{0.01in}j}f(q^{\hspace{0.01in}j}z). \label{JacIntMinInfZerNeg}%
\end{equation}

The $q$-in\-te\-grals in Eqs.~(\ref{ImpJacInt}) and (\ref{JacIntMinInfZerNeg})
can be combined to form a single $q$-in\-te\-gral over the interval $\left]
-\infty,\infty\right[  $:%
\begin{align}
\int_{-z.\infty}^{\hspace{0.01in}z.\infty}\text{d}_{q}x\hspace{0.01in}f(x)  &
=\int_{0}^{\hspace{0.01in}z.\infty}\text{d}_{q}x\hspace{0.01in}f(x)+\int
_{-\infty.z}^{\hspace{0.01in}0}\text{d}_{q}x\hspace{0.01in}f(x)\nonumber\\
&  =(1-q)\hspace{0.01in}z\sum_{j\hspace{0.01in}=\hspace{0.01in}-\infty
}^{\infty}q^{\hspace{0.01in}j}\left[  f(q^{\hspace{0.01in}j}z)+f(-q^{\hspace
{0.01in}j}z)\right]  . \label{UneJackIntAll}%
\end{align}
Eq.~(\ref{UneJackIntAll}) implies:%
\begin{equation}
\int_{-x.\infty}^{\hspace{0.01in}x.\infty}\text{d}_{q}z\hspace{0.01in}%
f(z)=\int_{-x.\infty}^{\hspace{0.01in}x.\infty}\text{d}_{q}z\hspace
{0.01in}f(-z). \label{JacInvYSpi}%
\end{equation}

A look at Eq.~(\ref{ImpJacInt}) or Eq.~(\ref{JacIntMinInfZerNeg}) shows that
only the lattice points $\left\{  q^{k}x\,|\,k\in\mathbb{Z}\right\}  $
contribute to the $q$-in\-te\-gral over $\left]  0,\infty\right[  $ or
$\left]  -\infty,0\right[  $. Thus, the integrand $f$ has to satisfy the
following boundary condition for the $q$-in\-te\-gral over $\left]
0,\infty\right[  $ or $\left]  -\infty,0\right[  $ to converge ($q<1$):%
\begin{equation}
\lim_{k\hspace{0.01in}\rightarrow-\infty}f(\pm\hspace{0.01in}q^{k}x)=0.
\label{RanBedFunEin}%
\end{equation}
If $f$ is continuous at $0$, the boundary condition above implies:%
\begin{align}
\int\nolimits_{0}^{\hspace{0.01in}z.\infty}\text{d}_{q}x\hspace{0.02in}%
D_{q}f(x)  &  =(1-q)\sum_{j\hspace{0.01in}=\hspace{0.01in}-\infty}^{\infty
}zq^{\hspace{0.01in}-j}\,\frac{f(q^{\hspace{0.01in}-j}z)-f(q^{\hspace
{0.01in}-j+1}z)}{(1-q)\hspace{0.01in}q^{\hspace{0.01in}-j}z}\nonumber\\
&  =\lim_{k\hspace{0.01in}\rightarrow\infty}\sum_{j=-k}^{k}\left[
f(q^{\hspace{0.01in}-j}z)-f(q^{\hspace{0.01in}-j+1}z)\right] \nonumber\\
&  =\lim_{k\hspace{0.01in}\rightarrow\infty}\left[  f(q^{-k}z)-f(q^{k+1}%
z)\right]  =-f(0). \label{PosParIntAblq1Dim}%
\end{align}
A similar reasoning shows:%
\begin{equation}
\int\nolimits_{-\infty.z}^{\hspace{0.01in}0}\text{d}_{q}x\hspace{0.02in}%
D_{q}f(x)=f(0). \label{NegParIntAblq1Dim}%
\end{equation}
We can combine Eqs.~(\ref{PosParIntAblq1Dim}) and (\ref{NegParIntAblq1Dim}) to
give the following result ($m\in\mathbb{N}$):%
\begin{align}
\int\nolimits_{-z.\infty}^{\hspace{0.01in}z.\infty}\text{d}_{q}x\hspace
{0.02in}D_{q}^{m}f(x)  &  =\int\nolimits_{0}^{\hspace{0.01in}z.\infty}%
\text{d}_{q}x\hspace{0.02in}D_{q}^{m}f(x)+\int\nolimits_{-\infty.z}%
^{\hspace{0.01in}0}\text{d}_{q}x\hspace{0.02in}D_{q}^{m}f(x)\nonumber\\
&  =-D_{q}^{m-1}f(0)+D_{q}^{m-1}f(0)=0. \label{UnIntJacAbl}%
\end{align}

If the two functions $f$ and $g$ vanish at infinity,\footnote{The weaker
condition $\lim_{k\rightarrow-\infty}f(\pm q^{k}a)\hspace{0.01in}g(\pm
q^{k}a)=0$ is also sufficient.} Eq.~(\ref{UnIntJacAbl}) and
Eq.~(\ref{ProRegJacAblBas}) of\ Chap.~\ref{KapqAnaExpTriFkt} imply the
following $q$\textbf{-de\-formed rule for integration by parts}:%
\begin{equation}
\int\nolimits_{-z.\infty}^{\hspace{0.01in}z.\infty}\text{d}_{q}x\hspace
{0.01in}f(x)\hspace{0.01in}D_{q}\hspace{0.01in}g(x)=-\int\nolimits_{-z.\infty
}^{\hspace{0.01in}z.\infty}\text{d}_{q}x\hspace{0.01in}\left[  D_{q}%
f(x)\right]  g(qx).
\end{equation}
Applying this rule several times, we get the following identity
\cite{olshanetsky1998q.alg} ($k\in\mathbb{N}$):%
\begin{equation}
\int\nolimits_{-z.\infty}^{\hspace{0.01in}z.\infty}\text{d}_{q}x\hspace
{0.01in}f(x)\hspace{0.01in}D_{q}^{k}g(x)=(-1)^{k}q^{-k(k\hspace{0.01in}%
-1)/2}\int\nolimits_{-z.\infty}^{\hspace{0.01in}z.\infty}\text{d}_{q}%
x\hspace{0.01in}\left[  D_{q}^{k}f(x)\right]  g(q^{k}%
x).\label{RegParIntBraLinK}%
\end{equation}

Using Eq.~(\ref{UnIntJacAbl}), we can show that the $q$\textit{-in\-te\-gral
over }$\left]  -\infty,\infty\right[  $ \textit{is invariant under
translation} \cite{Chryssomalakos:1993zm,Kempf:1994yd,1996q.alg.....8008K}:%
\begin{align}
&  \int\nolimits_{-z.\infty}^{\hspace{0.01in}z.\infty}\text{d}_{q}%
x\hspace{0.02in}f(x\,\bar{\oplus}\,a)=\int\nolimits_{-z.\infty}^{\hspace
{0.01in}z.\infty}\text{d}_{q}x\sum_{k\hspace{0.01in}=\hspace{0.01in}0}%
^{\infty}\frac{1}{[[k]]_{q}!}\left[  D_{q}^{k}f(x)\right]  a^{k}\nonumber\\
&  \qquad=\sum_{k\hspace{0.01in}=\hspace{0.01in}0}^{\infty}\frac{1}%
{[[k]]_{q}!}\int\nolimits_{-z.\infty}^{\hspace{0.01in}z.\infty}\text{d}%
_{q}x\left[  D_{q}^{k}f(x)\right]  a^{k}=\int\nolimits_{-z.\infty}%
^{\hspace{0.01in}z.\infty}\text{d}_{q}x\hspace{0.02in}f(x).
\label{qTraInvEinJacIntUnk}%
\end{align}
In the first step, we used Eq.~(\ref{qTraVerGerKonUnk2}) of
Chap.~\ref{KapqAnaExpTriFkt}. In the last step, the summands to $k>0$ vanish
due to Eq.~(\ref{UnIntJacAbl}). Similiarly, it holds:%
\begin{equation}
\int\nolimits_{-z.\infty}^{\hspace{0.01in}z.\infty}\text{d}_{q}x\hspace
{0.02in}f(a\,\bar{\oplus}\,x)=\int\nolimits_{-z.\infty}^{\hspace
{0.01in}z.\infty}\text{d}_{q}x\hspace{0.02in}f(x).
\end{equation}

Using the Jackson integral in Eq.~(\ref{ImpJacInt}), we can introduce the
function $\Theta_{q}(z)$ ($0<q<1$):%
\begin{equation}
\Theta_{q}(z)=\int\nolimits_{0}^{\hspace{0.01in}z.\infty}\text{d}_{q}%
x\,x^{-1}\sin_{q}(x)=(1-q)\sum_{m\hspace{0.01in}=\hspace{0.01in}-\infty
}^{\infty}\sin_{q}(q^{m}z). \label{DefTheZ}%
\end{equation}
Due to the definition of the improper $q$-in\-te\-gral in Eq.~(\ref{ImpJacInt}%
), it holds $(k\in\mathbb{Z)}$:%
\begin{equation}
\Theta_{q}(q^{k}z)=\Theta_{q}(z). \label{PerTheFkt}%
\end{equation}
Next, we show that $\Theta_{q}(z)$ is related to the function%
\begin{equation}
\mathbf{Q}(z\hspace{0.01in};q)=(1-q)\sum_{m\hspace{0.01in}=\hspace
{0.01in}-\infty}^{\infty}\frac{1}{z\hspace{0.01in}q^{m}+z^{-1}q^{-m}}.
\label{DefQzqFkt}%
\end{equation}
To this end, we rewrite the series in Eq.~(\ref{DefTheZ}) using the partial
fraction decompositions from Eq.~(\ref{ParZerTriFktqDef}) of the previous
chapter \cite{olshanetsky1998q.alg}:%
\begin{align}
\Theta_{q}(z)  &  =\frac{1-q}{(q\hspace{0.01in};q)_{\infty}}\sum
_{m\hspace{0.01in}=\hspace{0.01in}-\infty}^{\infty}q^{m}\sum_{k\hspace
{0.01in}=\hspace{0.01in}0}^{\infty}\frac{(-1)^{k}\,q^{k(k+3)/2}(1-q)z}%
{(q\hspace{0.01in};q)_{k}\,\left[  1+(1-q)^{2}z^{2}q^{2(m+k)}\right]
}\nonumber\\
&  =\frac{1-q}{(q\hspace{0.01in};q)_{\infty}}\sum_{k\hspace{0.01in}%
=\hspace{0.01in}0}^{\infty}\frac{(-1)^{k}\,q^{k(k+3)/2}}{(q\hspace
{0.01in};q)_{k}}\sum_{m\hspace{0.01in}=\hspace{0.01in}-\infty}^{\infty}%
\frac{(1-q)z\hspace{0.01in}q^{m}}{1+(1-q)^{2}z^{2}q^{2(m+k)}}\nonumber\\
&  =\frac{1-q}{(q\hspace{0.01in};q)_{\infty}}\sum_{k\hspace{0.01in}%
=\hspace{0.01in}0}^{\infty}\frac{(-1)^{k}\,q^{k(k\hspace{0.01in}+1)/2}%
}{(q\hspace{0.01in};q)_{k}}\sum_{m\hspace{0.01in}=\hspace{0.01in}-\infty
}^{\infty}\frac{(1-q)z\hspace{0.01in}q^{m}}{1+(1-q)^{2}z^{2}q^{2m}}\nonumber\\
&  =(1-q)\,\frac{\operatorname{E}_{q}(-q)}{(q\hspace{0.01in};q)_{\infty}}%
\sum_{m\hspace{0.01in}=\hspace{0.01in}-\infty}^{\infty}\frac{(1-q)z\hspace
{0.01in}q^{m}}{1+(1-q)^{2}z^{2}q^{2m}}=\mathbf{Q}((1-q)z;q).
\label{IdeTheQUmr}%
\end{align}
In the second step, we changed the order of summation, and in the next step,
we replaced $m+k$ by $m$. The two final steps are a consequence of
Eq.~(\ref{qExpLit}) in Chap.~\ref{KapqAnaExpTriFkt}\ and Eq.~(\ref{DefQzqFkt}).

Next, we show that the following limit holds $(z\in\mathbb{R}^{+})$:%
\begin{equation}
\lim_{q\hspace{0.01in}\rightarrow1^{-}}\mathbf{Q}(z\hspace{0.01in}%
;q)=\frac{\pi}{2}. \label{KlaGreQFkt}%
\end{equation}
In Eq.~(\ref{DefQzqFkt}), we replace the deformation parameter $q$ by
$\operatorname{e}^{-h}$ ($h>0$):%
\begin{align}
\mathbf{Q}(z\hspace{0.01in};\operatorname{e}^{-h})  &  =(1-\operatorname{e}%
^{-h})\sum_{m\hspace{0.01in}=\hspace{0.01in}-\infty}^{\infty}\frac
{1}{z\operatorname{e}^{-hm}+\hspace{0.01in}z^{-1}\operatorname{e}^{hm}%
}\nonumber\\
&  =h\sum_{m\hspace{0.01in}=\hspace{0.01in}-\infty}^{\infty}\frac
{z^{-1}\operatorname{e}^{hm}}{1+(z^{-1}\operatorname{e}^{hm})^{2}}+O(h^{2}).
\end{align}
We take the limit $h\rightarrow0^{+}$ that corresponds to the limit
$q\rightarrow1$:%
\begin{align}
&  \lim_{h\hspace{0.01in}\rightarrow\hspace{0.01in}0^{+}}\mathbf{Q}%
(z\hspace{0.01in};\operatorname{e}^{-h})=\lim_{h\hspace{0.01in}\rightarrow
\hspace{0.01in}0^{+}}\sum_{m\hspace{0.01in}=\hspace{0.01in}-\infty}^{\infty
}\frac{z^{-1}\operatorname{e}^{hm}}{1+(z^{-1}\operatorname{e}^{hm})^{2}%
}\hspace{0.01in}h=\int\nolimits_{-\infty}^{\hspace{0.01in}\infty}\frac
{z^{-1}\operatorname{e}^{x}}{1+(z^{-1}\operatorname{e}^{x})^{2}}%
\,\text{d}\hspace{0.01in}x\nonumber\\
&  \qquad=\int\nolimits_{-\infty}^{\hspace{0.01in}\infty}\frac{\text{d}%
(z^{-1}\operatorname{e}^{x})}{1+(z^{-1}\operatorname{e}^{x})^{2}}%
=\int\nolimits_{0}^{\hspace{0.01in}\infty}\frac{\text{d}\hspace{0.01in}%
y}{1+y^{2}}=\frac{1}{2}\int\nolimits_{-\infty}^{\hspace{0.01in}\infty}%
\frac{\text{d}\hspace{0.01in}y}{1+y^{2}}=\frac{\pi}{2}.
\end{align}
In the above calculation, the series becomes an improper integral for
$h\rightarrow0^{+}$. We solved this improper integral using the substitution
rule. Since the limit in Eq.~(\ref{KlaGreQFkt}) is independent of $z$, we also
have:%
\begin{equation}
\lim_{q\hspace{0.01in}\rightarrow1^{-}}\Theta_{q}(z)=\lim_{q\hspace
{0.01in}\rightarrow1^{-}}\mathbf{Q}((1-q)z\hspace{0.01in};q)=\frac{\pi}{2}.
\label{TheQFktKlaGre}%
\end{equation}
Accordingly, Eq.~(\ref{DefTheZ}) implies ($z\in\mathbb{R}^{+}$):%
\begin{equation}
\lim_{q\hspace{0.01in}\rightarrow1}\Theta_{q}(z)=\lim_{q\hspace{0.01in}%
\rightarrow1}\int\nolimits_{0}^{\hspace{0.01in}z.\infty}\text{d}_{q}%
\hspace{0.01in}x\,x^{-1}\sin_{q}(x)=\int\nolimits_{0}^{\hspace{0.01in}\infty
}\text{d}\hspace{0.01in}x\,x^{-1}\sin(x)=\frac{\pi}{2}.
\end{equation}

Using the power series expansions in Eq.~(\ref{TilDefqDefTriFkt}) of the
previous chapter, we can compute $q$\textbf{-in\-te\-grals of }$q$%
\textbf{-trigo\-nomet\-ric functions}. For example, we find
\cite{olshanetsky1998q.alg} ($0<q<1$, $M\in\mathbb{N}$):%
\begin{align}
&  \int\nolimits_{0}^{\hspace{0.01in}zq^{-M}}\text{d}_{q}\hspace
{0.01in}x\,\cos_{q}(x)=(1-q)z\sum_{m\hspace{0.01in}=\hspace{0.01in}-M}%
^{\infty}q^{m}\cos_{q}(q^{m}z)\nonumber\\
&  \qquad=\sum_{m\hspace{0.01in}=\hspace{0.01in}-M}^{\infty}(1-q)z\hspace
{0.01in}q^{m}\sum_{k\hspace{0.01in}=\hspace{0.01in}0}^{\infty}\frac
{(-1)^{k}\,(1-q)^{2k}q^{2mk}z^{2k}}{(q\hspace{0.01in};q)_{2k}}\nonumber\\
&  \qquad=\sum_{k\hspace{0.01in}=\hspace{0.01in}0}^{\infty}\frac
{(-1)^{k}\,(1-q)^{2k\hspace{0.01in}+1}z^{2k\hspace{0.01in}+1}}{(q\hspace
{0.01in};q)_{2k}}\sum_{m=-M}^{\infty}q^{m(2k\hspace{0.01in}+1)}\nonumber\\
&  \qquad=\sum_{k\hspace{0.01in}=\hspace{0.01in}0}^{\infty}\frac
{(-1)^{k}\,(1-q)^{2k\hspace{0.01in}+1}z^{2k\hspace{0.01in}+1}}{(q\hspace
{0.01in};q)_{2k}}\,\frac{(q^{-M})^{2k\hspace{0.01in}+1}}{1-q^{2k\hspace
{0.01in}+1}}\nonumber\\
&  \qquad=\sum_{k\hspace{0.01in}=\hspace{0.01in}0}^{\infty}\frac
{(-1)^{k}\,(1-q)^{2k\hspace{0.01in}+1}(q^{-M}z)^{2k\hspace{0.01in}+1}%
}{(q\hspace{0.01in};q)_{2k\hspace{0.01in}+1}}=\sin_{q}(q^{-M}z).
\label{SumForqTri1}%
\end{align}
In the third step of the above calculation, we used the following identities:%
\begin{align}
\sum_{m\hspace{0.01in}=\hspace{0.01in}-M}^{\infty}q^{m(2k\hspace{0.01in}+1)}
&  =\sum_{m\hspace{0.01in}=\hspace{0.01in}0}^{\infty}q^{(m\hspace
{0.01in}-M)(2k\hspace{0.01in}+1)}=q^{-M(2k\hspace{0.01in}+1)}\sum
_{m\hspace{0.01in}=\hspace{0.01in}0}^{\infty}q^{m(2k\hspace{0.01in}%
+1)}\nonumber\\
&  =\frac{q^{-M(2k\hspace{0.01in}+1)}}{1-q^{2k\hspace{0.01in}+1}}.
\end{align}
By similar reasoning, we can prove the following identity
\cite{olshanetsky1998q.alg}:%
\begin{equation}
\int\nolimits_{0}^{\hspace{0.01in}zq^{-M}}\text{d}_{q}x\,\sin{}_{q}%
(x)=(1-q)z\sum_{m\hspace{0.01in}=\hspace{0.01in}-M}^{\infty}q^{m}\sin
_{q}(q^{m}z)=1-\cos_{q}(q^{-M}z). \label{SumForqTri2}%
\end{equation}
If $M$ becomes infinite, Eqs.~(\ref{SumForqTri1}) and (\ref{SumForqTri2})
imply%
\begin{align}
\int_{0}^{\hspace{0.01in}\infty}\text{d}_{q}z\,\cos_{q}(z)  &  =(1-q)z\lim
_{M\rightarrow\hspace{0.01in}\infty}\sum_{m\hspace{0.01in}=\hspace{0.01in}%
-M}^{\infty}q^{m}\cos_{q}(q^{m}z)\nonumber\\
&  =\lim_{M\rightarrow\hspace{0.01in}\infty}\sin_{q}(q^{-M}z)=0
\end{align}
and%
\begin{align}
\int_{0}^{\hspace{0.01in}\infty}\text{d}_{q}z\,\sin_{q}(z)  &  =(1-q)z\lim
_{M\rightarrow\hspace{0.01in}\infty}\sum_{m\hspace{0.01in}=\hspace{0.01in}%
-M}^{\infty}q^{m}\sin_{q}(q^{m}z)\nonumber\\
&  =1-\lim_{M\rightarrow\hspace{0.01in}\infty}\cos_{q}(q^{-M}z)=1.
\label{IntSinQFktEin}%
\end{align}
In the above calculations, the final identity holds since the $q$%
-trigo\-nomet\-ric functions $\sin_{q}$ and $\cos_{q}$ vanish at infinity [cf.
the inequalities in (\ref{AbsqTriSinCos}) of the previous chapter].

\section{Distributions\label{KapQDisEin}}

$\mathcal{M}(x)$ denotes the set of functions with the domain being the
following $q$-lat\-tice:%
\begin{equation}
\mathbb{G}_{q,\hspace{0.01in}x_{0}}\mathbb{=}\left\{  \pm\hspace{0.01in}%
x_{0}\hspace{0.01in}q^{m}|\,m\in\mathbb{Z}\right\}  \cup\{0\}.
\end{equation}
Using the $q$-in\-te\-gral over $\left]  -\infty,\infty\right[  $ [cf.
Eq.~(\ref{UneJackIntAll}) of Chap.~\ref{KapQIntTrig}] we can introduce a
$q$\textbf{-scalar product} on $\mathcal{M}(x)$ \cite{Ruffing:1996ph}
[$\hspace{0.01in}f,g\in\mathcal{M}(x)$]:%
\begin{equation}
\left\langle f,g\right\rangle _{q}=\int\nolimits_{-x.\infty}^{\hspace
{0.01in}x.\infty}\text{d}_{q}z\,\overline{f(z)}\,g(z)=\sum_{\varepsilon
\hspace{0.01in}=\hspace{0.01in}\pm}\,\sum_{j=-\infty}^{\infty}(1-q)\hspace
{0.01in}q^{\hspace{0.01in}j}x\,\overline{f(\varepsilon q^{\hspace{0.01in}j}%
x)}\,g(\varepsilon q^{\hspace{0.01in}j}x). \label{SkaProqDefEinBraLin}%
\end{equation}
Here, $\bar{f}$ denotes the complex conjugate of $f$. The $q$-scalar product
induces a $q$\textbf{-norm} on $\mathcal{M}(x)$%
\ \cite{Jambor:2004ph,Lavagno:2009vg}:%
\begin{equation}
\left\Vert f\right\Vert _{q}^{2}=\left\langle f,\hspace{0.01in}f\right\rangle
_{q}=\sum_{\varepsilon\hspace{0.01in}=\hspace{0.01in}\pm}\,\sum_{j=-\infty
}^{\infty}(1-q)\hspace{0.01in}q^{\hspace{0.01in}j}x\left\vert f(\varepsilon
q^{\hspace{0.01in}j}x)\right\vert ^{2}.
\end{equation}
The sum in the above equation takes the value zero if and only if
$f(\pm\hspace{0.01in}q^{m}x)=0$ for all $m\in\mathbb{Z}$. Any function that
vanishes on the $q$-lat\-tice can be identified with the zero function. This
observations ensure that the $q$-scalar product is positive definite. The set
of functions with finite $q$-norm forms a $q$\textbf{-ver\-sion of the Hilbert
space of square-integrable functions}
\cite{Cerchiai:1999,Hinterding:2000ph,Wess:math-ph9910013}:%
\begin{equation}
L_{q}^{2}=\{f\,|\,\left\langle f,\hspace{0.01in}f\right\rangle _{q}<\infty\}.
\end{equation}

Next, we consider all functions of $\mathcal{M}(x)$ that are continuous at the
origin of the lattice $\mathbb{G}_{q,\hspace{0.01in}x_{0}}$ ($0<q<1$):%
\begin{equation}
\lim_{m\hspace{0.01in}\rightarrow\hspace{0.01in}\infty}f(\pm\hspace
{0.01in}x_{0}\hspace{0.01in}q^{m})=f(0) \label{SteUrsGit}%
\end{equation}
In addition, there is a set of positive constants $C_{k,l}(q)$ such that the
following conditions hold ($k,l\in\mathbb{N}_{0}$):%
\begin{equation}
\left\vert \,x^{k}D_{q}^{l}f(x)\right\vert \leq C_{k,l}(q)\quad\text{for}\quad
x\in\mathbb{G}_{q,\hspace{0.01in}x_{0}}. \label{SchAbfFktBed}%
\end{equation}
These conditions imply:%
\begin{equation}
\lim_{m\hspace{0.01in}\rightarrow\hspace{0.01in}\infty}(D_{q}^{l}f)(\pm
\hspace{0.01in}x\hspace{0.01in}q^{-m})=0. \label{LimQTesFkt}%
\end{equation}
The functions of $\mathcal{M}(x)$ with the above properties make up the new
set $\mathcal{S}(x)$, which is a $q$\textbf{-ver\-sion of the space of test
functions for tempered distributions}. For these $q$-de\-formed test
functions, we introduce a family of norms:%
\begin{equation}
\left\Vert f\right\Vert _{k,l}=\max_{x\hspace{0.01in}\in\hspace{0.01in}%
\mathbb{G}_{q,x_{0}}}\left\vert x^{k}D_{q}^{l}f(x)\right\vert .
\end{equation}
This family of norms induces a topology on the set $\mathcal{S}(x)$
\cite{olshanetsky1998q.alg}.

A continuous linear functional $l:\mathcal{S}(x)\rightarrow\mathbb{C}$ is
called a $q$\textbf{-dis\-tri\-bu\-tion}. A $q$-dis\-tri\-bu\-tion $l$ is
called regular if there is a function $f\in\mathcal{M}(x)$ with the following
property:%
\begin{equation}
l(\hspace{0.01in}g)=\left\langle f,g\right\rangle _{q}=\int\nolimits_{-\infty
}^{\hspace{0.01in}\infty}\text{d}_{q}x\,\overline{f(x)}\,g(x).
\label{DefRegDisqDef}%
\end{equation}

We can use regular $q$-dis\-tri\-bu\-tions to obtain new $q$%
-dis\-tri\-bu\-tions. Let $f_{n}(x)$ be a sequence of functions that determine
a sequence of regular $q$-dis\-tri\-bu\-tions. If the limit%
\begin{equation}
l(\hspace{0.01in}g)=\lim_{n\hspace{0.01in}\rightarrow\hspace{0.01in}\infty
}\left\langle f_{n},g\right\rangle _{q}%
\end{equation}
exists for all $g\in\mathcal{S}(x)$, the mapping $l:g\mapsto l(\hspace
{0.01in}g)$ will be a $q$-dis\-tri\-bu\-tion \cite{olshanetsky1998q.alg}.

Let $l$ be a regular $q$-dis\-tri\-bu\-tion represented by a function $f$. We
can calculate the Jackson derivatives of $l$ as follows ($k\in\mathbb{N},$
$g\in\mathcal{S}(x)$):%
\begin{align}
(D_{q}^{k}\hspace{0.01in}l)(\hspace{0.01in}g)  &  =\int\nolimits_{-\infty
}^{\hspace{0.01in}\infty}\text{d}_{q}x\,\overline{D_{q}^{k}f(x)}%
\,g(x)=\int\nolimits_{-\infty}^{\hspace{0.01in}\infty}\text{d}_{q}%
x\,\Big[D_{q}^{k}\hspace{0.01in}\overline{f(x)}\,\Big]\,g(x)\nonumber\\
&  =(-1)^{k}q^{-k(k\hspace{0.01in}-1)/2}\int\nolimits_{-\infty}^{\hspace
{0.01in}\infty}\text{d}_{q}x\,\overline{f(q^{k}x)}\,D_{q}^{k}g(x)\nonumber\\
&  =(-1)^{k}q^{-k(k\hspace{0.01in}+1)/2}\int\nolimits_{-\infty}^{\hspace
{0.01in}\infty}\text{d}_{q}x\,\overline{f(x)}\,(D_{q}^{k}g)(q^{-k}x).
\label{AblqDisHer}%
\end{align}
In the second last step, we used Eq.~(\ref{RegParIntBraLinK}) of
Chap.~\ref{KapQIntTrig}, and in the final step, we rescaled the integration
variable. We can simplify the result of Eq.~(\ref{AblqDisHer}) by applying the
following identity:%
\begin{equation}
(D_{q}^{k}g)(q^{-k}x)=q^{k(k\hspace{0.01in}-1)/2}D_{q^{-1}}^{k}g(x).
\end{equation}
This way, we obtain:%
\begin{equation}
(D_{q}^{k}\hspace{0.01in}l)(\hspace{0.01in}g)=\int\nolimits_{-\infty}%
^{\hspace{0.01in}\infty}\text{d}_{q}x\,\overline{f(x)}\big((-q^{-1}D_{q^{-1}%
})^{k}g\big)(x)=l\big((-q^{-1}D_{q^{-1}})^{k}g\big). \label{JacAblKOrdDis}%
\end{equation}

By introducing the functions \cite{olshanetsky1998q.alg}%
\begin{equation}
\phi_{m}^{\hspace{0.01in}\varepsilon}(x)=\left\{
\begin{array}
[c]{ll}%
1 & \text{for\quad}x=\varepsilon q^{m}x_{0},\\
0 & \text{for\quad}x\neq\varepsilon q^{m}x_{0},
\end{array}
\right.  \label{BasVekFktBraLin}%
\end{equation}
we can define \textbf{step functions on the }$q$\textbf{-lat\-tice}:%
\begin{equation}
\theta_{q}^{\pm}(z)=\sum_{n\hspace{0.01in}=\hspace{0.01in}-\infty}^{\infty
}\phi_{n}^{\pm}(z) \label{DefTheQDefEin}%
\end{equation}
The functions $\theta_{q}^{\hspace{0.01in}\varepsilon}$ determine regular
$q$-dis\-tri\-bu\-tions ($\hspace{0.01in}g(x)\in\mathcal{S}(x)$, $0<q<1$):%
\begin{equation}
\theta_{q}^{\hspace{0.01in}\varepsilon}(\hspace{0.01in}g)=\left\langle
\theta_{q}^{\hspace{0.01in}\varepsilon},g\right\rangle _{q}=\int_{0}%
^{\hspace{0.01in}\infty}\text{d}_{q}z\,g(\varepsilon z)=(1-q)\hspace
{0.01in}x_{0}\sum_{m\hspace{0.01in}=\hspace{0.01in}-\infty}^{\infty}%
q^{m}g(\varepsilon x_{0}\hspace{0.01in}q^{m}). \label{QStuFktDis}%
\end{equation}

Next, we define a $q$\textbf{-ver\-sion of the delta distribution}:%
\begin{equation}
\delta_{q}(\hspace{0.01in}g)=\left\langle \delta_{q},g\right\rangle _{q}%
=\lim_{m\hspace{0.01in}\rightarrow\hspace{0.01in}\infty}\frac{g(x_{0}%
\hspace{0.01in}q^{m})+g(-x_{0}\hspace{0.01in}q^{m})}{2}=g(0).
\label{DefDelQEinDim}%
\end{equation}
As the function $g$ is continuous at the origin, it holds%
\begin{equation}
2\hspace{0.01in}\delta_{q}(\hspace{0.01in}g)=\delta_{q}^{+}(\hspace
{0.01in}g)+\delta_{q}^{-}(\hspace{0.01in}g) \label{ZerQDelDis}%
\end{equation}
where ( $\varepsilon\in\{+,-\}$)%
\begin{equation}
\delta_{q}^{\hspace{0.01in}\varepsilon}(\hspace{0.01in}g)=\left\langle
\delta_{q}^{\hspace{0.01in}\varepsilon},g\right\rangle _{q}=\lim
_{m\hspace{0.01in}\rightarrow\hspace{0.01in}\infty}g(\varepsilon x_{0}%
\hspace{0.01in}q^{m})=g(0^{\varepsilon}).
\end{equation}

Using Eq.~(\ref{JacAblKOrdDis}), we can calculate the \textbf{Jackson
derivative for the }$q$\textbf{-dis\-tri\-bu\-tion} $\theta_{q}^{\pm}$
($0<q<1$, $g\in\mathcal{S}(x)$, $\varepsilon\in\{+,-\}$):%
\begin{align}
(D_{q}\hspace{0.01in}\theta_{q}^{\hspace{0.01in}\varepsilon})(\hspace
{0.01in}g)  &  =-q^{-1}\langle\theta_{q}^{\hspace{0.01in}\varepsilon
},D_{q^{-1}}g\rangle_{q}\nonumber\\
&  =-(q^{-1}-1)\,x_{0}\sum_{m\hspace{0.01in}=\hspace{0.01in}-\infty}^{\infty
}q^{m}\,\frac{g(\varepsilon x_{0}\hspace{0.01in}q^{m})-g(\varepsilon
x_{0}\hspace{0.01in}q^{m\hspace{0.01in}-1})}{(1-q^{-1})\,\varepsilon
x_{0}\hspace{0.01in}q^{m}}\nonumber\\
&  =\varepsilon\sum_{m\hspace{0.01in}=\hspace{0.01in}-\infty}^{\infty}\left[
\hspace{0.01in}g(\varepsilon x_{0}\hspace{0.01in}q^{m})-g(\varepsilon
x_{0}\hspace{0.01in}q^{m\hspace{0.01in}-1})\right] \nonumber\\
&  =\varepsilon\lim_{M\rightarrow\hspace{0.01in}\infty}\sum_{m\hspace
{0.01in}=\hspace{0.01in}-M}^{M}\left[  \hspace{0.01in}g(\varepsilon
x_{0}\hspace{0.01in}q^{m})-g(\varepsilon x_{0}\hspace{0.01in}q^{m\hspace
{0.01in}-1})\right] \nonumber\\
&  =\varepsilon\lim_{M\rightarrow\hspace{0.01in}\infty}\left[  \hspace
{0.01in}g(\varepsilon x_{0}\hspace{0.01in}q^{M})-g(\varepsilon x_{0}%
\hspace{0.01in}q^{-M-1})\right]  =\varepsilon\hspace{0.01in}g(0^{\varepsilon
}). \label{AblDisqTheFkt}%
\end{align}
In the last step of the calculation above, we used the boundary conditions in
Eqs.~(\ref{SteUrsGit}) and (\ref{LimQTesFkt}). Due to $g(0^{\varepsilon
})=\delta_{q}^{\hspace{0.01in}\varepsilon}(\hspace{0.01in}g)$, the result in
Eq.~(\ref{AblDisqTheFkt}) gives us the following derivative rule:%
\begin{equation}
(D_{q}\hspace{0.01in}\theta_{q}^{\hspace{0.01in}\varepsilon})(\hspace
{0.01in}g)=\varepsilon\hspace{0.01in}\delta_{q}^{\hspace{0.01in}\varepsilon
}(\hspace{0.01in}g). \label{AblTheFktEin}%
\end{equation}
With this relation, we can represent both parts of the $q$-delta distribution
as the limit of regular $q$-dis\-tri\-bu\-tions. To this end, we proceed as
follows:%
\begin{equation}
\delta_{q}^{\hspace{0.01in}\varepsilon}(\hspace{0.01in}g)=\varepsilon
\hspace{0.01in}(D_{q}\hspace{0.01in}\theta_{q}^{\hspace{0.01in}\varepsilon
})(\hspace{0.01in}g)=\left\langle \varepsilon D_{q}\hspace{0.01in}\theta
_{q}^{\hspace{0.01in}\varepsilon},g\right\rangle _{q}=\int\nolimits_{-\infty
}^{\hspace{0.01in}\infty}\text{d}_{q}x\,\overline{\varepsilon D_{q}%
\hspace{0.01in}\theta_{q}^{\hspace{0.01in}\varepsilon}(x)}\,g(x).
\label{GreQDelDis}%
\end{equation}
It remains to calculate the Jackson derivative for $\varepsilon\hspace
{0.01in}\theta_{q}^{\hspace{0.01in}\varepsilon}(z)$:%
\begin{align}
\varepsilon D_{q}\hspace{0.01in}\theta_{q}^{\hspace{0.01in}\varepsilon}(x)  &
=\varepsilon\sum_{n\hspace{0.01in}=\hspace{0.01in}-\infty}^{\infty}%
D_{q}\hspace{0.01in}\phi_{n}^{\hspace{0.01in}\varepsilon}(x)=\varepsilon
\sum_{n\hspace{0.01in}=\hspace{0.01in}-\infty}^{\infty}\frac{\phi_{n}%
^{\hspace{0.01in}\varepsilon}(x)-\phi_{n}^{\hspace{0.01in}\varepsilon}%
(qx)}{(1-q)\hspace{0.01in}x}\nonumber\\
&  =\varepsilon\hspace{0.01in}(1-q)^{-1}x^{-1}\sum_{n\hspace{0.01in}%
=\hspace{0.01in}-\infty}^{\infty}\left(  \phi_{n}^{\hspace{0.01in}\varepsilon
}(x)-\phi_{n\hspace{0.01in}-1}^{\hspace{0.01in}\varepsilon}(x)\right)
\nonumber\\
&  =\varepsilon\hspace{0.01in}(1-q)^{-1}x^{-1}\lim_{N\rightarrow
\hspace{0.01in}\infty}\sum_{n\hspace{0.01in}=\hspace{0.01in}-N}^{N}\left(
\phi_{n}^{\hspace{0.01in}\varepsilon}(x)-\phi_{n\hspace{0.01in}-1}%
^{\hspace{0.01in}\varepsilon}(x)\right) \nonumber\\
&  =\varepsilon\hspace{0.01in}(1-q)^{-1}x^{-1}\lim_{N\rightarrow
\hspace{0.01in}\infty}\left(  \phi_{N}^{\hspace{0.01in}\varepsilon}%
(x)-\phi_{-N-1}^{\hspace{0.01in}\varepsilon}(x)\right)  .
\end{align}
Plugging this result into Eq.~(\ref{GreQDelDis}), we finally get:%
\begin{align}
\delta_{q}^{\hspace{0.01in}\varepsilon}(\hspace{0.01in}g)  &  =\lim
_{N\rightarrow\hspace{0.01in}\infty}\int\nolimits_{-\infty}^{\hspace
{0.01in}\infty}\text{d}_{q}x\,\frac{\varepsilon\left(  \phi_{N}^{\hspace
{0.01in}\varepsilon}(x)-\phi_{-N-1}^{\hspace{0.01in}\varepsilon}(x)\right)
}{(1-q)\hspace{0.01in}x}\,g(x)\nonumber\\
&  =\lim_{N\rightarrow\hspace{0.01in}\infty}\int\nolimits_{-\infty}%
^{\hspace{0.01in}\infty}\text{d}_{q}x\,\frac{\varepsilon\hspace{0.01in}%
\phi_{N}^{\hspace{0.01in}\varepsilon}(x)}{(1-q)\hspace{0.01in}x}%
\,g(x)-\lim_{N\rightarrow\infty}\int\nolimits_{-\infty}^{\hspace{0.01in}%
\infty}\text{d}_{q}x\,\frac{\varepsilon\hspace{0.01in}\phi_{-N-1}%
^{\hspace{0.01in}\varepsilon}(x)}{(1-q)\hspace{0.01in}x}\,g(x)\nonumber\\
&  =\lim_{N\rightarrow\hspace{0.01in}\infty}\int\nolimits_{-\infty}%
^{\hspace{0.01in}\infty}\text{d}_{q}x\,\frac{\varepsilon\hspace{0.01in}%
\phi_{N}^{\hspace{0.01in}\varepsilon}(x)}{(1-q)\hspace{0.01in}x}\,g(x).
\label{GreRegDisqDis}%
\end{align}
In the final step, we have omitted the second $q$-in\-te\-gral since it
vanishes due to Eq.~(\ref{LimQTesFkt}).

Next, we want to find the \textbf{Jackson derivatives of the }$q$%
\textbf{-delta distribution}. To this end, we perform the following
calculation ($k\in\mathbb{N}$):%
\begin{align}
(D_{q}^{k}\delta_{q}^{\hspace{0.01in}\varepsilon})(\hspace{0.01in}g)  &
=\varepsilon D_{q}^{k\hspace{0.01in}+1}\theta_{q}^{\hspace{0.01in}\varepsilon
}(\hspace{0.01in}g)=\varepsilon\hspace{0.01in}(-q^{-1})^{k}(D_{q}%
\hspace{0.01in}\theta_{q}^{\hspace{0.01in}\varepsilon})\big(D_{q^{-1}}%
^{k}g\big)\nonumber\\
&  =(-q^{-1})^{k}\hspace{0.01in}\delta_{q}^{\hspace{0.01in}\varepsilon
}\big(D_{q^{-1}}^{k}g\big)=(-q^{-1})^{k}\big(D_{q^{-1}}^{k}%
g\big)(0^{\varepsilon}). \label{AblQDelDis1Dim}%
\end{align}
The first two identities follow from Eq.~(\ref{GreQDelDis}) or
Eq.~(\ref{JacAblKOrdDis}). Due to Eq.~(\ref{ZerQDelDis}), the above result
yields the following formulas:%
\begin{equation}
(D_{q}^{k}\delta_{q})(\hspace{0.01in}g)=(-q^{-1})^{k}\hspace{0.01in}\delta
_{q}\big(D_{q^{-1}}^{k}g\big)=(-q^{-1})^{k}\big(D_{q^{-1}}^{k}g\big)(0).
\end{equation}

We want to study how the $q$-delta distribution behaves under $q$%
-trans\-la\-tions. We assume that we have a function $\delta_{q}(x)$ with the
following property:%
\begin{equation}
\delta_{q}(\hspace{0.01in}g)=\int\nolimits_{-\infty}^{\hspace{0.01in}\infty
}\text{d}_{q}\hspace{0.01in}x\,\delta_{q}(x)\,g(x)=g(0).
\label{ChaIdeqDelQKal1dim}%
\end{equation}
Using the $q$-Taylor formula in Eqs.~(\ref{QTayForTyp2}) and
(\ref{QTayForTyp3}) of Chap.~\ref{KapqAnaExpTriFkt}, we get:%
\begin{align}
\delta_{q}(x\,\bar{\oplus}\,(\bar{\ominus}\,a))  &  =\sum_{k\hspace
{0.01in}=\hspace{0.01in}0}^{\infty}\frac{q^{k(k\hspace{0.01in}-1)/2}%
}{[[k]]_{q}!}\left(  D_{q}^{k}\hspace{0.01in}\delta_{q}(x)\right)
\hspace{0.01in}(-a)^{k},\nonumber\\
\delta_{q}((\bar{\ominus}\,a)\,\bar{\oplus}\,x)  &  =\sum_{k\hspace
{0.01in}=\hspace{0.01in}0}^{\infty}\frac{q^{k(k\hspace{0.01in}-1)/2}%
}{[[k]]_{q}!}\,(-a)^{k}\hspace{0.01in}D_{q}^{k}\hspace{0.01in}\delta_{q}(x).
\label{DefVerDelFktQ1Dim}%
\end{align}
Using the expressions above, we do the following calculation :%
\begin{align}
&  \int_{-\infty}^{\hspace{0.01in}\infty}\text{d}_{q}\hspace{0.01in}%
x\,\delta_{q}((\bar{\ominus}\,qa)\,\bar{\oplus}\,x)\,f(x)=\sum_{k\hspace
{0.01in}=\hspace{0.01in}0}^{\infty}\frac{q^{k(k\hspace{0.01in}-1)/2}%
}{[[k]]_{q}!}\,(-qa)^{k}\int_{-\infty}^{\hspace{0.01in}\infty}\text{d}%
_{q}x\,(D_{q}^{k}\delta_{q}(x))\hspace{0.01in}f(x)\nonumber\\
&  \qquad=\sum_{k\hspace{0.01in}=\hspace{0.01in}0}^{\infty}\frac
{q^{k(k\hspace{0.01in}-1)/2}}{[[k]]_{q}!}\,(-qa)^{k}\,(-1)^{k}\,q^{-k(k\hspace
{0.01in}-1)/2}\int_{-\infty}^{\hspace{0.01in}\infty}\text{d}_{q}x\,\delta
_{q}(q^{k}x)\,D_{q}^{k}f(x)\nonumber\\
&  \qquad=\sum_{k\hspace{0.01in}=\hspace{0.01in}0}^{\infty}\frac{1}%
{[[k]]_{q}!}\,a^{k}\int_{-\infty}^{\hspace{0.01in}\infty}\text{d}_{q}%
x\,\delta_{q}(x)\,(D_{q}^{k}f)(q^{-k}x)\nonumber\\
&  \qquad=\sum_{k\hspace{0.01in}=\hspace{0.01in}0}^{\infty}\frac{1}%
{[[k]]_{q}!}\,a^{k}(D_{q}^{k}f)(0)=f(a). \label{HerChaIdeDel1Dim}%
\end{align}
In the second step, we employed the rules for integration by parts [cf.
Eq.~(\ref{RegParIntBraLinK}) of the previous chapter]. In the third step, we
changed the scaling of the integration variable. The second last step follows
from Eq.~(\ref{ChaIdeqDelQKal1dim}), and the final step results from
Eqs.~(\ref{qTraVerGerUnKon}) and (\ref{NeuEleQTra}) of
Chap.~\ref{KapqAnaExpTriFkt}. In the same way, we can derive the following
identity:%
\begin{equation}
\int_{-\infty}^{\hspace{0.01in}\infty}\text{d}_{q}x\hspace{0.01in}%
f(x)\,\delta_{q}(x\,\bar{\oplus}\,(\bar{\ominus}\,qa))=f(a).
\label{ChaIdeDelQFkt1DimTra}%
\end{equation}

The $q$-delta functions are related to the functions $\phi_{m}^{\pm}(x)$ in
Eq.~(\ref{BasVekFktBraLin}). Concretely, the identities in
Eqs.~(\ref{HerChaIdeDel1Dim}) and (\ref{ChaIdeDelQFkt1DimTra}) are satisfied
for $a\in\mathbb{G}_{q,\hspace{0.01in}x_{0}}$ if it holds ($0<q<1$):%
\begin{equation}
\left.  \delta_{q}((\bar{\ominus}\,qa)\,\bar{\oplus}\,x)\right\vert
_{a\hspace{0.01in}=\hspace{0.01in}\pm\hspace{0.01in}x_{0}q^{m}}=\left.
\delta_{q}(x\,\bar{\oplus}\,(\bar{\ominus}\,qa))\right\vert _{a\hspace
{0.01in}=\hspace{0.01in}\pm\hspace{0.01in}x_{0}q^{m}}=\frac{\phi_{m}^{\pm}%
(x)}{(1-q)\hspace{0.01in}x_{0}\hspace{0.01in}q^{m}}. \label{DarQDelKroDel1Dim}%
\end{equation}
The functions $\phi_{m}^{\pm}(x)$ can be represented by infinitely
differentiable functions. To this end, we introduce the following function:%
\begin{equation}
\psi(x,a)=\left\{
\begin{array}
[c]{ll}%
0 & \text{if\ }x<-a,\\
\operatorname*{e}\nolimits^{-\operatorname*{e}\nolimits^{x/(x^{2}%
-\hspace{0.01in}a^{2})}}a^{2}\operatorname*{e} & \text{if\ }-a\leq x\leq a,\\
a^{2}\operatorname*{e} & \text{if\ }a<x.
\end{array}
\right.
\end{equation}
We calculate the partial derivative of $\psi(x,a)$ with respect to $x$:%
\begin{equation}
\frac{\partial}{\partial x}\psi(x,a)=\left\{
\begin{array}
[c]{ll}%
0 & \text{if }x<-a,\\
\operatorname*{e}\nolimits^{x/(x^{2}-\hspace{0.01in}a^{2})-\operatorname*{e}%
\nolimits^{x/(x^{2}-\hspace{0.01in}a^{2})}}\frac{x^{2}+\hspace{0.01in}a^{2}%
}{(x^{2}-\hspace{0.01in}a^{2})^{2}}\,a^{2}\operatorname*{e} & \text{if }-a\leq
x\leq a,\\
0 & \text{if }a<x.
\end{array}
\right.
\end{equation}
This function is infinitely differentiable with respect to $x$ and has a
bell-shaped graph that is symmetrical to the axis $x=0$
\cite{Grossmann:1988ah}. Accordingly, $\partial_{x}\psi(x,a)$ takes its
maximum value at $x=0$, namely $\partial_{x}\psi(0,a)=1$. Furthermore, we
define a function that assigns to $y$ half the distance from $y$ to the next
lattice point $q\hspace{0.01in}y$ or $q^{-1}y$ :%
\begin{equation}
b_{q}(\hspace{0.01in}y)=\frac{1}{2}\min\left\{  |(1-q^{-1})\hspace
{0.01in}y|,|(1-q)\hspace{0.01in}y|\right\}  .
\end{equation}
Using $\partial_{x}\psi(x,a)$ and $b_{q}(\hspace{0.01in}y)$, we can introduce
the following function:%
\begin{equation}
\phi_{q}(x,y)=\frac{\partial}{\partial x}\psi(x-y,b_{q}(\hspace{0.01in}y)).
\label{DefKonDelFkt1Dim}%
\end{equation}
This function has the properties%
\begin{equation}
\phi_{q}(x,x)=1\quad\text{for}\quad x\in\mathbb{R}%
\end{equation}
and%
\begin{equation}
\phi_{q}(x,y)=0\quad\text{for}\quad x\in\mathbb{R}\backslash\left[
\hspace{0.01in}y-b_{q}(\hspace{0.01in}y),y+b_{q}(\hspace{0.01in}y)\right]  .
\end{equation}
This enables us to identify the function $\phi_{m}^{\hspace{0.01in}%
\varepsilon}(x)$ on the one-di\-men\-sion\-al $q$-lat\-tice $\mathbb{G}%
_{q,\hspace{0.01in}x_{0}}$ with an infinitely differentiable function:%
\begin{equation}
\phi_{m}^{\hspace{0.01in}\varepsilon}(x)=\phi_{q}(x,\varepsilon x_{0}q^{m}).
\end{equation}
Finally, the result above leads to the following representations [cf.
Eq.~(\ref{DarQDelKroDel1Dim})]:%
\begin{equation}
\delta_{q}((\bar{\ominus}\,qa)\,\bar{\oplus}\,x)=\delta_{q}(x\,\bar{\oplus
}\,(\bar{\ominus}\,qa))=\frac{\phi_{q}(x,a)}{(1-q)\hspace{-0.01in}\left\vert
a\right\vert }=\frac{\phi_{q}(x,a)}{(1-q)\hspace{-0.01in}\left\vert
x\right\vert }. \label{DarQDelKroDelAlg1Dim}%
\end{equation}

The step functions $\theta_{q}^{\varepsilon}$ given in
Eq.~(\ref{DefTheQDefEin}) are $q$-pe\-ri\-od\-ic, as we can seen from the
following calculation ($\varepsilon\in\{+,-\}$):%
\begin{equation}
\theta_{q}^{\hspace{0.01in}\varepsilon}(qz)=\sum_{n\hspace{0.01in}%
=\hspace{0.01in}-\infty}^{\infty}\phi_{n}^{\hspace{0.01in}\varepsilon
}(qz)=\sum_{n\hspace{0.01in}=\hspace{0.01in}-\infty}^{\infty}\phi
_{n\hspace{0.01in}-1}^{\hspace{0.01in}\varepsilon}(z)=\sum_{n\hspace
{0.01in}=\hspace{0.01in}-\infty}^{\infty}\phi_{n}^{\hspace{0.01in}\varepsilon
}(z)=\theta_{q}^{\hspace{0.01in}\varepsilon}(z). \label{qPerTheFkt}%
\end{equation}
If we introduce a scaling operator $\Lambda$ with%
\begin{equation}
\Lambda\hspace{0.01in}f(x)=f(qx), \label{ActScaOpe}%
\end{equation}
Eq.~(\ref{qPerTheFkt}) implies that the step functions $\theta_{q}%
^{\varepsilon}$ are invariant under the action of $\Lambda$ ($n\in\mathbb{N}$,
$\varepsilon\in\{+,-\}$):%
\begin{equation}
\Lambda\hspace{-0.01in}^{n}\hspace{0.01in}\theta_{q}^{\hspace{0.01in}%
\varepsilon}(x)=\theta_{q}^{\hspace{0.01in}\varepsilon}(q^{n}x)=\theta
_{q}^{\hspace{0.01in}\varepsilon}(x).
\end{equation}

Under the action of the scaling operator, the $q$-delta function behaves as
follows:%
\begin{equation}
\Lambda\hspace{-0.01in}^{n}\hspace{0.01in}\delta_{q}(x)=\delta_{q}%
(q^{n}x)=q^{-n}\delta_{q}(x). \label{SkaDelFktEinQ}%
\end{equation}
We prove these identities by the following calculation:%
\begin{gather}
\int_{-\infty}^{\hspace{0.01in}\infty}\text{d}_{q}x\left[  \hspace
{0.01in}\Lambda\hspace{-0.01in}^{n}\delta_{q}(x)\right]  \hspace
{0.01in}g(x)=\int_{-\infty}^{\hspace{0.01in}\infty}\text{d}_{q}\hspace
{0.01in}x\,\delta_{q}(q^{n}x)\,g(x)\nonumber\\
=q^{-n}\int_{-\infty}^{\hspace{0.01in}\infty}\text{d}_{q}x\,\delta
_{q}(x)\,g(q^{-n}x)=q^{-n}g(0)=\int_{-\infty}^{\hspace{0.01in}\infty}%
\text{d}_{q}p\,q^{-n}\delta_{q}(x)\,g(x).
\end{gather}
The result above, which holds for any function $g\in\mathcal{S}(x)$, confirms
the identity in Eq.~(\ref{SkaDelFktEinQ}).

We regain the shifted $q$-delta functions given in
Eq.~(\ref{DarQDelKroDelAlg1Dim}) by applying a Jackson derivative to a step
function on the $q$-lat\-tice. To show this, we perform the following
calculation [$\hspace{0.01in}g(x)\in\mathcal{S}(x),$ $\varepsilon\in\{+,-\},$
$m\in\mathbb{Z}$]:%
\begin{align}
&  \int_{-\infty}^{\hspace{0.01in}\infty}\text{d}_{q}x\,g(x)\,D_{q}%
\sum_{n\hspace{0.01in}=\hspace{0.01in}-\infty}^{m}\phi_{n}^{\hspace
{0.01in}\varepsilon}(x)=\int_{-\infty}^{\hspace{0.01in}\infty}\text{d}%
_{q}x\,g(x)\sum_{n\hspace{0.01in}=\hspace{0.01in}-\infty}^{m}\frac{\phi
_{n}^{\hspace{0.01in}\varepsilon}(x)-\phi_{n}^{\hspace{0.01in}\varepsilon
}(qx)}{(1-q)\hspace{0.01in}x}\nonumber\\
&  \qquad=\int_{-\infty}^{\hspace{0.01in}\infty}\text{d}_{q}x\,\frac
{g(x)}{(1-q)\hspace{0.01in}x}\sum_{n\hspace{0.01in}=\hspace{0.01in}-\infty
}^{m}\left[  \phi_{n}^{\hspace{0.01in}\varepsilon}(x)-\phi_{n-1}%
^{\hspace{0.01in}\varepsilon}(x)\right] \nonumber\\
&  \qquad=\int_{-\infty}^{\hspace{0.01in}\infty}\text{d}_{q}x\,g(x)\,\frac
{\phi_{m}^{\hspace{0.01in}\varepsilon}(x)}{(1-q)\hspace{0.01in}x}%
-\lim_{k\hspace{0.01in}\rightarrow\hspace{0.01in}-\infty}\int_{-\infty
}^{\hspace{0.01in}\infty}\text{d}_{q}x\,g(x)\,\frac{\phi_{k}^{\hspace
{0.01in}\varepsilon}(x)}{(1-q)\hspace{0.01in}x}\nonumber\\
&  \qquad=\varepsilon\hspace{0.01in}g(\varepsilon x_{0}\hspace{0.01in}%
q^{m})-\varepsilon\lim_{k\hspace{0.01in}\rightarrow\hspace{0.01in}-\infty
}g(\varepsilon x_{0}\hspace{0.01in}q^{k})=\varepsilon\hspace{0.01in}%
g(\varepsilon x_{0}\hspace{0.01in}q^{m}).
\end{align}
Moreover, Eq.~(\ref{ChaIdeDelQFkt1DimTra}) implies:%
\begin{equation}
\int_{-\infty}^{\hspace{0.01in}\infty}\text{d}_{q}x\,g(x)\hspace
{-0.01in}\left.  \delta_{q}(x\,\bar{\oplus}\,(\bar{\ominus}\,qa))\right\vert
_{a\hspace{0.01in}=\hspace{0.01in}\varepsilon x_{0}q^{m}}=g(\varepsilon
x_{0}\hspace{0.01in}q^{m}).
\end{equation}
Thus, the following identity holds in the set of regular $q$%
-dis\-tri\-bu\-tions:%
\begin{equation}
D_{q}\sum_{n\hspace{0.01in}=\hspace{0.01in}-\infty}^{m}\phi_{n}^{\hspace
{0.01in}\varepsilon}(x)=\varepsilon\left.  \delta_{q}(x\,\bar{\oplus}%
\,(\bar{\ominus}\,qa))\right\vert _{a\hspace{0.01in}=\hspace{0.01in}%
\varepsilon x_{0}q^{m}}. \label{AblStuDel1Dim}%
\end{equation}
We can identify the sum on the left-hand side of the above equation with a
shifted step function on the $q$-lattice:%
\begin{equation}
\sum_{n\hspace{0.01in}=\hspace{0.01in}-\infty}^{m}\phi_{n}^{\hspace
{0.01in}\varepsilon}(x)=\left.  \theta_{q}^{\hspace{0.01in}\varepsilon
}(x\,\bar{\oplus}\,(\bar{\ominus}\,qa))\right\vert _{a\hspace{0.01in}%
=\hspace{0.01in}\varepsilon x_{0}q^{m}}. \label{DarVerStu1Dim}%
\end{equation}
Combining Eq.~(\ref{AblStuDel1Dim}) and Eq.~(\ref{DarVerStu1Dim}), we finally
obtain:%
\begin{equation}
\left.  D_{q,\hspace{0.01in}x}\hspace{0.01in}\theta_{q}^{\hspace
{0.01in}\varepsilon}(x\,\bar{\oplus}\,(\bar{\ominus}\,qa))\right\vert
_{a\hspace{0.01in}=\hspace{0.01in}\varepsilon x_{0}q^{m}}=\left.
\varepsilon\hspace{0.01in}\delta_{q}(x\,\bar{\oplus}\,(\bar{\ominus
}\,qa))\right\vert _{a\hspace{0.01in}=\hspace{0.01in}\varepsilon x_{0}q^{m}}.
\end{equation}

To verify the identity in Eq.~(\ref{DarVerStu1Dim}), we do the following
calculation:%
\begin{align}
&  \int\nolimits_{-x_{0}.\infty}^{\hspace{0.01in}x_{0}.\infty}\text{d}%
_{q}x\sum_{n\hspace{0.01in}=\hspace{0.01in}-\infty}^{m}\phi_{n}^{\hspace
{0.01in}\varepsilon}(x)\,g(x)=(1-q)\sum_{j\hspace{0.01in}=\hspace
{0.01in}-\infty}^{\infty}x_{0}\hspace{0.01in}q^{\hspace{0.01in}j}%
\sum_{n\hspace{0.01in}=\hspace{0.01in}-\infty}^{m}\phi_{n}^{\hspace
{0.01in}\varepsilon}(\varepsilon x_{0}\hspace{0.01in}q^{\hspace{0.01in}%
j})\,g(\varepsilon x_{0}\hspace{0.01in}q^{\hspace{0.01in}j})\nonumber\\
&  \qquad=(1-q)\sum_{n\hspace{0.01in}=\hspace{0.01in}-\infty}^{m}x_{0}%
\hspace{0.01in}q^{n}\hspace{0.01in}g(\varepsilon x_{0}\hspace{0.01in}%
q^{n})=(1-q)\sum_{n\hspace{0.01in}=\hspace{0.01in}-\infty}^{0}x_{0}%
\hspace{0.01in}q^{m+n}\hspace{0.01in}g(\varepsilon x_{0}\hspace{0.01in}%
q^{m+n})\nonumber\\
&  \qquad=(1-q)\sum_{n\hspace{0.01in}=1}^{\infty}\,x_{0}\hspace{0.01in}%
q^{m\hspace{0.01in}+1-\hspace{0.01in}n}\hspace{0.01in}g(\varepsilon
x_{0}\hspace{0.01in}q^{m\hspace{0.01in}+1-\hspace{0.01in}n})=\varepsilon
\int\nolimits_{\varepsilon x_{0}q^{m+1}}^{\hspace{0.01in}\varepsilon
x_{0}.\infty}\text{d}_{q}x\,g(x). \label{HerZusSteDel1}%
\end{align}
Similarly to Eq.~(\ref{HerChaIdeDel1Dim}), we have:%
\begin{align}
&  \int\nolimits_{-x_{0}.\infty}^{\hspace{0.01in}x_{0}.\infty}\text{d}%
_{q}x\,g(x)\,\theta_{q}^{\hspace{0.01in}\varepsilon}(x\,\bar{\oplus}%
\,(\bar{\ominus}\,qa))=\nonumber\\
&  \qquad=\sum_{k\hspace{0.01in}=\hspace{0.01in}0}^{\infty}\frac
{q^{k(k\hspace{0.01in}-1)/2}}{[[k]]_{q}!}\,(-qa)^{k}\int\nolimits_{-x_{0}%
.\infty}^{\hspace{0.01in}x_{0}.\infty}\text{d}_{q}x\,g(x)\,D_{q}^{k}\theta
_{q}^{\hspace{0.01in}\varepsilon}(x)\nonumber\\
&  \qquad=\sum_{k\hspace{0.01in}=\hspace{0.01in}0}^{\infty}\frac{(qa)^{k}%
}{[[k]]_{q}!}\int\nolimits_{-x_{0}.\infty}^{\hspace{0.01in}x_{0}.\infty
}\text{d}_{q}\hspace{0.01in}x\,\theta_{q}^{\hspace{0.01in}\varepsilon}%
(q^{k}x)\,D_{q}^{k}g(x)\nonumber\\
&  \qquad=\int\nolimits_{-x_{0}.\infty}^{\hspace{0.01in}x_{0}.\infty}%
\text{d}_{q}x\,\theta_{q}^{\hspace{0.01in}\varepsilon}(x)\,g(x\,\bar{\oplus
}\,qa)=\varepsilon\int\nolimits_{0}^{\hspace{0.01in}\varepsilon x_{0}.\infty
}\text{d}_{q}x\,g(x\,\bar{\oplus}\,qa). \label{HerZusSteDel2}%
\end{align}
We rewrite the last expression of Eq.~(\ref{HerZusSteDel2}) as follows
\cite{1996q.alg.....8008K}:%
\begin{equation}
\varepsilon\int\nolimits_{0}^{\hspace{0.01in}\varepsilon x_{0}.\infty}%
\text{d}_{q}x\,g(x\,\bar{\oplus}\,qa)=\varepsilon\int\nolimits_{0\,\bar
{\oplus}\,qa}^{\hspace{0.01in}\varepsilon x_{0}.\infty}\text{d}_{q}%
x\,g(x)=\varepsilon\int\nolimits_{qa}^{\hspace{0.01in}\varepsilon x_{0}%
.\infty}\text{d}_{q}x\,g(x). \label{HerZusSteDel3}%
\end{equation}
Putting the results in Eqs.~(\ref{HerZusSteDel1})-(\ref{HerZusSteDel3})
together, we get the identity from which we can read off
Eq.~(\ref{DarVerStu1Dim}):%
\begin{align}
&  \int\nolimits_{-x_{0}.\infty}^{\hspace{0.01in}x_{0}.\infty}\text{d}%
_{q}x\sum_{n\hspace{0.01in}=\hspace{0.01in}-\infty}^{m}\phi_{n}^{\hspace
{0.01in}\varepsilon}(x)\,g(x)=\varepsilon\int\nolimits_{\varepsilon
x_{0}q^{m+1}}^{\hspace{0.01in}\varepsilon x_{0}.\infty}\text{d}_{q}%
x\,g(x)\nonumber\\
&  \qquad\qquad=\int\nolimits_{-x_{0}.\infty}^{\hspace{0.01in}x_{0}.\infty
}\text{d}_{q}x\,g(x)\left.  \theta_{q}^{\hspace{0.01in}\varepsilon}%
(x\,\bar{\oplus}\,(\bar{\ominus}\,qa))\right\vert _{a\hspace{0.01in}%
=\hspace{0.01in}\varepsilon x_{0}q^{m}}.
\end{align}

\section{Fourier transforms\label{KapEinDimQFouTraN}}

In this chapter, we consider $q$-ver\-sions of one-dimensional Fourier
transforms. The $q$-de\-formed Fourier transforms are mappings between the
function spaces $\mathcal{M}_{R}(x)$ and $\mathcal{M}(\hspace{0.01in}p)$. In
this respect, $x$ and $p$ play the role of position and momentum variables.

We use the $q$-de\-formed exponential in Eq.~(\ref{AusQExpEin}) of
Chap.~\ref{KapqAnaExpTriFkt}\ to introduce $q$\textbf{-ver\-sions of
one-di\-men\-sion\-al plane waves} \cite{Kempf:1994yd}:%
\begin{equation}
\exp_{q}(x|\text{i}p)=\sum_{k\hspace{0.01in}=\hspace{0.01in}0}^{\infty}%
\frac{x^{k}(\text{i}p)^{k}}{[[k]]_{q}!}. \label{DefEinQEbeWel}%
\end{equation}
With the help of the Jackson derivative $D_{q}$, we introduce the momentum
operator%
\begin{equation}
\hat{P}=\text{i}^{-1}D_{q}. \label{ImpOpe1DimUnk}%
\end{equation}
The $q$-de\-formed plane waves in Eq.~(\ref{DefEinQEbeWel}) are eigenfunctions
of $\hat{P}$:%
\begin{align}
\text{i}^{-1}D_{q,x}\triangleright\exp_{q}(x|\text{i}p)  &  =\sum
_{k\hspace{0.01in}=1}^{\infty}\frac{[[k]]}{[[k]]_{q}!}\,x^{k-1}(\text{i}%
p)^{k\hspace{0.01in}-1}p\\
&  =\exp_{q}(x|\text{i}p)\hspace{0.01in}p. \label{EigFktImp1Dim}%
\end{align}
Additionally, we need the following function:%
\begin{align}
\exp_{q^{-1}}(-\text{i}p|\hspace{0.01in}x)  &  =\sum_{k\hspace{0.01in}%
=\hspace{0.01in}0}^{\infty}\frac{q^{k(k\hspace{0.01in}-1)/2}(-\text{i}%
p)^{k}x^{k}}{[[k]]_{q}!}\nonumber\\
&  =\sum_{k\hspace{0.01in}=\hspace{0.01in}0}^{\infty}\frac{(-\text{i}%
p)^{k}\hspace{0.01in}x^{k}}{[[k]]_{q^{-1}}!}. \label{ReiEntqExpInv}%
\end{align}
The momentum operator acts on this function as follows:%
\begin{align}
\text{i}^{-1}D_{q,x}\triangleright\exp_{q^{-1}}(-\text{i}p|x)  &
=\text{i}^{-1}\sum_{k\hspace{0.01in}=1}^{\infty}\frac{(-\text{i}%
)^{k}q^{k(k\hspace{0.01in}-1)/2}}{[[k-1]]_{q}!}\hspace{0.01in}p^{k}%
x^{k\hspace{0.01in}-1}\nonumber\\
&  =-\exp_{q^{-1}}(-\text{i}p|qx)\hspace{0.01in}p. \label{AblJacImpInvExp}%
\end{align}
As was shown in Ref.~\cite{olshanetsky1998q.alg}, the $q$-de\-formed plane
waves satisfy \textbf{completeness relations and orthogonality conditions},
given by%
\begin{align}
\int\limits_{-x_{0}.\infty}^{x_{0}.\infty}\text{d}_{q}p\,\exp_{q}%
(x_{0}|\text{i}p)\hspace{0.01in}\exp_{q^{-1}}(-\text{i}p|qx)  &  =\left\{
\begin{array}
[c]{cl}%
\frac{4\Theta_{0}}{(1-q)x_{0}} & \text{f\"{u}r\quad}x=x_{0},\\
0 & \text{f\"{u}r\quad}x\neq x_{0},
\end{array}
\right. \label{IntExpExpInv}\\[0.06in]
\int\limits_{-p_{0}.\infty}^{p_{0}.\infty}\text{d}_{q}x\,\exp_{q^{-1}%
}(-\text{i}p_{0}|qx)\hspace{0.01in}\exp_{q}(x|\text{i}p)  &  =\left\{
\begin{array}
[c]{cl}%
\frac{4\Theta_{0}}{(1-q)p_{0}} & \text{f\"{u}r\quad}p=p_{0},\\
0 & \text{f\"{u}r\quad}p\neq p_{0},
\end{array}
\right.  \label{IntExpExpInv2}%
\end{align}
with [cf. Eq.~(\ref{DefTheZ}) in Chap.~\ref{KapQIntTrig}]%
\begin{equation}
\Theta_{0}=\Theta_{q}(x_{0}\hspace{0.01in}p_{0}). \label{DefTheNul}%
\end{equation}

Next, we introduce $q$-Fourier transforms between $\mathcal{M}(x)$ and
$\mathcal{M}(\hspace{0.01in}p)$ \cite{Kempf:1994yd,olshanetsky1998q.alg}:%
\begin{align}
\mathcal{F}_{q}(\phi)(\hspace{0.01in}p)  &  =\frac{1}{2\sqrt{\Theta_{0}}}%
\int_{-x_{0}.\infty}^{\hspace{0.01in}x_{0}.\infty}\text{d}_{q}x\,\phi
(x)\hspace{0.01in}\exp_{q}(x|\text{i}p),\nonumber\\[0.04in]
\mathcal{F}_{q}^{\hspace{0.01in}-1}(\psi)(x)  &  =\frac{1}{2\sqrt{\Theta_{0}}%
}\int_{-p_{0}.\infty}^{p_{0}.\infty}\text{d}_{q}p\,\psi(\hspace{0.01in}%
p)\hspace{0.01in}\exp_{q^{-1}}(-\text{i}p|qx). \label{DefQFouTraEin}%
\end{align}
We show that the mappings $\mathcal{F}_{q}$\textit{ and }$\mathcal{F}%
_{q}^{\hspace{0.01in}-1}$\textit{ are inverse to each other}. Let $\phi
_{n}^{\pm}$ and $\psi_{n}^{\pm}$ with $n\in\mathbb{Z}$ denote the basis
functions of $\mathcal{M}(x)$ and $\mathcal{M}(\hspace{0.01in}p)$,
respectively\ [cf. Eq.~(\ref{BasVekFktBraLin}) of the previous chapter]. It
suffices to prove the following identities:%
\begin{align}
(\mathcal{F}_{q}^{\hspace{0.01in}-1}\circ\mathcal{F}_{q})(\phi_{n}^{\pm})  &
=\phi_{n}^{\pm},\nonumber\\
(\mathcal{F}_{q}\circ\mathcal{F}_{q}^{\hspace{0.01in}-1})(\psi_{n}^{\pm})  &
=\psi_{n}^{\pm}. \label{IdeInvFouTraBasFktqEin}%
\end{align}
To prove the first of the two identities above, we do the following
calculation:%
\begin{align}
&  (\mathcal{F}_{q}^{\hspace{0.01in}-1}\circ\mathcal{F}_{q})(\phi_{n}^{\pm
}(x))=\nonumber\\
&  \qquad=\frac{1}{4\Theta_{0}}\int_{-p_{0}.\infty}^{p_{0}.\infty}\text{d}%
_{q}p\,\int_{-x_{0}.\infty}^{\hspace{0.01in}x_{0}.\infty}\text{d}_{q}\xi
\,\phi_{n}^{\pm}(\xi)\hspace{0.01in}\exp_{q}(\xi|\text{i}p)\hspace{0.01in}%
\exp_{q^{-1}}(-\text{i}p|qx)\nonumber\\
&  \qquad=\frac{(1-q)\hspace{0.01in}x_{0}\hspace{0.01in}q^{n}}{4\Theta_{0}%
}\int_{-p_{0}.\infty}^{\hspace{0.01in}p_{0}.\infty}\text{d}_{q}p\,\exp_{q}%
(\pm\hspace{0.01in}x_{0}\hspace{0.01in}q^{n}|\hspace{0.01in}\text{i}%
p)\hspace{0.01in}\exp_{q^{-1}}(-\text{i}p|qx)\nonumber\\
&  \qquad=\frac{(1-q)\hspace{0.01in}x_{0}}{4\Theta_{0}}\int_{-p_{0}.\infty
}^{\hspace{0.01in}p_{0}.\infty}\text{d}_{q}p\,\exp_{q}(x_{0}|\text{i}%
p)\hspace{0.01in}\exp_{q^{-1}}(-\text{i}p\hspace{0.01in}|q\hspace{0.01in}%
(\pm\hspace{0.01in}q^{-n})\hspace{0.01in}x).
\end{align}
This result, together with Eq.~(\ref{IntExpExpInv}), implies:%
\begin{equation}
(\mathcal{F}_{q}^{\hspace{0.01in}-1}\circ\mathcal{F}_{q})(\phi_{n}^{\pm
}(x))=\left\{
\begin{array}
[c]{cc}%
1 & \text{for\quad}x=\pm\hspace{0.01in}x_{0}\hspace{0.01in}q^{n},\\
0 & \text{for\quad}x\neq\pm\hspace{0.01in}x_{0}\hspace{0.01in}q^{n}.
\end{array}
\right.  \label{InvFourIdeEinQ1}%
\end{equation}
We can verify the second identity of Eq.~(\ref{IdeInvFouTraBasFktqEin}) in the
same manner:%
\begin{align}
&  (\mathcal{F}_{q}\circ\mathcal{F}_{q}^{\hspace{0.01in}-1})(\psi_{n}^{\pm
}(\hspace{0.01in}p))=\nonumber\\
&  \qquad=\frac{1}{4\Theta_{0}}\int_{-x_{0}.\infty}^{\hspace{0.01in}%
x_{0}.\infty}\text{d}_{q}x\,\int_{-p_{0}.\infty}^{\hspace{0.01in}p_{0}.\infty
}\text{d}_{q}p^{\prime}\,\psi_{n}^{\pm}(\hspace{0.01in}p^{\prime}%
)\hspace{0.01in}\exp_{q^{-1}}(-\text{i}p^{\prime}|qx)\hspace{0.01in}\exp
_{q}(x|\text{i}p)\nonumber\\
&  \qquad=\frac{(1-q)\hspace{0.01in}p_{0}\hspace{0.01in}q^{n}}{4\Theta_{0}%
}\int_{-x_{0}.\infty}^{\hspace{0.01in}x_{0}.\infty}\text{d}_{q}x\,\exp
_{q^{-1}}(-\text{i}\hspace{0.01in}(\pm\hspace{0.01in}q^{n})\hspace
{0.01in}p_{0}|qx)\hspace{0.01in}\exp_{q}(x|\text{i}p)\nonumber\\
&  \qquad=\frac{(1-q)\hspace{0.01in}p_{0}}{4\Theta_{0}}\int_{-x_{0}.\infty
}^{\hspace{0.01in}x_{0}.\infty}\text{d}_{q}x\,\exp_{q^{-1}}(-\text{i}%
p_{0}|qx)\hspace{0.01in}\exp_{q}(x|\hspace{-0.02in}\pm\hspace{-0.02in}%
\text{i}q^{-n}p).
\end{align}
Taking into account Eq.~(\ref{IntExpExpInv2}), we end up with:%
\begin{equation}
(\mathcal{F}_{q}\circ\mathcal{F}_{q}^{\hspace{0.01in}-1})(\psi_{n}^{\pm
})(\hspace{0.01in}p)=\left\{
\begin{array}
[c]{cc}%
1 & \text{for\quad}p=\pm\hspace{0.01in}p_{0}\hspace{0.01in}q^{n},\\
0 & \text{for\quad}p\neq\pm\hspace{0.01in}p_{0}\hspace{0.01in}q^{n}.
\end{array}
\right.  \label{InvFourIdeEinQ2}%
\end{equation}

For $x\in\mathbb{G}_{q,\hspace{0.01in}x_{0}}$ and $p\in\mathbb{G}%
_{q,\hspace{0.01in}p_{0}}$, Eqs.~(\ref{IntExpExpInv}) and (\ref{IntExpExpInv2}%
) become ($n,m\in\mathbb{Z}$ and $\varepsilon,\varepsilon^{\hspace
{0.01in}\prime}=\pm1$)%
\begin{align}
&  \int_{-p_{0}.\infty}^{\hspace{0.01in}p_{0}.\infty}\text{d}_{q}p\,\exp
_{q}(\varepsilon x_{0}\hspace{0.01in}q^{n}|\text{i}p)\hspace{0.01in}%
\exp_{q^{-1}}(-\text{i}p|q\varepsilon^{\hspace{0.01in}\prime}x_{0}%
\hspace{0.01in}q^{m})=\nonumber\\
&  \qquad\qquad=q^{-n}\int_{-p_{0}.\infty}^{\hspace{0.01in}p_{0}.\infty
}\text{d}_{q}p\,\exp_{q}(x_{0}|\text{i}p)\hspace{0.01in}\exp_{q^{-1}%
}(-\text{i}p|q\varepsilon\varepsilon^{\hspace{0.01in}\prime}x_{0}%
\hspace{0.01in}q^{m-n})\nonumber\\
&  \qquad\qquad=\frac{4\Theta_{0}}{(1-q)\hspace{0.01in}x_{0}\hspace
{0.01in}q^{n}}\,\delta_{nm}\hspace{0.01in}\delta_{\varepsilon\varepsilon
^{\prime}} \label{PreVolRelExp1Dim}%
\end{align}
or%
\begin{align}
&  \int_{-x_{0}.\infty}^{\hspace{0.01in}x_{0}.\infty}\text{d}_{q}%
x\,\exp_{q^{-1}}(-\text{i}\varepsilon p_{0}\hspace{0.01in}q^{n}|qx)\hspace
{0.01in}\exp_{q}(x|\text{i}\varepsilon^{\hspace{0.01in}\prime}p_{0}%
\hspace{0.01in}q^{m})=\nonumber\\
&  \qquad\qquad=q^{-n}\int_{-x_{0}.\infty}^{\hspace{0.01in}x_{0}.\infty
}\text{d}_{q}x\,\exp_{q^{-1}}(-\text{i}p_{0}|qx)\hspace{0.01in}\exp
_{q}(x|\text{i}\varepsilon\varepsilon^{\hspace{0.01in}\prime}p_{0}%
\hspace{0.01in}q^{m-n})\nonumber\\
&  \qquad\qquad=\frac{4\Theta_{0}}{(1-q)\hspace{0.01in}p_{0}\hspace
{0.01in}q^{n}}\,\delta_{nm}\hspace{0.01in}\delta_{\varepsilon\varepsilon
^{\prime}}. \label{PreOrtRelExp1Dim}%
\end{align}
We get the second integral expression in the equations above by rescaling the
integration variable. Next, we rewrite the above identities using the
following functions:%
\begin{equation}
u_{p}(x)=\frac{1}{2\sqrt{\Theta_{0}}}\exp_{q}(x|\text{i}p),\qquad u_{p}^{\ast
}(x)=\frac{1}{2\sqrt{\Theta_{0}}}\exp_{q^{-1}}(-\text{i}p|qx).
\label{ImpEigQDefEin}%
\end{equation}
The \textit{completeness relation} in Eq.~(\ref{PreVolRelExp1Dim}) yields%
\begin{align}
&  \int_{-p_{0}.\infty}^{\hspace{0.01in}p_{0}.\infty}\text{d}_{q}%
p\,u_{p}(\varepsilon x_{0}\hspace{0.01in}q^{n})\,u_{p}^{\ast}(\varepsilon
^{\hspace{0.01in}\prime}x_{0}\hspace{0.01in}q^{m})=\frac{\delta_{nm}%
\hspace{0.01in}\delta_{\varepsilon\varepsilon^{\prime}}}{(1-q)\hspace
{0.01in}x_{0}\hspace{0.01in}q^{n}}\nonumber\\
&  \qquad\qquad=\left.  \delta_{q}(x\,\bar{\oplus}\,(\bar{\ominus}%
\,q\hspace{0.01in}y))\right\vert _{x\hspace{0.01in}=\hspace{0.01in}\varepsilon
x_{0}q^{n};\hspace{0.01in}y\hspace{0.01in}=\hspace{0.01in}\varepsilon^{\prime
}x_{0}q^{m}}. \label{VolRelExp1Dim}%
\end{align}
Similarly, the \textit{orthonormality relation} in Eq.~(\ref{PreOrtRelExp1Dim}%
) takes on the following form:%
\begin{align}
&  \int_{-x_{0}.\infty}^{\hspace{0.01in}x_{0}.\infty}\text{d}_{q}%
x\,u_{\varepsilon p_{0}q^{n}}^{\ast}(x)\,u_{\varepsilon^{\prime}p_{0}q^{m}%
}(x)=\frac{\delta_{nm}\hspace{0.01in}\delta_{\varepsilon\varepsilon^{\prime}}%
}{(1-q)\hspace{0.01in}p_{0}\hspace{0.01in}q^{n}}\nonumber\\
&  \qquad\qquad=\left.  \delta_{q}((\bar{\ominus}\,qp)\,\bar{\oplus
}\,p^{\prime})\right\vert _{p\hspace{0.01in}=\hspace{0.01in}\varepsilon
p_{0}q^{n};\hspace{0.01in}p^{\prime}=\hspace{0.01in}\varepsilon^{\prime}%
p_{0}q^{m}}. \label{ONRel1DimImpFkt}%
\end{align}
In writing down the above relations, we used the expressions for the $q$-delta
functions in Eq.~(\ref{DarQDelKroDel1Dim}) of the previous chapter and took
into account the following identities [cf. Eq.~(\ref{BasVekFktBraLin}) of the
previous chapter]:%
\begin{equation}
\phi_{m}^{\hspace{0.01in}\varepsilon}(\varepsilon^{\hspace{0.01in}\prime}%
x_{0}\hspace{0.01in}q^{n})=\delta_{mn}\hspace{0.01in}\delta_{\varepsilon
\varepsilon^{\prime}},\qquad\psi_{m}^{\hspace{0.01in}\varepsilon}%
(\varepsilon^{\hspace{0.01in}\prime}p_{0}\hspace{0.01in}q^{n})=\delta
_{mn}\hspace{0.01in}\delta_{\varepsilon\varepsilon^{\prime}}.
\end{equation}

Due to the completeness relation in Eq.~(\ref{VolRelExp1Dim}), we can expand a
function $f(x)$ in position space $\mathcal{M}(x)$ in terms of the functions
$u_{p}^{\ast}(x)$:%
\begin{equation}
f(x)=\int_{-p_{0}.\infty}^{\hspace{0.01in}p_{0}.\infty}\text{d}_{q}%
p\,a_{p}\hspace{0.01in}u_{p}^{\ast}(x)\quad\text{with}\quad a_{p}%
=\mathcal{F}_{q}(f)(\hspace{0.01in}p). \label{FouDarFkt1DimQ}%
\end{equation}
Similarly, we can expand a function $g(\hspace{0.01in}p)$ in position space
$\mathcal{M}(\hspace{0.01in}p)$ in terms of the functions $u_{p}(x)$:%
\begin{equation}
g(\hspace{0.01in}p)=\int_{-x_{0}.\infty}^{\hspace{0.01in}x_{0}.\infty}%
\text{d}_{q}x\,b_{x}\hspace{0.01in}u_{p}(x)\quad\text{with}\quad
b_{x}=\mathcal{F}_{q}^{\hspace{0.01in}-1}(\hspace{0.01in}g)(x).
\end{equation}

Next, we consider a scaling operator that acts on functions in position or
momentum space as follows:%
\begin{equation}
\Lambda\hspace{0.01in}f(x)=f(qx),\qquad\Lambda\hspace{0.01in}g(\hspace
{0.01in}p)=g(q^{-1}p),
\end{equation}
We show that the scaling operator $\Lambda$ commutes with $q$-Fourier
transforms as follows \cite{olshanetsky1998q.alg}:%
\begin{equation}
\mathcal{F}_{q}\hspace{0.01in}\Lambda=q^{-1}\Lambda\hspace{0.01in}%
\mathcal{F}_{q},\qquad\mathcal{F}_{q}^{\hspace{0.01in}-1}\Lambda
=q\hspace{0.01in}\Lambda\hspace{0.01in}\mathcal{F}_{q}^{\hspace{0.01in}-1}.
\label{VerSkaFou1DimQ}%
\end{equation}
We prove the first relation by the following calculation:%
\begin{align}
\mathcal{F}_{q}(\Lambda\phi(x))(\hspace{0.01in}p)  &  =\mathcal{F}_{q}%
(\phi(qx))(\hspace{0.01in}p)\nonumber\\
&  =\frac{1}{2\sqrt{\Theta_{0}}}\int_{-x_{0}.\infty}^{\hspace{0.01in}%
x_{0}.\infty}\text{d}_{q}x\,\phi(qx)\hspace{0.01in}\exp_{q}(x|\text{i}%
p)\nonumber\\
&  =\frac{q^{-1}}{2\sqrt{\Theta_{0}}}\int_{-x_{0}.\infty}^{\hspace
{0.01in}x_{0}.\infty}\text{d}_{q}x\,\phi(x)\hspace{0.01in}\exp_{q}%
(x|\text{i}q^{-1}p)\nonumber\\
&  =q^{-1}(\mathcal{F}_{q}\hspace{0.01in}\phi)(q^{-1}p)=q^{-1}\Lambda
(\mathcal{F}_{q}\hspace{0.01in}\phi)(\hspace{0.01in}p).
\end{align}
Similar considerations lead to the second relation in
Eq.~(\ref{VerSkaFou1DimQ}):%
\begin{align}
\mathcal{F}_{q}^{\hspace{0.01in}-1}(\Lambda\psi(\hspace{0.01in}p))(x)  &
=\mathcal{F}_{q}^{\hspace{0.01in}-1}(\psi(q^{-1}p))(x)\nonumber\\
&  =\frac{1}{2\sqrt{\Theta_{0}}}\int_{-p_{0}.\infty}^{\hspace{0.01in}%
p_{0}.\infty}\text{d}_{q}p\,\psi(q^{-1}p)\hspace{0.01in}\exp_{q^{-1}%
}(-\text{i}p|qx)\nonumber\\
&  =\frac{q}{2\sqrt{\Theta_{0}}}\int_{-p_{0}.\infty}^{\hspace{0.01in}%
p_{0}.\infty}\text{d}_{q}p\,\psi(\hspace{0.01in}p)\hspace{0.01in}\exp_{q^{-1}%
}(-\text{i}p|qqx)\nonumber\\
&  =q\hspace{0.01in}\mathcal{F}_{q}^{\hspace{0.01in}-1}(\psi)(qx)=q\hspace
{0.01in}\Lambda\hspace{0.01in}\mathcal{F}_{q}^{\hspace{0.01in}-1}(\psi)(x).
\end{align}

As is well known, ordinary Fourier transforms turn derivative operators and
multiplication operators into each other. The $q$-de\-formed Fourier
transforms show this property, as well \cite{olshanetsky1998q.alg}, since we
have%
\begin{equation}
\mathcal{F}_{q}(D_{q,\hspace{0.01in}x}\hspace{0.01in}\phi)=-\text{i}%
\Lambda\hspace{0.01in}p\hspace{0.01in}\mathcal{F}_{q}(\phi),\qquad
\mathcal{F}_{q}(x\phi)=-\text{i}D_{q,\hspace{0.01in}p}\hspace{0.01in}%
\mathcal{F}_{q}(\phi)\label{VerFouQ1DimAblOrt1}%
\end{equation}
and%
\begin{equation}
\mathcal{F}_{q}^{\hspace{0.01in}-1}(D_{q,\hspace{0.01in}p}\hspace{0.01in}%
\psi)=\text{i}\hspace{0.01in}\mathcal{F}_{q}^{\hspace{0.01in}-1}%
(\psi)\,x,\qquad\mathcal{F}_{q}^{\hspace{0.01in}-1}(\hspace{0.01in}%
p\hspace{0.01in}\psi)=\text{i}\hspace{0.01in}q^{-1}\Lambda^{-1}D_{q,\hspace
{0.01in}x}\hspace{0.01in}\mathcal{F}_{q}^{\hspace{0.01in}-1}(\psi
).\label{VerFouQ1DimAblOrt2}%
\end{equation}
We only prove Eq.~(\ref{VerFouQ1DimAblOrt2}) since
Eq.~(\ref{VerFouQ1DimAblOrt1}) can be proven by similar reasoning. The first
identity in Eq.~(\ref{VerFouQ1DimAblOrt2}) can be verified as follows:%
\begin{align}
\mathcal{F}_{q}^{\hspace{0.01in}-1}(D_{q,\hspace{0.01in}p}\hspace{0.01in}\psi)
&  =\frac{1}{2\sqrt{\Theta_{0}}}\int_{-p_{0}.\infty}^{\hspace{0.01in}%
p_{0}.\infty}\text{d}_{q}p\,[D_{q,\hspace{0.01in}p}\psi(\hspace{0.01in}%
p)]\exp_{q^{-1}}(-\text{i}p|qx)\nonumber\\
&  =-\frac{1}{2\sqrt{\Theta_{0}}}\int_{-p_{0}.\infty}^{\hspace{0.01in}%
p_{0}.\infty}\text{d}_{q}p\,\psi(qp)\,D_{q,\hspace{0.01in}p}\exp_{q^{-1}%
}(-\text{i}p|qx)\nonumber\\
&  =\frac{\text{i}q}{2\sqrt{\Theta_{0}}}\int_{-p_{0}.\infty}^{\hspace
{0.01in}p_{0}.\infty}\text{d}_{q}p\,\psi(qp)\hspace{0.01in}\exp_{q^{-1}%
}(-\text{i}qp|qx)\,x\nonumber\\
&  =\frac{\text{i}}{2\sqrt{\Theta_{0}}}\int_{-p_{0}.\infty}^{\hspace
{0.01in}p_{0}.\infty}\text{d}_{q}p\,\psi(\hspace{0.01in}p)\hspace{0.01in}%
\exp_{q^{-1}}(-\text{i}p|qx)\,x\nonumber\\
&  =\text{i}\hspace{0.01in}\mathcal{F}_{q}^{\hspace{0.01in}-1}(\psi)\,x.
\end{align}
To check the second identity in Eq.~(\ref{VerFouQ1DimAblOrt2}), we proceed as
follows:%
\begin{align}
\mathcal{F}_{q}^{\hspace{0.01in}-1}(\hspace{0.01in}p\hspace{0.01in}\psi) &
=\frac{1}{2\sqrt{\Theta_{0}}}\int_{-p_{0}.\infty}^{\hspace{0.01in}p_{0}%
.\infty}\text{d}_{q}p\,p\hspace{0.01in}\psi(\hspace{0.01in}p)\hspace
{0.01in}\exp_{q^{-1}}(-\text{i}p|qx)\nonumber\\
&  =\frac{\text{i}}{2\sqrt{\Theta_{0}}}\int_{-p_{0}.\infty}^{\hspace
{0.01in}p_{0}.\infty}\text{d}_{q}p\,\psi(\hspace{0.01in}p)\,D_{q,x}%
\exp_{q^{-1}}(-\text{i}p|x)\nonumber\\
&  =\text{i}D_{q,\hspace{0.01in}x}\hspace{0.01in}\mathcal{F}_{q}%
^{\hspace{0.01in}-1}(\psi)(q^{-1}x)=\text{i}q^{-1}\Lambda^{-1}D_{q,\hspace
{0.01in}x}\hspace{0.01in}\mathcal{F}_{q}^{\hspace{0.01in}-1}(\psi)(x).
\end{align}
In the two computations above, we have applied Eq.~(\ref{AblJacImpInvExp}).

The ordinary Fourier transform is a unitary operator. Thus, the scalar product
of two given functions has the same value as the scalar product of their
Fourier transforms. In this respect, the $q$-de\-formed Fourier transforms in
Eq.~(\ref{DefQFouTraEin}) are unitary maps on the Hilbert space $L_{q}^{2}$
with the scalar product in Eq.~(\ref{SkaProqDefEinBraLin}) of the
previous\ chapter. However, the $q$-de\-formed Fourier transforms of the two
factors of the $q$-scalar product slightly differ from each other. Concretely,
we have%
\begin{equation}
\left\langle \hspace{0.01in}g\hspace{0.01in},f\right\rangle _{q,x}%
=\left\langle \psi\hspace{0.01in},\phi\right\rangle _{q,p}
\label{ParIdeFour1dimQ}%
\end{equation}
with%
\begin{align}
\phi(\hspace{0.01in}p)  &  =\mathcal{F}_{q}(f)(\hspace{0.01in}p)=\frac
{1}{2\sqrt{\Theta_{0}}}\int_{-x_{0}.\infty}^{\hspace{0.01in}x_{0}.\infty
}\text{d}_{q}x\,\exp_{q}(\text{i}p|x)\hspace{0.01in}f(x),\nonumber\\[0.04in]
\psi(\hspace{0.01in}p)  &  =\mathcal{\tilde{F}}_{q}(\hspace{0.01in}%
g)(\hspace{0.01in}p)=\frac{1}{2\sqrt{\Theta_{0}}}\int_{-x_{0}.\infty}%
^{\hspace{0.01in}x_{0}.\infty}\text{d}_{q}x\,\exp_{q^{-1}}(\text{i}%
p|qx)\,g(x),
\end{align}
and%
\begin{align}
f(x)  &  =\mathcal{F}_{q}^{\hspace{0.01in}-1}(\phi)(x)=\frac{1}{2\sqrt
{\Theta_{0}}}\int_{-p_{0}.\infty}^{\hspace{0.01in}p_{0}.\infty}\text{d}%
_{q}p\,\phi(\hspace{0.01in}p)\exp_{q^{-1}}(-\text{i}p|qx),\nonumber\\[0.04in]
g(x)  &  =\mathcal{\tilde{F}}_{q}^{\hspace{0.01in}-1}(\psi)(x)=\frac{1}%
{2\sqrt{\Theta_{0}}}\int_{-p_{0}.\infty}^{\hspace{0.01in}p_{0}.\infty}%
\text{d}_{q}p\,\exp_{q}(x|\hspace{-0.02in}-\hspace{-0.02in}\text{i}%
p)\,\psi(\hspace{0.01in}p). \label{InvFourTraInd}%
\end{align}
To prove Eq.~(\ref{ParIdeFour1dimQ}), we insert the expressions for the
$q$-Fourier transforms into one of the two scalar products and apply
Eq.~(\ref{PreVolRelExp1Dim}) or Eq.~(\ref{PreOrtRelExp1Dim}).

Our aim is to consider the $q$\textbf{-Fourier transforms of elementary
functions and }$q$\textbf{-dis\-tri\-bu\-tions}. First, we calculate the
$q$\textbf{-Fourier transform of }$x^{-1}$:%
\begin{align}
\mathcal{F}_{q}(x^{-1})(\hspace{0.01in}p) &  =\frac{1}{2\sqrt{\Theta_{0}}}%
\int_{-x_{0}.\infty}^{\hspace{0.01in}x_{0}.\infty}\text{d}_{q}x\,x^{-1}%
\exp_{q}(x|\text{i}p)=\nonumber\\
&  =\frac{1-q}{2\sqrt{\Theta_{0}}}\sum_{\varepsilon\hspace{0.01in}%
=\hspace{0.01in}\pm}\,\sum_{m\hspace{0.01in}=\hspace{0.01in}-\infty}^{\infty
}\varepsilon x_{0}\hspace{0.01in}q^{m}x_{0}^{-1}q^{-m}\exp_{q}(\varepsilon
x_{0}\hspace{0.01in}q^{m}|\text{i}p)\nonumber\\
&  =\frac{1-q}{2\sqrt{\Theta_{0}}}\sum_{m\hspace{0.01in}=\hspace
{0.01in}-\infty}^{\infty}\left[  \exp_{q}(\text{i}x_{0}\hspace{0.01in}%
q^{m}p)-\exp_{q}(-\text{i}x_{0}\hspace{0.01in}q^{m}p)\right]  \nonumber\\
&  =\frac{\text{i}(1-q)}{\sqrt{\Theta_{0}}}\sum_{m\hspace{0.01in}%
=\hspace{0.01in}-\infty}^{\infty}\sin_{q}(x_{0}\hspace{0.01in}q^{m}%
p)=\frac{\text{i}}{\sqrt{\Theta_{0}}}\,\Theta_{q}(x_{0}\hspace{0.01in}%
p).\label{FouTrazInvEinQ1}%
\end{align}
For the final step in the above computation, we applied Eq.~(\ref{DefTheZ}) of
Chap.~\ref{KapQIntTrig}. We can simplify the final expression in
Eq.~(\ref{FouTrazInvEinQ1}) since $p$ can only take the values $\pm
\hspace{0.01in}p_{0}\hspace{0.01in}q^{m}$ with $m\in\mathbb{Z}$:%
\begin{align}
\frac{\text{i}}{\sqrt{\Theta_{0}}}\,\Theta_{q}(\pm\hspace{0.01in}x_{0}%
\hspace{0.01in}p_{0}\hspace{0.01in}q^{m}) &  =\pm\frac{\text{i}}{\sqrt
{\Theta_{0}}}\,\Theta_{q}(x_{0}\hspace{0.01in}p_{0}\hspace{0.01in}q^{m}%
)=\pm\frac{\text{i}}{\sqrt{\Theta_{0}}}\,\Theta_{q}(x_{0}\hspace{0.01in}%
p_{0})\nonumber\\
&  =\pm\frac{\text{i}}{\sqrt{\Theta_{0}}}\,\Theta_{0}=\pm\hspace
{0.01in}\text{i}\sqrt{\Theta_{0}}.\label{UmfFouTrafInv1}%
\end{align}
In the first step, we took into account that $\Theta_{q}$ is an odd function
due to its definition [cf. Eq.~(\ref{DefTheZ}) of Chap.~\ref{KapQIntTrig}]:%
\begin{equation}
\Theta_{q}(-z)=-\hspace{0.01in}\Theta_{q}(z).
\end{equation}
In the second step, we used the fact that $\Theta_{q}$ is also a
$q$-pe\-ri\-od\-ic function [cf.\ Eq.~(\ref{PerTheFkt}) of
Chap.~\ref{KapQIntTrig}]. The penultimate step follows from
Eq.~(\ref{DefTheNul}). Combining the results of Eqs.~(\ref{FouTrazInvEinQ1})
and (\ref{UmfFouTrafInv1}), we finally obtain \cite{olshanetsky1998q.alg}:%
\begin{equation}
\mathcal{F}_{q}(x^{-1})(\hspace{0.01in}p)=\text{i}\sqrt{\Theta_{0}%
}\operatorname*{sgn}(\hspace{0.01in}p)=\text{i}\sqrt{\Theta_{0}}\left[
\theta_{q}^{+}(\hspace{0.01in}p)-\theta_{q}^{-}(\hspace{0.01in}p)\right]
.\label{FouTraInvXEin}%
\end{equation}

Next, we calculate the $q$\textbf{-Fourier transform of }$1$%
.\ Eq.~(\ref{FouTraInvXEin})\ together with Eq.~(\ref{AblTheFktEin}) of the
previous chapter and Eq.~(\ref{VerFouQ1DimAblOrt1}) implies
\cite{Kempf:1994yd,olshanetsky1998q.alg}:%
\begin{align}
\mathcal{F}_{q}(1)(\hspace{0.01in}p)  &  =\mathcal{F}_{q}(xx^{-1}%
)(\hspace{0.01in}p)=-\text{i}D_{q,\hspace{0.01in}p}\hspace{0.01in}%
\mathcal{F}_{q}(x^{-1})(\hspace{0.01in}p)\nonumber\\
&  =\sqrt{\Theta_{0}}\,D_{q,\hspace{0.01in}p}\left[  \theta_{q}^{+}%
(\hspace{0.01in}p)-\theta_{q}^{-}(\hspace{0.01in}p)\right]  =2\sqrt{\Theta
_{0}}\,\delta_{q}(\hspace{0.01in}p). \label{FouTra1EinQ}%
\end{align}
Using this result and Eq.~(\ref{JacInvYSpi}) of Chap.~\ref{KapQIntTrig}, we
can show that $\delta_{q}(\hspace{0.01in}p)$ is even:%
\begin{align}
\delta_{q}(-p)  &  =\frac{1}{2\sqrt{\Theta_{0}}}\,\mathcal{F}_{q}%
(1)(-p)=\frac{1}{4\Theta_{0}}\int_{-x_{0}.\infty}^{\hspace{0.01in}x_{0}%
.\infty}\text{d}_{q}x\,\exp_{q}(x|\hspace{-0.02in}-\hspace{-0.02in}%
\text{i}p)\nonumber\\
&  =\frac{1}{4\Theta_{0}}\int_{-x_{0}.\infty}^{\hspace{0.01in}x_{0}.\infty
}\text{d}_{q}x\,\exp_{q}(\text{i}p|\hspace{-0.02in}-\hspace{-0.02in}%
x)=\nonumber\\
&  =\frac{1}{4\Theta_{0}}\int_{-x_{0}.\infty}^{\hspace{0.01in}x_{0}.\infty
}\text{d}_{q}x\,\exp_{q}(\text{i}p|x)=\delta_{q}(\hspace{0.01in}p).
\label{GerDelFkt1Dim}%
\end{align}
Eq.~(\ref{FouTra1EinQ}) implies that the shifted $q$-delta function can be
written as the $q$\textbf{-Fourier transform of a }$q$%
\textbf{-ex\-po\-nen\-tial function}:%
\begin{gather}
\delta_{q}(a\,\bar{\oplus}\,x)=\frac{1}{2\sqrt{\Theta_{0}}}\,\mathcal{F}%
_{q}(1)(a\,\bar{\oplus}\,x)=\frac{1}{4\Theta_{0}}\int_{-p_{0}.\infty}%
^{\hspace{0.01in}p_{0}.\infty}\text{d}_{q}p\,\exp_{q}(\text{i}p|a\,\bar
{\oplus}\,x)\nonumber\\
=\frac{1}{4\Theta_{0}}\int_{-p_{0}.\infty}^{\hspace{0.01in}p_{0}.\infty
}\text{d}_{q}p\,\exp_{q}(\text{i}p|a)\hspace{0.01in}\exp_{q}(\text{i}%
p|x)=\frac{1}{2\sqrt{\Theta_{0}}}\,\mathcal{F}_{q}(\exp_{q}(\text{i}p|a))(x).
\end{gather}
In the second last step of the above computation, we applied the addition
theorem for $q$-ex\-po\-nen\-tials [cf. Eq.~(\ref{AddTheqExp1Dim}) of
Chap.~\ref{KapqAnaExpTriFkt}].

Next, we compute the $q$-\textbf{Fourier transform of }$\theta_{q}^{+}%
-\hspace{0.01in}\theta_{q}^{-}$. Using Eq.~(\ref{DefTheQDefEin}) of the
previous chapter and Eq.~(\ref{DefqTriSinCosExp1}) of Chap.~\ref{KapQTrigFkt},
we get \cite{olshanetsky1998q.alg}:%
\begin{align}
\mathcal{F}_{q}(\theta_{q}^{+}-\hspace{0.01in}\theta_{q}^{-})(\hspace
{0.01in}p)  &  =\frac{1}{2\sqrt{\Theta_{0}}}\int\limits_{0}^{x_{0}.\infty
}\text{d}_{q}x\,\exp_{q}(\text{i}x|p)-\frac{1}{2\sqrt{\Theta_{0}}}%
\int\limits_{-x_{0}.\infty}^{0}\text{d}_{q}x\,\exp_{q}(\text{i}x|p)\nonumber\\
&  =\frac{(1-q)\hspace{0.01in}x_{0}}{2\sqrt{\Theta_{0}}}\sum_{m\hspace
{0.01in}=\hspace{0.01in}-\infty}^{\infty}q^{m}\left[  \exp_{q}(\text{i}%
x_{0}\hspace{0.01in}q^{m}p)-\exp_{q}(-\text{i}x_{0}\hspace{0.01in}%
q^{m}p)\right] \nonumber\\
&  =\frac{\text{i}(1-q)\hspace{0.01in}x_{0}}{\sqrt{\Theta_{0}}}\sum
_{m\hspace{0.01in}=\hspace{0.01in}-\infty}^{\infty}q^{m}\sin_{q}(x_{0}%
\hspace{0.01in}q^{m}p)\nonumber\\
&  =\frac{\text{i}}{\sqrt{\Theta_{0}}}\,p^{-1}\lim_{M\rightarrow
\hspace{0.01in}\infty}\,(1-q)\hspace{0.01in}x_{0}\hspace{0.01in}%
p\sum_{m\hspace{0.01in}=\hspace{0.01in}-M}^{\infty}q^{m}\sin_{q}(x_{0}%
\hspace{0.01in}q^{m}p)\nonumber\\
&  =\frac{\text{i}}{\sqrt{\Theta_{0}}}\,p^{-1}\lim_{M\rightarrow
\hspace{0.01in}\infty}\,\left[  1-\cos_{q}(x_{0}\hspace{0.01in}q^{-M}%
p)\right]  =\frac{\text{i}}{\sqrt{\Theta_{0}}}\,p^{-1}.
\label{FouTraThePluMinEinQ}%
\end{align}
The above computation follows the same logic as the one in
Eq.~(\ref{IntSinQFktEin}) of Chap.~\ref{KapQTrigFkt}.

From the $q$-Fourier transform of $\theta_{q}^{+}-\hspace{0.01in}\theta
_{q}^{-}$, we derive the $q$\textbf{-Fourier transform of the }$q$%
\textbf{-delta function}. Using Eq.~(\ref{AblTheFktEin}) of the previous
chapter and Eq.~(\ref{VerFouQ1DimAblOrt1}), we get \cite{olshanetsky1998q.alg}%
:%
\begin{align}
\mathcal{F}_{q}(\delta_{q})(\hspace{0.01in}p)  &  =\frac{1}{2}\mathcal{F}%
_{q}\left(  D_{q}(\theta_{q}^{+}-\hspace{0.01in}\theta_{q}^{-})\right)
(\hspace{0.01in}p)\nonumber\\
&  =-\frac{\text{i}}{2}\hspace{0.01in}\Lambda\hspace{0.01in}p\hspace
{0.01in}\mathcal{F}_{q}(\theta_{q}^{+}-\hspace{0.01in}\theta_{q}^{-}%
)(\hspace{0.01in}p)=\frac{1}{2\sqrt{\Theta_{0}}}. \label{FouTraQDelEinDim}%
\end{align}
A direct computation using Eq.~(\ref{DefDelQEinDim}) of the previous chapter
and Eq.~(\ref{Bed2qExp}) of Chap.~\ref{KapqAnaExpTriFkt} gives the same
result:%
\begin{align}
\mathcal{F}_{q}(\delta_{q})(\hspace{0.01in}p)  &  =\frac{1}{2\sqrt{\Theta_{0}%
}}\int_{-x_{0}.\infty}^{\hspace{0.01in}x_{0}.\infty}\text{d}_{q}x\,\delta
_{q}(x)\hspace{0.01in}\exp_{q}(x|\text{i}p)\nonumber\\
&  =\frac{1}{2\sqrt{\Theta_{0}}}\exp_{q}(x|\text{i}p)|_{x\hspace
{0.01in}=\hspace{0.01in}0}=\frac{1}{2\sqrt{\Theta_{0}}}.
\end{align}

We can use the results from Eqs.~(\ref{FouTra1EinQ}) and
(\ref{FouTraThePluMinEinQ}) to find the $q$\textbf{-Fourier transforms of
}$\theta_{q}^{+}$\textbf{ and }$\theta_{q}^{-}$. Due to $1=\theta_{q}%
^{+}+\hspace{0.01in}\theta_{q}^{-}$, we have \cite{olshanetsky1998q.alg}:%
\begin{align}
\mathcal{F}_{q}(\theta_{q}^{+})  &  =\frac{1}{2}\mathcal{F}_{q}(\theta_{q}%
^{+}-\hspace{0.01in}\theta_{q}^{-}+1)=\frac{\text{i}}{2\sqrt{\Theta_{0}}%
}\hspace{0.01in}p^{-1}+\sqrt{\Theta_{0}}\,\delta_{q}(\hspace{0.01in}%
p),\nonumber\\
\mathcal{F}_{q}(\theta_{q}^{-})  &  =\frac{1}{2}\mathcal{F}_{q}(\theta_{q}%
^{-}-\hspace{0.01in}\theta_{q}^{+}+1)=\frac{-\text{i}}{2\sqrt{\Theta_{0}}%
}\hspace{0.01in}p^{-1}+\sqrt{\Theta_{0}}\,\delta_{q}(\hspace{0.01in}p).
\end{align}

Furthermore, we can calculate the $q$\textbf{-Fourier transform of a power
function with an integer exponent}. As we show below, it holds
\cite{olshanetsky1998q.alg} ($n\in\mathbb{N}_{0}$)%
\begin{equation}
\mathcal{F}_{q}(x^{n})=\text{i}^{n}\hspace{0.01in}2\sqrt{\Theta_{0}%
}\,q^{-n(n\hspace{0.01in}+1)/2}[[\hspace{0.01in}n]]_{q}!\,p^{-n}\delta
_{q}(\hspace{0.01in}p) \label{FouTraGanPotEindQ}%
\end{equation}
and%
\begin{equation}
\mathcal{F}_{q}(x^{-n\hspace{0.01in}-1})=\text{i}^{n\hspace{0.01in}+1}%
\sqrt{\Theta_{0}}\,([[\hspace{0.01in}n]]_{q}!)^{-1}p^{n}\operatorname*{sgn}%
(\hspace{0.01in}p). \label{FouTraGanPotEindQ2}%
\end{equation}
From Eq.~(\ref{FouTraGanPotEindQ}), we obtain a formula for calculating the
$q$-Fourier transform of a function that can be represented as a power series:%
\begin{equation}
\mathcal{F}_{q}(\hspace{0.01in}g)=2\sqrt{\Theta_{0}}\sum_{n\hspace
{0.01in}=\hspace{0.01in}0}^{\infty}\text{i}^{n}q^{-n(n\hspace{0.01in}%
+1)/2}p^{-n}\delta_{q}(\hspace{0.01in}p)\hspace{0.01in}(D_{q}^{n}%
g)(0)\hspace{0.02in}. \label{QFouTraAlgFktEin}%
\end{equation}

To prove Eq.~(\ref{FouTraGanPotEindQ}), we do the following calculation:%
\begin{align}
\delta_{q}(\hspace{0.01in}p)  &  =\frac{1}{2\sqrt{\Theta_{0}}}\hspace
{0.01in}\mathcal{F}_{q}(1)(\hspace{0.01in}p)=\frac{1}{2\sqrt{\Theta_{0}}%
}\hspace{0.01in}\frac{\mathcal{F}_{q}(D_{q}^{n}x^{n})(\hspace{0.01in}%
p)}{[[\hspace{0.01in}n]]_{q}!}\nonumber\\
&  =\frac{1}{2\sqrt{\Theta_{0}}}\hspace{0.01in}\frac{(-\text{i})^{n}%
}{[[\hspace{0.01in}n]]_{q}!}\hspace{0.01in}(\Lambda\hspace{0.01in}%
p)^{n}\hspace{0.01in}\mathcal{F}_{q}(x^{n})(\hspace{0.01in}p).
\label{SolfForFou}%
\end{align}
The first identity follows from Eq.~(\ref{FouTra1EinQ}), the second one from
the identity%
\begin{equation}
1=([[\hspace{0.01in}n]]_{q}!)^{-1}D_{q}^{n}x^{n}, \label{DarIdeJaclKoo}%
\end{equation}
and the last one from Eq.~(\ref{VerFouQ1DimAblOrt1}). Solving
Eq.~(\ref{SolfForFou}) for $\mathcal{F}_{q}(x^{n})$ finally gives:%
\begin{align}
\mathcal{F}_{q}(x^{n})(\hspace{0.01in}p)  &  =2\sqrt{\Theta_{0}}%
\hspace{0.01in}\text{i}^{n}[[\hspace{0.01in}n]]_{q}!\hspace{0.01in}%
(\hspace{0.01in}p^{-1}\Lambda\hspace{-0.02in}^{-1})^{n}\delta_{q}%
(\hspace{0.01in}p)\nonumber\\
&  =2\sqrt{\Theta_{0}}\hspace{0.01in}\text{i}^{n}q^{-n(n\hspace{0.01in}%
+1)/2}[[\hspace{0.01in}n]]_{q}!\,p^{-n}\delta_{q}(\hspace{0.01in}p).
\end{align}
In the last step, we applied Eq.~(\ref{SkaDelFktEinQ}) of the previous chapter
and the following identities:%
\begin{equation}
(\hspace{0.01in}p^{-1}\Lambda\hspace{-0.02in}^{-1})^{n}=q^{-1-2-\ldots
-(n\hspace{0.01in}-1)}p^{-n}\Lambda\hspace{-0.02in}^{-n}=q^{-n(n\hspace
{0.01in}-1)/2}\hspace{0.01in}p^{-n}\Lambda\hspace{-0.02in}^{-n}.
\end{equation}

Next, we prove the formula in Eq.~(\ref{FouTraGanPotEindQ2}):%
\begin{align}
\mathcal{F}_{q}(x^{-n\hspace{0.01in}-1})(\hspace{0.01in}p)  &  =(-1)^{n}%
\hspace{0.01in}\frac{q^{n(n\hspace{0.01in}+1)/2}}{[[\hspace{0.01in}n]]_{q}%
!}\hspace{0.01in}\mathcal{F}_{q}(D_{q}^{n}x^{-1})(\hspace{0.01in}p)\nonumber\\
&  =\text{i}^{n}\hspace{0.01in}\frac{q^{n(n\hspace{0.01in}+1)/2}}%
{[[\hspace{0.01in}n]]_{q}!}\hspace{0.01in}(\Lambda\hspace{0.01in}p)^{n}%
\hspace{0.01in}\mathcal{F}_{q}(x^{-1})(\hspace{0.01in}p)\nonumber\\
&  =\text{i}^{n\hspace{0.01in}+1}\hspace{0.01in}\frac{\sqrt{\Theta_{0}}%
}{[[\hspace{0.01in}n]]_{q}!}\,p^{n}\operatorname*{sgn}(\hspace{0.01in}p).
\label{EinDimFourXInv}%
\end{align}
In the first step, we applied the result of the following calculation:%
\begin{align}
D_{q}^{k}z^{-1}  &  =[[-1]]_{q}\hspace{0.01in}[[-2]]_{q}\ldots\lbrack
\lbrack-k]]_{q}\,z^{-k\hspace{0.01in}-1}\nonumber\\
&  =(-1)^{k}q^{-k(k\hspace{0.01in}+1)/2}\hspace{0.01in}[[k]]_{q}%
!\,z^{-k\hspace{0.01in}-1}. \label{DarNegPotJacAblKooEinDim}%
\end{align}
The second step in Eq.~(\ref{EinDimFourXInv}) follows from
Eq.~(\ref{VerFouQ1DimAblOrt1}) and the last step from Eq.~(\ref{FouTraInvXEin}%
) if we take into account%
\begin{equation}
(\Lambda\hspace{0.01in}p)^{n}=q^{-n(n\hspace{0.01in}+1)/2}\hspace{0.01in}%
p^{n}\Lambda^{n}%
\end{equation}
and%
\begin{equation}
\Lambda\operatorname*{sgn}(\hspace{0.01in}p)=\operatorname*{sgn}%
(\hspace{0.01in}p).
\end{equation}

Finally, we calculate the $q$\textbf{-Fourier transform of the function
}$\Theta_{q}(z)$ \cite{olshanetsky1998q.alg}:%
\begin{align}
\mathcal{F}_{q}(\Theta_{q})  &  =\frac{1}{2\sqrt{\Theta_{0}}}\int
_{-x_{0}.\infty}^{\hspace{0.01in}x_{0}.\infty}\text{d}_{q}x\,\Theta_{q}%
(x)\exp_{q}(x|\text{i}p)\nonumber\\
&  =\frac{(1-q)\hspace{0.01in}x_{0}}{2\sqrt{\Theta_{0}}}\sum_{\varepsilon
\hspace{0.01in}=\hspace{0.01in}\pm}\hspace{0.02in}\sum_{m\hspace
{0.01in}=\hspace{0.01in}-\infty}^{\infty}q^{m}\hspace{0.01in}\Theta
_{q}(\varepsilon x_{0}q^{m})\exp_{q}(\varepsilon\hspace{0.01in}\text{i}%
\hspace{0.01in}x_{0}\hspace{0.01in}q^{m}p)\nonumber\\
&  =\frac{(1-q)\hspace{0.01in}x_{0}\hspace{0.01in}\Theta_{q}(x_{0})}%
{2\sqrt{\Theta_{0}}}\sum_{m\hspace{0.01in}=\hspace{0.01in}-\infty}^{\infty
}q^{m}\left[  \exp_{q}(\text{i}x_{0}\hspace{0.01in}q^{m}p)-\exp_{q}%
(-\text{i}x_{0}\hspace{0.01in}q^{m}p)\right] \nonumber\\
&  =\frac{\text{i}\hspace{0.01in}\Theta_{q}(x_{0})}{\sqrt{\Theta_{0}}}%
\hspace{0.01in}p^{-1}.
\end{align}
For the third step, we took into account that the function $\Theta_{q}(z)$ is
even and $q$-pe\-ri\-od\-ic [cf. Eq.~(\ref{PerTheFkt}) of
Chap.~\ref{KapQIntTrig}]. The last step results from the same considerations
we used in Eq.~(\ref{FouTraThePluMinEinQ}).

We take another look at Eq.~(\ref{ChaIdeqDelQKal1dim}) of the previous
chapter. This identity holds for functions being continuous at the origin. In
particular, we can apply it to functions having a power series expansion
around the origin. We can remove this restriction if we establish the
following identities for $n\in\mathbb{Z}$:%
\begin{equation}
\int\nolimits_{-\infty}^{\hspace{0.01in}\infty}\text{d}_{q}x\,\delta
_{q}(x)\hspace{0.02in}x^{n}=\left\{
\begin{tabular}
[c]{ll}%
$1,$ & f\"{u}r $n=0,$\\
$0,$ & f\"{u}r $n\neq0.$%
\end{tabular}
\ \right.  \label{ErwChaIdeEindQDelt}%
\end{equation}
This way, we can apply Eq.~(\ref{ChaIdeqDelQKal1dim}) of the previous chapter
to a function for which we have a Laurent series centered at the origin:%
\begin{equation}
\int\nolimits_{-\infty}^{\hspace{0.01in}\infty}\text{d}_{q}x\,\delta
_{q}(x)\hspace{0.01in}g(x)=g_{0}\quad\text{with}\quad g(z)=\sum_{k\hspace
{0.01in}=\hspace{0.01in}-\infty}^{\infty}g_{k}\hspace{0.01in}z^{k}.
\end{equation}
We aim to prove Eq.~(\ref{ErwChaIdeEindQDelt}). To this end, we compute the
$n$-th \textbf{Jackson derivative of }$\delta_{q}(x)$. Using
Eq.~(\ref{MehAnwJacAblFkt}) of Chap.~\ref{KapqAnaExpTriFkt} and
Eq.~(\ref{SkaDelFktEinQ}) of the previous chapter, we obtain ($n\in\mathbb{N}%
$):%
\begin{align}
D_{q}^{n}\hspace{0.01in}\delta_{q}(x)  &  =(-1)^{n}\hspace{0.01in}%
\frac{q^{-n(n\hspace{0.01in}+1)/2}}{(1-q)^{n}\hspace{0.01in}x^{n}}%
\sum_{m\hspace{0.01in}=\hspace{0.01in}0}^{n}%
\genfrac{[}{]}{0pt}{}{n}{m}%
_{q}(-q)^{m}q^{m(m\hspace{0.01in}-1)/2}\hspace{0.01in}\delta_{q}(x)\nonumber\\
&  =(-1)^{n}\hspace{0.01in}\frac{q^{n(n\hspace{0.01in}+1)/2}}{(1-q)^{n}%
\hspace{0.01in}x^{n}}\,(q\,;q)_{n}\,\delta_{q}(x)\nonumber\\
&  =(-1)^{n}q^{-n(n\hspace{0.01in}+1)/2}\hspace{0.01in}[[\hspace
{0.01in}n]]_{q}!\,x^{-n}\hspace{0.01in}\delta_{q}(x). \label{JacAblQDelN}%
\end{align}
In the second step, we used the Gaussian binomial formula \cite{Kac:2002eb}:%
\begin{equation}
(z\hspace{0.01in};q)_{n}=\sum_{m\hspace{0.01in}=\hspace{0.01in}0}^{n}%
\genfrac{[}{]}{0pt}{}{n}{m}%
_{q}(-z)^{m}q^{m(m\hspace{0.01in}-1)/2}.
\end{equation}
The final step in Eq.~(\ref{JacAblQDelN}) follows from
Eq.~(\ref{ZusqFakPocSym}) of Chap.~\ref{KapqAnaExpTriFkt}. Using the result in
Eq.~(\ref{JacAblQDelN}) and the rule in Eq.~(\ref{RegParIntBraLinK}) of
Chap.~\ref{KapQIntTrig}, we can perform the following computation:%
\begin{align}
\int_{-x_{0}.\infty}^{\hspace{0.01in}x_{0}.\infty}\text{d}_{q}x\,x^{-n}%
\hspace{0.01in}\delta_{q}(x)  &  =\frac{(-1)^{n}q^{n(n\hspace{0.01in}+1)/2}%
}{[[\hspace{0.01in}n]]_{q}!}\int_{-x_{0}.\infty}^{\hspace{0.01in}x_{0}.\infty
}\text{d}_{q}x\,D_{q}^{n}\delta_{q}(x)\nonumber\\
&  =\frac{1}{[[\hspace{0.01in}n]]_{q}!}\int_{-x_{0}.\infty}^{\hspace
{0.01in}x_{0}.\infty}\text{d}_{q}x\,\delta_{q}(q^{n}x)\cdot D_{q}^{n}1=0.
\end{align}

The results above are consistent with our earlier considerations. For example,
we can derive the expression for $\mathcal{F}_{q}(x^{n})$ in
Eq.~(\ref{FouTraGanPotEindQ}) using the derivative rule in
Eq.~(\ref{JacAblQDelN}) ($n\in\mathbb{N}_{0}$):%
\begin{align}
\mathcal{F}_{q}(x^{n})  &  =(-\text{i})^{n}D_{q,\hspace{0.01in}p}^{n}%
\hspace{0.01in}\mathcal{F}_{q}(1)(\hspace{0.01in}p)=2\sqrt{\Theta_{0}}%
\hspace{0.01in}(-\text{i})^{n}D_{q,\hspace{0.01in}p}^{n}\hspace{0.01in}%
\delta_{q}(\hspace{0.01in}p)\nonumber\\
&  =2\sqrt{\Theta_{0}}\hspace{0.01in}\text{i}^{n}q^{-n(n\hspace{0.01in}%
+1)/2}\hspace{0.01in}[[\hspace{0.01in}n]]_{q}!\,p^{-n}\delta_{q}%
(\hspace{0.01in}p). \label{QFouXn1Dim}%
\end{align}
With Eq.~(\ref{ErwChaIdeEindQDelt}) and the series expansions in
Eq.~(\ref{ReiEntqExpInv}), we obtain:%
\begin{align}
\mathcal{F}_{q}^{\hspace{0.01in}-1}(\hspace{0.01in}p^{-n}\delta_{q}%
(\hspace{0.01in}p))  &  =\frac{1}{2\sqrt{\Theta_{0}}}\int_{-p_{0}.\infty
}^{\hspace{0.01in}p_{0}.\infty}\text{d}_{q}p\,\delta_{q}(\hspace
{0.01in}p)\,p^{-n}\exp_{q^{-1}}(-\text{i}p|qx)\nonumber\\
&  =\frac{1}{2\sqrt{\Theta_{0}}}\hspace{0.01in}\frac{(-\text{i})^{n}%
q^{n(n\hspace{0.01in}+1)/2}x^{n}}{[[\hspace{0.01in}n]]_{q}!}.
\label{QFouXn1DimInv}%
\end{align}
Using the results of Eqs.~(\ref{QFouXn1Dim}) and (\ref{QFouXn1DimInv}), we
immediately check the following identity:%
\begin{equation}
(\mathcal{F}_{q}^{\hspace{0.01in}-1}\circ\mathcal{F}_{q})(x^{n})=x^{n}.
\label{InvFouEinXN}%
\end{equation}

Furthermore, we can rewrite the expression in Eq.~(\ref{DefVerDelFktQ1Dim}) of
the previous chapter by applying the derivative rule in Eq.~(\ref{JacAblQDelN}%
):%
\begin{equation}
\delta_{q}((\bar{\ominus}\,qa)\,\bar{\oplus}\,x)=\sum_{k\hspace{0.01in}%
=\hspace{0.01in}0}^{\infty}\frac{q^{k(k\hspace{0.01in}+1)/2}}{[[k]]_{q}%
!}(-a)^{k}\hspace{0.01in}(D_{q})^{k}\delta_{q}(x)=\sum_{k\hspace
{0.01in}=\hspace{0.01in}0}^{\infty}a^{k}\delta_{q}(x)\,x^{-k}.
\label{AltDar1DimTraDelFkt0}%
\end{equation}
In the same manner, we obtain:%
\begin{equation}
\delta_{q}(x\,\bar{\oplus}\,(\bar{\ominus}\,qa))=\sum_{k\hspace{0.01in}%
=\hspace{0.01in}0}^{\infty}x^{-k}\delta_{q}(x)\hspace{0.01in}a^{k}.
\label{AltDar1DimTraDelFkt}%
\end{equation}

The right-hand side of Eq.~(\ref{AltDar1DimTraDelFkt0}) or
Eq.~(\ref{AltDar1DimTraDelFkt}) satisfies Eq.~(\ref{ChaIdeDelQFkt1DimTra}) of
the previous chapter if the function $f$ has a power series expansion around
the origin. If should Eq.~(\ref{ChaIdeDelQFkt1DimTra}) of the last chapter
should hold when $f$ has a pole at the origin, we need to extend the summation
on the right-hand side of Eq.~(\ref{AltDar1DimTraDelFkt0}) or
Eq.~(\ref{AltDar1DimTraDelFkt}) to negative integers [cf.
Eq.~(\ref{ErwChaIdeEindQDelt})]:
\begin{align}
\delta_{q}((\bar{\ominus}\,qa)\,\bar{\oplus}\,x)  &  =\sum_{k\hspace
{0.01in}=\hspace{0.01in}-\infty}^{\infty}a^{k}\delta_{q}(x)\,x^{-k}%
,\nonumber\\
\delta_{q}(x\,\bar{\oplus}\,(\bar{\ominus}\,qa))  &  =\sum_{k\hspace
{0.01in}=\hspace{0.01in}-\infty}^{\infty}x^{-k}\delta_{q}(x)\hspace
{0.01in}a^{k}. \label{DelFktMeroFktEinDim}%
\end{align}
Using these expressions and taking into account Eq.~(\ref{ErwChaIdeEindQDelt}%
), we can perform the following calculation ($n\in\mathbb{N}$):%
\begin{gather}
\int_{-x_{0}.\infty}^{\hspace{0.01in}x_{0}.\infty}\text{d}_{q}x\,\delta
_{q}((\bar{\ominus}\,qa)\,\bar{\oplus}\,x)\,x^{-n}=\int_{-x_{0}.\infty
}^{\hspace{0.01in}x_{0}.\infty}\text{d}_{q}x\sum_{k\hspace{0.01in}%
=\hspace{0.01in}-\infty}^{\infty}a^{k}\hspace{0.01in}\delta_{q}(x)\,x^{-k}%
\hspace{0.01in}x^{-n}\nonumber\\
=\sum_{k\hspace{0.01in}=\hspace{0.01in}-\infty}^{\infty}a^{k}\int
_{-x_{0}.\infty}^{\hspace{0.01in}x_{0}.\infty}\text{d}_{q}x\,\delta
_{q}(x)\,x^{-(k+n)}=\sum_{k\hspace{0.01in}=\hspace{0.01in}-\infty}^{\infty
}a^{k}\hspace{0.01in}\delta_{k,-n}=a^{-n}.
\end{gather}

For completeness, we calculate the $q$\textbf{-Fourier transform of the
function }$u_{p}^{\ast}(x)$ [cf. Eq.~(\ref{ImpEigQDefEin})], using
Eq.~(\ref{QFouTraAlgFktEin}) and Eq.~(\ref{DelFktMeroFktEinDim}):%
\begin{align}
\mathcal{F}_{q}(u_{p}^{\ast}(x))  &  =\frac{1}{2\sqrt{\Theta_{0}}}%
\hspace{0.02in}\mathcal{F}_{q}(\exp_{q^{-1}}(-\text{i}p^{\prime}%
|\hspace{0.01in}qx))(\hspace{0.01in}p)\nonumber\\
&  =\sum_{n\hspace{0.01in}=\hspace{0.01in}0}^{\infty}\hspace{0.01in}%
\text{i}^{n}\frac{q^{-n(n\hspace{0.01in}+1)/2}}{p^{n}}\hspace{0.02in}%
\delta_{q}(\hspace{0.01in}p)\left.  \big(D_{q,\hspace{0.01in}x}^{n}%
\exp_{q^{-1}}(-\text{i}p^{\prime}|\hspace{0.01in}qx)\big)\right\vert
_{x\hspace{0.01in}=\hspace{0.01in}0}\nonumber\\
&  =\sum_{n\hspace{0.01in}=\hspace{0.01in}0}^{\infty}p^{-n}\delta_{q}%
(\hspace{0.01in}p)\hspace{0.01in}p^{\prime\hspace{0.01in}n}=\delta_{q}%
(\hspace{0.01in}p\,\bar{\oplus}\,(\bar{\ominus}\,qp^{\prime})).
\label{EinDimFourEbeWell}%
\end{align}
In the third step of the calculation above, we used the result of the
following calculation:%
\begin{align}
&  \left.  \big(D_{q,\hspace{0.01in}x}^{n}\exp_{q^{-1}}(-\text{i}p^{\prime
}|\hspace{0.01in}qx)\big)\right\vert _{x\hspace{0.01in}=\hspace{0.01in}0}%
=\sum_{k\hspace{0.01in}=\hspace{0.01in}0}^{\infty}\frac{q^{k(k\hspace
{0.01in}+1)/2}}{[[k]]_{q}!}\hspace{-0.02in}\left.  (-\text{i}p^{\prime}%
)^{k}\big(D_{q,\hspace{0.01in}x}^{n}\hspace{0.01in}x^{k}\big)\right\vert
_{x\hspace{0.01in}=\hspace{0.01in}0}\nonumber\\
&  \qquad\qquad=\sum_{k\hspace{0.01in}=\hspace{0.01in}0}^{\infty}%
\hspace{0.01in}\frac{q^{k(k\hspace{0.01in}+1)/2}}{[[k]]_{q}!}\hspace
{0.01in}(-\text{i}p^{\prime})^{k}\hspace{0.01in}[[\hspace{0.01in}n]]_{q}!%
\genfrac{[}{]}{0pt}{}{k}{n}%
_{q}\left.  x^{k-n}\right\vert _{x\hspace{0.01in}=\hspace{0.01in}0}\nonumber\\
&  \qquad\qquad=q^{n(n\hspace{0.01in}+1)/2}(-\text{i}p^{\prime})^{n}.
\end{align}

We can apply the reasoning above to the \textbf{inverse }$q$\textbf{-Fourier
transform} $\mathcal{F}_{q}^{-1}$. This way, we find:%
\begin{equation}
\mathcal{F}_{q}^{\hspace{0.01in}-1}(\delta_{q}(\hspace{0.01in}p))=\frac
{1}{2\sqrt{\Theta_{0}}}\int_{-p_{0}.\infty}^{\hspace{0.01in}p_{0}.\infty
}\text{d}_{q}p\,\delta_{q}(\hspace{0.01in}p)\exp_{q^{-1}}(-\text{i}%
p|qx)=\frac{1}{2\sqrt{\Theta_{0}}}.
\end{equation}
We rewrite the inverse $q$-Fourier transform by using $q$-in\-ver\-sions [cf.
Eq.~(\ref{ExpForInvQExp}) of Chap.~\ref{KapqAnaExpTriFkt}]:%
\begin{align}
\mathcal{F}_{q}^{\hspace{0.01in}-1}(f)  &  =\frac{1}{2\sqrt{\Theta_{0}}}%
\int\nolimits_{-p_{0}.\infty}^{p_{0}.\infty}\text{d}_{q}\hspace{0.01in}%
p\,\exp_{q}(\text{i}p|\bar{\ominus}\,qx)=\nonumber\\
&  =\frac{1}{2\sqrt{\Theta_{0}}}\,\underline{\bar{S}}_{\hspace{0.01in}%
x}\hspace{-0.03in}\left(  \int\nolimits_{-p_{0}.\infty}^{p_{0}.\infty}%
\text{d}_{q}\hspace{0.01in}p\,\exp_{q}(\text{i}p|qx)\right)  .
\end{align}
With the help of Eq.~(\ref{KonAntXn}) of Chap.~\ref{KapqAnaExpTriFkt} and a
procedure similar to that in Eq.~(\ref{FouTraThePluMinEinQ}), we get:%
\begin{align}
\mathcal{F}_{q}^{\hspace{0.01in}-1}(\theta_{q}^{+}-\hspace{0.01in}\theta
_{q}^{-})  &  =\frac{1}{2\sqrt{\Theta_{0}}}\,\underline{\bar{S}}%
_{\hspace{0.01in}x}\hspace{-0.03in}\left(  \int\limits_{0}^{p_{0}.\infty
}\text{d}_{q}\hspace{0.01in}p\,\exp_{q}(\text{i}p|qx)-\int\limits_{-p_{0}%
.\infty}^{0}\text{d}_{q}\hspace{0.01in}p\,\exp_{q}(\text{i}p|qx)\right)
\nonumber\\
&  =\frac{\text{i}(1-q)\hspace{0.01in}p_{0}}{q\sqrt{\Theta_{0}}}%
\,\underline{\bar{S}}_{\hspace{0.01in}x}\hspace{-0.03in}\left(  \sum
_{m\hspace{0.01in}=\hspace{0.01in}-\infty}^{\infty}q^{m}\sin_{q}%
(\hspace{0.01in}p_{0}\hspace{0.01in}q^{m}x)\right) \nonumber\\
&  =\underline{\bar{S}}_{\hspace{0.01in}x}\hspace{-0.03in}\left(
\frac{\text{i}}{q\sqrt{\Theta_{0}}\hspace{0.01in}x}\right)  =-\frac{\text{i}%
}{\sqrt{\Theta_{0}}}\hspace{0.01in}x^{-1}. \label{InvFourSgn}%
\end{align}
\textbf{ }In the same manner, we obtain [also cf. Eq.~(\ref{FouTrazInvEinQ1}%
)]:%
\begin{align}
\mathcal{F}_{q}^{\hspace{0.01in}-1}(\hspace{0.01in}p^{-1})  &  =\frac
{1}{2\sqrt{\Theta_{0}}}\,\underline{\bar{S}}_{\hspace{0.01in}x}\hspace
{-0.03in}\left(  \,\int\limits_{-p_{0}.\infty}^{p_{0}.\infty}\text{d}%
_{q}\hspace{0.01in}p\,p^{-1}\exp_{q}(\text{i}p|qx)\right) \nonumber\\
&  =\frac{\text{i}\hspace{0.01in}(1-q)}{\sqrt{\Theta_{0}}}\,\underline{\bar
{S}}_{\hspace{0.01in}x}\hspace{-0.03in}\left(  \sum_{m\hspace{0.01in}%
=\hspace{0.01in}-\infty}^{\infty}\sin_{q}(\hspace{0.01in}p_{0}\hspace
{0.01in}q^{m\hspace{0.01in}+1}x)\right) \nonumber\\
&  =\frac{\text{i}}{\sqrt{\Theta_{0}}}\,\underline{\bar{S}}_{\hspace{0.01in}%
x}(\Theta_{q}(\hspace{0.01in}p_{0}\hspace{0.01in}x))=-\text{i}\sqrt{\Theta
_{0}}\operatorname*{sgn}(x). \label{FouInvPinv}%
\end{align}
We apply this result to calculate $\mathcal{F}_{q}^{\hspace{0.01in}-1}(1)$
[also cf. Eq.~(\ref{FouTra1EinQ})]:%
\begin{align}
\mathcal{F}_{q}^{\hspace{0.01in}-1}(1)(x)  &  =\mathcal{F}_{q}^{\hspace
{0.01in}-1}(\hspace{0.01in}p\hspace{0.01in}p^{-1})(x)=\text{i}\hspace
{0.01in}q^{-1}\Lambda^{-1}D_{q,\hspace{0.01in}x}\hspace{0.01in}\mathcal{F}%
_{q}^{\hspace{0.01in}-1}(\hspace{0.01in}p^{-1})(x)\nonumber\\
&  =q^{-1}\sqrt{\Theta_{0}}\hspace{0.01in}\Lambda^{-1}D_{q,\hspace{0.01in}%
x}\operatorname*{sgn}(x)=2\sqrt{\Theta_{0}}\hspace{0.01in}\delta_{q}(x).
\label{EinInvFourEin}%
\end{align}
In the last step of the above calculation, we used Eqs.~(\ref{AblTheFktEin})
and (\ref{SkaDelFktEinQ}) from Chap.~\ref{KapQDisEin}. From
Eq.~(\ref{EinInvFourEin}), we derive an expression for $\mathcal{F}%
_{q}^{\hspace{0.01in}-1}(\hspace{0.01in}p^{n})$ with $n\in\mathbb{N}$, taking
into account Eq.~(\ref{VerFouQ1DimAblOrt2}):%
\begin{align}
\delta_{q}(x)  &  =\frac{1}{2\sqrt{\Theta_{0}}}\,\mathcal{F}_{q}%
^{\hspace{0.01in}-1}(1)(x)=\frac{1}{2\sqrt{\Theta_{0}}}\,\frac{\mathcal{F}%
_{q}^{\hspace{0.01in}-1}(D_{q,\hspace{0.01in}p}^{n}\hspace{0.01in}p^{n}%
)(x)}{[[\hspace{0.01in}n]]_{q}!}\nonumber\\
&  =\frac{1}{2\sqrt{\Theta_{0}}}\,\frac{\text{i}^{n}\mathcal{F}_{q}%
^{\hspace{0.01in}-1}(\hspace{0.01in}p^{n})(x)}{[[\hspace{0.01in}n]]_{q}%
!}\,x^{n}.
\end{align}
The result above leads to the following identity:%
\begin{equation}
\mathcal{F}_{q}^{\hspace{0.01in}-1}(\hspace{0.01in}p^{n})(x)=\text{i}%
^{-n}\hspace{0.01in}2\sqrt{\Theta_{0}}\,[[\hspace{0.01in}n]]_{q}%
!\,x^{-n}\hspace{0.01in}\delta_{q}(x). \label{InvFourEinPN}%
\end{equation}
For a function $f$ that has a power series expansion around the origin, the
formula above implies the following identity [also cf.
Eq.~(\ref{QFouTraAlgFktEin})]:%
\begin{equation}
\mathcal{F}_{q}^{\hspace{0.01in}-1}(f(\hspace{0.01in}p))(x)=2\sqrt{\Theta_{0}%
}\hspace{0.01in}\sum_{n\hspace{0.01in}=\hspace{0.01in}0}^{\infty}%
\hspace{0.01in}\text{i}^{-n}(D_{q,\hspace{0.01in}p}^{n}f)(0)\,\delta
_{q}(x)\hspace{0.01in}x^{-n}.
\end{equation}
Due to this identity, we can calculate the inverse $q$-Fourier transform of
the function $u_{p}(x)$:%
\begin{align}
\mathcal{F}_{q}^{\hspace{0.01in}-1}(u_{p}(x^{\prime}))(x)  &  =\sum
_{n\hspace{0.01in}=\hspace{0.01in}0}^{\infty}\hspace{0.01in}\text{i}%
^{-n}\left.  \big(D_{q,\hspace{0.01in}p}^{n}\exp_{q}(x^{\prime}|\hspace
{0.01in}\text{i}p)\big)\right\vert _{p=0}\delta_{q}(x)\hspace{0.01in}%
x^{-n}\nonumber\\
&  =\sum_{n\hspace{0.01in}=\hspace{0.01in}0}^{\infty}x^{\prime\hspace
{0.01in}n}\hspace{0.01in}\delta_{q}(x)\hspace{0.01in}x^{-n}=\delta_{q}%
((\bar{\ominus}\,qx^{\prime})\,\bar{\oplus}\,x).
\end{align}
Finally, we derive an expression for $\mathcal{F}_{q}^{\hspace{0.01in}%
-1}(\hspace{0.01in}p^{-n\hspace{0.01in}-1})$ with $n\in\mathbb{N}$. The
following calculation uses considerations similar to that in
Eq.~(\ref{EinDimFourXInv}):%
\begin{align}
\mathcal{F}_{q}^{\hspace{0.01in}-1}(\hspace{0.01in}p^{-n\hspace{0.01in}%
-1})(x)  &  =\frac{(-1)^{n}\hspace{0.01in}q^{n(n\hspace{0.01in}+1)/2}%
}{[[\hspace{0.01in}n]]_{q}!}\,\mathcal{F}_{q}^{\hspace{0.01in}-1}%
(D_{q,\hspace{0.01in}p}^{n}\hspace{0.01in}p^{-1})(x)\nonumber\\
&  =\frac{(-\text{i})^{n}\hspace{0.01in}q^{n(n\hspace{0.01in}+1)/2}}%
{[[\hspace{0.01in}n]]_{q}!}\,x^{n}\mathcal{F}_{q}^{\hspace{0.01in}-1}%
(\hspace{0.01in}p^{-1})(x)\nonumber\\
&  =\sqrt{\Theta_{0}}\,\frac{(-\text{i})^{n\hspace{0.01in}+1}\hspace
{0.01in}q^{n(n\hspace{0.01in}+1)/2}}{[[\hspace{0.01in}n]]_{q}!}\,x^{n}%
\operatorname*{sgn}(x). \label{InvFouPMinEin}%
\end{align}

With our findings, we can directly confirm the identities ($n\in\mathbb{N}$)%
\begin{align}
(\mathcal{F}_{q}\circ\mathcal{F}_{q}^{\hspace{0.01in}-1})(\delta_{q})  &
=\delta_{q},\nonumber\\
(\mathcal{F}_{q}\circ\mathcal{F}_{q}^{\hspace{0.01in}-1})(1)  &
=1,\nonumber\\
(\mathcal{F}_{q}\circ\mathcal{F}_{q}^{\hspace{0.01in}-1})(\hspace{0.01in}%
p^{n})  &  =p^{n},\nonumber\\
(\mathcal{F}_{q}\circ\mathcal{F}_{q}^{\hspace{0.01in}-1})(\hspace
{0.01in}p^{-n\hspace{0.01in}-1})  &  =p^{-n-1}, \label{NacInvFouDirEin}%
\end{align}
and%
\begin{align}
(\mathcal{F}_{q}^{\hspace{0.01in}-1}\circ\mathcal{F}_{q})(\delta_{q})  &
=\delta_{q},\nonumber\\
(\mathcal{F}_{q}^{\hspace{0.01in}-1}\circ\mathcal{F}_{q})(1)  &
=1,\nonumber\\
(\mathcal{F}_{q}^{\hspace{0.01in}-1}\circ\mathcal{F}_{q})(x^{n})  &
=x^{n},\nonumber\\
(\mathcal{F}_{q}^{\hspace{0.01in}-1}\circ\mathcal{F}_{q})(x^{-n\hspace
{0.01in}-1})  &  =x^{-n-1}. \label{NacInvFouDirEin2}%
\end{align}
To prove the last identity in Eq.~(\ref{NacInvFouDirEin})\ or
Eq.~(\ref{NacInvFouDirEin2}), we additionally need [cf.
Eq.~(\ref{FouTraThePluMinEinQ})]%
\begin{align}
\mathcal{F}_{q}(x^{n}\operatorname*{sgn}(x))(\hspace{0.01in}p)  &
=(-\text{i})^{n}D_{q,\hspace{0.01in}p}^{n}\hspace{0.01in}\mathcal{F}%
_{q}(\operatorname*{sgn}(x))(\hspace{0.01in}p)=\nonumber\\
&  =\frac{\text{i}^{n\hspace{0.01in}+1}\hspace{0.01in}[[\hspace{0.01in}%
n]]_{q}!}{\sqrt{\Theta_{0}}\hspace{0.01in}q^{n(n\hspace{0.01in}+1)/2}}%
\hspace{0.01in}p^{-n-1}%
\end{align}
or [cf. Eq.~(\ref{InvFourSgn})]%
\begin{align}
\mathcal{F}_{q}^{\hspace{0.01in}-1}(\hspace{0.01in}p^{n}\operatorname*{sgn}%
(\hspace{0.01in}p))(x)  &  =\text{i}^{n}q^{-n(n\hspace{0.01in}+1)/2}%
\Lambda^{-n}D_{q,\hspace{0.01in}x}^{n}\hspace{0.01in}\mathcal{F}_{q}%
^{\hspace{0.01in}-1}(\operatorname*{sgn}(\hspace{0.01in}p))(x)\nonumber\\
&  =\frac{(-\text{i})^{n\hspace{0.01in}+1}\hspace{0.01in}[[\hspace
{0.01in}n]]_{q}!}{\sqrt{\Theta_{0}}}\hspace{0.01in}x^{-n\hspace{0.01in}-1}.
\label{FourInvPnSgn}%
\end{align}
Note that we have confirmed $(\mathcal{F}_{q}^{\hspace{0.01in}-1}%
\circ\mathcal{F}_{q})(x^{n})=x^{n}$ using Eqs.~(\ref{QFouXn1Dim}) and
(\ref{QFouXn1DimInv}). To prove $(\mathcal{F}_{q}\circ\mathcal{F}_{q}%
^{\hspace{0.01in}-1})(\hspace{0.01in}p^{n})=p^{n}$ in the same way, we
additionally need:%
\begin{align}
\mathcal{F}_{q}(x^{-n}\hspace{0.01in}\delta_{q}(x))(\hspace{0.01in}p)  &
=\frac{1}{2\sqrt{\Theta_{0}}}\int_{-x_{0}.\infty}^{\hspace{0.01in}x_{0}%
.\infty}\text{d}_{q}x\,x^{-n}\delta_{q}(x)\hspace{0.01in}\exp_{q}%
(x|\text{i}p)\nonumber\\
&  =\frac{1}{2\sqrt{\Theta_{0}}}\hspace{0.01in}\frac{(\text{i}p)^{n}%
}{[[\hspace{0.01in}n]]_{q}!}.
\end{align}

\section{Fourier transforms on $\mathbb{G}_{q^{1/2},\,x_{0}}$%
\label{KapForqHalb}}

To obtain an improper $q$-de\-formed integral over the lattice $\mathbb{G}%
_{q^{1/2},\,x_{0}}$, we can add an improper $q$-in\-te\-gral over the lattice
$\mathbb{G}_{q,\hspace{0.01in}x_{0}}$ to an improper $q$-in\-te\-gral over the
lattice $\mathbb{G}_{q,\hspace{0.01in}q^{1/2}x_{0}}$ ($0<q<1$)
\cite{Cerchiai:1999,Hinterding:2000ph,Wess:math-ph9910013}:%
\begin{align}
&  \frac{1}{1+q^{1/2}}\left(  \int_{-x_{0}.\infty}^{\hspace{0.01in}%
x_{0}.\infty}\text{d}_{q}x\,f(x)+\int_{-x_{0}q^{1/2}.\infty}^{\hspace
{0.01in}x_{0}q^{1/2}.\infty}\text{d}_{q}x\,f(x)\right)  =\nonumber\\
&  \qquad\qquad=\frac{1-q}{1+q^{1/2}}\sum_{\varepsilon\hspace{0.01in}%
=\hspace{0.01in}\pm}^{\infty}\hspace{0.01in}\sum_{m\hspace{0.01in}%
=\hspace{0.01in}-\infty}^{\infty}\left[  x_{0}\hspace{0.01in}q^{m}%
f(\varepsilon x_{0}\hspace{0.01in}q^{m})+x_{0}\hspace{0.01in}q^{m\hspace
{0.01in}+1/2}f(\varepsilon x_{0}\hspace{0.01in}q^{m\hspace{0.01in}%
+1/2})\right] \nonumber\\
&  \qquad\qquad=(1-q^{1/2})\sum_{m\hspace{0.01in}=\hspace{0.01in}-\infty
}^{\infty}x_{0}q^{m/2}\left[  f(x_{0}\hspace{0.01in}q^{m/2})+f(-x_{0}%
\hspace{0.01in}q^{m/2})\right] \nonumber\\
&  \qquad\qquad=\int_{-x_{0}.\infty}^{\hspace{0.01in}x_{0}.\infty}%
\text{d}_{q^{1/2}}x\,f(x). \label{VerJacInt1Dim}%
\end{align}
Note that we have multiplied the sum of the two $q$-in\-te\-grals by the
factor $(1+q^{1/2})^{-1}$ to regain the ordinary integral for $q\rightarrow1$.

The $q^{1/2}$-in\-te\-gral in Eq.~(\ref{VerJacInt1Dim}) consists of
$q$-in\-te\-grals invariant under $q$-trans\-la\-tions. Thus, under
$q$-trans\-la\-tions the $q^{1/2}$-in\-te\-gral is also invariant. By analogy
with Eq.~(\ref{UnIntJacAbl}) from Chap.~\ref{KapQIntTrig}, the following holds
($m\in\mathbb{N}_{0}$):%
\begin{align}
&  \int_{-x.\infty}^{\hspace{0.01in}x.\infty}\text{d}_{q^{1/2}}z\,D_{q}%
^{m}f(z)=\nonumber\\
&  \qquad=\frac{1}{1+q^{1/2}}\left(  \int_{-x_{0}.\infty}^{\hspace
{0.01in}x_{0}.\infty}\text{d}_{q}\hspace{0.01in}x\,D_{q}^{m}f(x)+\int
_{-x_{0}q^{1/2}.\infty}^{\hspace{0.01in}x_{0}q^{1/2}.\infty}\text{d}%
_{q}\hspace{0.01in}x\,D_{q}^{m}f(x)\right) \nonumber\\
&  \qquad=0.
\end{align}
This identity, together with Eq.~(\ref{qTraVerGerKonUnk2}) of
Chap.~\ref{KapqAnaExpTriFkt}, implies:%
\begin{align}
&  \int_{-x.\infty}^{\hspace{0.01in}x.\infty}\text{d}_{q^{1/2}}z\,f(z\,\bar
{\oplus}\,a)=\int_{-x.\infty}^{\hspace{0.01in}x.\infty}\text{d}_{q^{1/2}}%
z\sum_{k\hspace{0.01in}=\hspace{0.01in}0}^{\infty}\frac{1}{[[k]]_{q}!}%
\hspace{0.01in}[D_{q}^{k}f(z)]\,a^{k}\nonumber\\
&  \qquad=\sum_{k\hspace{0.01in}=\hspace{0.01in}0}^{\infty}\frac{1}%
{[[k]]_{q}!}\int_{-x.\infty}^{\hspace{0.01in}x.\infty}\text{d}_{q^{1/2}%
}z\,[D_{q}^{k}f(z)]\,a^{k}=\int_{-x.\infty}^{\hspace{0.01in}x.\infty}%
\text{d}_{q^{1/2}}z\,f(z).
\end{align}
By similar reasoning, we obtain%
\begin{equation}
\int_{-x.\infty}^{\hspace{0.01in}x.\infty}\text{d}_{q^{1/2}}z\,f(a\,\bar
{\oplus}\,z)=\int_{-x.\infty}^{\hspace{0.01in}x.\infty}\text{d}_{q^{1/2}%
}z\,f(z).
\end{equation}

The \textbf{rule for integration by parts} also holds for the improper
$q^{1/2}$-in\-te\-gral. It gets this property from the $q$-in\-te\-grals of
which it is composed [cf. Eq.~(\ref{RegParIntBraLinK}) in
Chap.~\ref{KapQIntTrig})]:%
\begin{equation}
\int\text{d}_{q^{1/2}}x\,f(x)\,D_{q}^{k}g(x)=(-1)^{k}q^{-k(k\hspace
{0.01in}-1)/2}\int\text{d}_{q^{1/2}}x\left[  D_{q}^{k}\,f(x)\right]
g(q^{k}x). \label{ParQHalInt}%
\end{equation}

Next, we show how to formulate $q$-Fourier transforms using the $q^{1/2}%
$-in\-te\-grals from Eq.~(\ref{VerJacInt1Dim}). Instead of
Eqs.~(\ref{IntExpExpInv}) and (\ref{IntExpExpInv2}) of the last chapter, it
holds%
\begin{equation}
\int\limits_{-p_{0}.\infty}^{p_{0}.\infty}\text{d}_{q^{1/2}}p\,\exp_{q}%
(x_{0}|\text{i}p)\exp_{q^{-1}}(-\text{i}p|qx)=\left\{
\begin{array}
[c]{cc}%
\frac{\operatorname*{vol}\nolimits_{q}}{(1-q^{1/2})\hspace{0.01in}x_{0}%
}\vspace*{0.08in} & \text{for\quad}x=x_{0},\\
0 & \text{for\quad}x\neq x_{0},
\end{array}
\right.  \label{qIntProExpEinDim}%
\end{equation}
and%
\begin{equation}
\int\limits_{-x_{0}.\infty}^{x_{0}.\infty}\text{d}_{q^{1/2}}x\,\exp_{q^{-1}%
}(-\text{i}p_{0}|qx)\exp_{q}(x|\text{i}p)=\left\{
\begin{array}
[c]{cc}%
\frac{\operatorname*{vol}\nolimits_{q}}{(1-q^{1/2})\hspace{0.01in}p_{0}%
}\vspace*{0.08in} & \text{for\quad}p=p_{0},\\
0 & \text{for\quad}p\neq p_{0},
\end{array}
\right.  \label{qIntOrtExpEinDim}%
\end{equation}
with%
\begin{equation}
\operatorname*{vol}\nolimits_{q}=\frac{8\left[  \Theta_{q}(x_{0}%
\hspace{0.01in}p_{0})+\Theta_{q}(q^{1/2}x_{0}\hspace{0.01in}p_{0})\right]
}{(1+q^{1/2})^{2}}. \label{qDefVolEleEinDim}%
\end{equation}
Remarkably, we obtain the identities in Eqs.~(\ref{qIntProExpEinDim}) and
(\ref{qIntOrtExpEinDim}) from the identities in Eqs.~(\ref{IntExpExpInv}) and
(\ref{IntExpExpInv2}) of the previous chapter by the following substitutions:%
\begin{align}
\int_{-x_{0}.\infty}^{\hspace{0.01in}x_{0}.\infty}\text{d}_{q}x\rightarrow
\int_{-x_{0}.\infty}^{\hspace{0.01in}x_{0}.\infty}\text{d}_{q^{1/2}}x,\qquad
&  \int_{-p_{0}.\infty}^{\hspace{0.01in}p_{0}.\infty}\text{d}_{q}%
p\rightarrow\int_{-p_{0}.\infty}^{\hspace{0.01in}p_{0}.\infty}\text{d}%
_{q^{1/2}}p,\nonumber\\
4\Theta_{0}\rightarrow\operatorname*{vol}\nolimits_{q},\qquad &
1-q\rightarrow1-q^{1/2}. \label{UegRegQHal}%
\end{align}
Accordingly, the defining expression for the $q$-de\-formed Fourier transform
also changes [cf. Eq.~(\ref{DefQFouTraEin}) of the previous chapter]:%
\begin{align}
\mathcal{F}_{q^{1/2}}(\phi)(\hspace{0.01in}p)  &  =\operatorname*{vol}%
\nolimits_{q}^{-1/2}\int_{-x_{0}.\infty}^{\hspace{0.01in}x_{0}.\infty}%
\text{d}_{q^{1/2}}x\,\phi(x)\hspace{0.01in}\exp_{q}(x|\text{i}%
p),\nonumber\\[0.04in]
\mathcal{F}_{q^{1/2}}^{\hspace{0.01in}-1}(\psi)(x)  &  =\operatorname*{vol}%
\nolimits_{q}^{-1/2}\int_{-p_{0}.\infty}^{\hspace{0.01in}p_{0}.\infty}%
\text{d}_{q^{1/2}}p\,\psi(\hspace{0.01in}p)\hspace{0.01in}\exp_{q^{-1}%
}(-\text{i}p|qx). \label{DefQFouEinDimQHal}%
\end{align}
For all $n\in\mathbb{Z}$, it holds [cf. Eq.~(\ref{BasVekFktBraLin}) of
Chap.~\ref{KapQDisEin}]%
\begin{equation}
(\mathcal{F}_{q^{1/2}}^{\hspace{0.01in}-1}\circ\mathcal{F}_{q^{1/2}}%
)(\phi_{n/2}^{\pm})(x)=\left\{
\begin{array}
[c]{cc}%
1 & \text{for\quad}x=\pm\hspace{0.01in}x_{0}\hspace{0.01in}q^{n/2},\\
0 & \text{for\quad}x\neq\pm\hspace{0.01in}x_{0}\hspace{0.01in}q^{n/2},
\end{array}
\right.  \label{HalWerFouIdeEinDim}%
\end{equation}
and%
\begin{equation}
(\mathcal{F}_{q^{1/2}}\circ\mathcal{F}_{q^{1/2}}^{\hspace{0.01in}-1}%
)(\psi_{n/2}^{\pm})(\hspace{0.01in}p)=\left\{
\begin{array}
[c]{cc}%
1 & \text{for\quad}p=\pm\hspace{0.01in}p_{0}\hspace{0.01in}q^{n/2},\\
0 & \text{for\quad}p\neq\pm\hspace{0.01in}p_{0}\hspace{0.01in}q^{n/2}.
\end{array}
\right.  \label{HalWerFouIdeEinDim2}%
\end{equation}
We can prove Eq.~(\ref{HalWerFouIdeEinDim}) by the following calculation:%
\begin{align}
&  \frac{1}{\operatorname*{vol}\nolimits_{q}}\int_{-p_{0}.\infty}%
^{\hspace{0.01in}p_{0}.\infty}\text{d}_{q^{1/2}}p\,\int_{-x_{0}.\infty
}^{\hspace{0.01in}x_{0}.\infty}\text{d}_{q^{1/2}}\xi\,\phi_{n/2}^{\pm}%
(\xi)\hspace{0.01in}\exp_{q}(\xi|\text{i}p)\,\exp_{q^{-1}}(-\text{i}%
p|qx)\nonumber\\
&  \qquad=\frac{1-q^{1/2}}{\operatorname*{vol}\nolimits_{q}}\,x_{0}%
\hspace{0.01in}q^{n/2}\int_{-p_{0}.\infty}^{\hspace{0.01in}p_{0}.\infty
}\text{d}_{q^{1/2}}p\,\exp_{q}(\pm\hspace{0.01in}x_{0}\hspace{0.01in}%
q^{n/2}|\text{i}p)\,\exp_{q^{-1}}(-\text{i}p|qx)\nonumber\\
&  \qquad=\frac{1-q^{1/2}}{\operatorname*{vol}\nolimits_{q}}\,x_{0}%
\int_{-p_{0}.\infty}^{\hspace{0.01in}p_{0}.\infty}\text{d}_{q^{1/2}}%
p\,\exp_{q}(x_{0}|\text{i}p)\,\exp_{q^{-1}}(-\text{i}p|\hspace{-0.02in}%
\pm\hspace{-0.02in}q^{-n/2}qx).
\end{align}
Due to Eq.~(\ref{qIntProExpEinDim}), the result above implies
Eq.~(\ref{HalWerFouIdeEinDim}). Eq.~(\ref{HalWerFouIdeEinDim2}) can be proved
in a similar way.

The substitutions in Eq.~(\ref{UegRegQHal}) also change the normalization of
the \textbf{momentum eigenfunctions}:%
\begin{equation}
u_{p}(x)=\operatorname*{vol}\nolimits_{q}^{-1/2}\exp_{q}(\text{i}x|p),\qquad
u_{p}^{\ast}(x)=\operatorname*{vol}\nolimits_{q}^{-1/2}\exp_{q^{-1}}%
(-\text{i}p|qx).
\end{equation}
These momentum eigenfunctions again satisfy \textbf{completeness relations and
orthogonality conditions}. Due to the identities in
Eqs.~(\ref{qIntProExpEinDim}) and (\ref{qIntOrtExpEinDim}), it holds
($n,m\in\mathbb{Z}$; $\varepsilon,\varepsilon^{\prime}=\pm1$)%
\begin{align}
&  \int_{-p_{0}.\infty}^{\hspace{0.01in}p_{0}.\infty}\text{d}_{q^{1/2}%
}p\,u_{p}(\varepsilon\hspace{0.01in}q^{n/2})\,u_{p}^{\ast}(\varepsilon
^{\prime}q^{m/2})=\frac{\delta_{\varepsilon\varepsilon^{\prime}}\delta_{nm}%
}{(1-q^{1/2})\hspace{0.01in}x_{0}\hspace{0.01in}q^{n/2}}\nonumber\\
&  \qquad\qquad=\delta_{q^{1/2}}((\bar{\ominus}\,qx)\,\bar{\oplus
}\,y)|_{x\hspace{0.01in}=\hspace{0.01in}\varepsilon\hspace{0.01in}x_{0}%
q^{n/2};\,y\hspace{0.01in}=\hspace{0.01in}\varepsilon^{\prime}x_{0}q^{m/2}}
\label{QVolRelNeuNorExp1Dim}%
\end{align}
and%
\begin{align}
&  \int_{-x_{0}.\infty}^{\,x_{0}.\infty}\text{d}_{q^{1/2}}x\,u_{\varepsilon
\hspace{0.01in}q^{n/2}}^{\ast}(x)\,u_{\varepsilon^{\prime}q^{m/2}}%
(x)=\frac{\delta_{\varepsilon\varepsilon^{\prime}}\delta_{nm}}{(1-q^{1/2}%
)\hspace{0.01in}p_{0}\hspace{0.01in}q^{n/2}}\nonumber\\
&  \qquad\qquad=\delta_{q^{1/2}}((\bar{\ominus}\,qp)\,\bar{\oplus}\,p^{\prime
})|_{p\hspace{0.01in}=\hspace{0.01in}\varepsilon\hspace{0.01in}p_{0}%
q^{n/2};\,p^{\prime}=\hspace{0.01in}\varepsilon^{\prime}p_{0}q^{m/2}},
\label{QOrtRelNeuNorExp1Dim}%
\end{align}
where $\delta_{q^{1/2}}$ denotes the delta distribution for the $q^{1/2}%
$-in\-te\-gral:%
\begin{equation}
\int_{-x_{0}.\infty}^{\hspace{0.01in}x_{0}.\infty}\text{d}_{q^{1/2}}%
x\,\delta_{q^{1/2}}(x)\,f(x)=f(0). \label{ChaIdeQHalDelFktEinDim}%
\end{equation}
Note that the last identity in Eq.~(\ref{QVolRelNeuNorExp1Dim}) or
Eq.~(\ref{QOrtRelNeuNorExp1Dim}) holds due to Eq.~(\ref{VerQDelqHalBasFkt}).

We can apply the substitutions in Eq.~(\ref{UegRegQHal}) to
Eqs.~(\ref{FouDarFkt1DimQ})-(\ref{InvFourTraInd}) of the previous chapter.The
derivations of these identities remain valid if we change the improper Jackson
integrals and the normalization factors according to Eq.~(\ref{UegRegQHal}).

For completeness, we calcuate the \textbf{classical limit of the }%
$q$\textbf{-de\-formed volume element} introduced in
Eq.~(\ref{qDefVolEleEinDim}). Using Eq.~(\ref{TheQFktKlaGre}) of
Chap.~\ref{KapQIntTrig}, we find:%
\begin{align}
\lim_{q\hspace{0.01in}\rightarrow1^{-}}\operatorname*{vol}\nolimits_{q}  &
=\lim_{q\hspace{0.01in}\rightarrow1^{-}}\frac{8\left[  \Theta_{q}(x_{0}%
\hspace{0.01in}p_{0})+\Theta_{q}(q^{1/2}x_{0}\hspace{0.01in}p_{0})\right]
}{(1+q^{1/2})^{2}}\nonumber\\
&  =\frac{8\left(  \frac{\pi}{2}+\frac{\pi}{2}\right)  }{4}=2\pi.
\label{KlaGreVolEle1Dim}%
\end{align}

We summarize some more \textbf{properties of the }$q^{1/2}$\textbf{-delta
function}. First, we calculate the $q$-derivatives of $\delta_{q^{1/2}}(z)$ by
using Eq.~(\ref{ParQHalInt}):%
\begin{align}
\big(D_{q}^{k}\delta_{q^{1/2}}\big)(\hspace{0.01in}g)  &  =\int_{-x_{0}%
.\infty}^{\hspace{0.01in}x_{0}.\infty}\text{d}_{q^{1/2}}x\,\big[D_{q}%
^{k}\delta_{q^{1/2}}(x)\big]g(x)\nonumber\\
&  =(-1)^{k}q^{-k(k\hspace{0.01in}+1)/2}\int_{-x_{0}.\infty}^{\hspace
{0.01in}x_{0}.\infty}\text{d}_{q^{1/2}}x\,\delta_{q^{1/2}}(x)\,(D_{q}%
^{k}g)(q^{-k}x)\nonumber\\
&  =(-1)^{k}q^{-k(k\hspace{0.01in}+1)/2}(D_{q}^{k}g)(0).
\end{align}
Furthermore, $\delta_{q^{1/2}}$ is again related to the $q$-derivative of a
step function on the $q^{1/2}$-lattice:%
\begin{equation}
D_{q}\hspace{0.01in}\theta_{q^{1/2}}^{\pm}(x)=\frac{\pm\hspace{0.01in}%
2}{1+q^{1/2}}\,\delta_{q^{1/2}}(x)\quad\text{with}\quad\theta_{q^{1/2}}^{\pm
}(x)=\sum_{n\hspace{0.01in}=\hspace{0.01in}-\infty}^{\infty}\phi_{n/2}^{\pm
}(x). \label{ZusTheDelQEinDimqHal}%
\end{equation}
Note that Eq.~(\ref{ZusTheDelQEinDimqHal}) is a distributional identity. If
$g$ is a test function that is continuous at the origin of the lattice
$\mathbb{G}_{q^{1/2},\,x_{0}}$, it holds:%
\begin{equation}
\big(D_{q}\hspace{0.01in}\theta_{q^{1/2}}^{\pm}\big)(\hspace{0.01in}%
g)=\frac{\pm\hspace{0.01in}2}{1+q^{1/2}}\,g(0^{\pm})=\frac{\pm\hspace
{0.01in}2}{1+q^{1/2}}\,\delta_{q^{1/2}}(g). \label{AblStuFktQHal}%
\end{equation}
We prove this relation by the following calculation (we restrict ourselves to
$\theta_{q^{1/2}}^{+}$; there is a similar calculation for $\theta_{q^{1/2}%
}^{-}$):%
\begin{align}
&  \big(D_{q}\hspace{0.01in}\theta_{q^{1/2}}^{+}\big)(\hspace{0.01in}%
g)=\int_{-x_{0}.\infty}^{\hspace{0.01in}x_{0}.\infty}\text{d}_{q^{1/2}%
}x\,\big[D_{q}\hspace{0.01in}\theta_{q^{1/2}}^{+}(x)\big]g(x)\nonumber\\
&  \qquad\quad=\frac{-q^{-1}}{1+q^{1/2}}\left(  \int_{0}^{\hspace{0.01in}%
x_{0}.\infty}\text{d}_{q}x\,(D_{q}\hspace{0.01in}g)(q^{-1}x)+\int_{0}%
^{\hspace{0.01in}x_{0}q^{1/2}.\infty}\text{d}_{q}x\,(D_{q}\hspace
{0.01in}g)(q^{-1}x)\right) \nonumber\\
&  \qquad\quad=\frac{-1}{1+q^{1/2}}\left(  -g(0^{+})-g(0^{+})\right)
=\frac{2}{1+q^{1/2}}\,g(0^{+}).
\end{align}
In the second step of the calculation above, we used Eq.~(\ref{AblqDisHer})
from Chap.~\ref{KapQDisEin}\ and Eq.~(\ref{VerJacInt1Dim}). The penultimate
step follows from Eq.~(\ref{PosParIntAblq1Dim}) of Chap.~\ref{KapQIntTrig} and
the fact that $g$ vanishes at infinity.

From Eq.~(\ref{ZusTheDelQEinDimqHal}), we obtain a representation of
$\delta_{q^{1/2}}^{\pm}(z)$ as the limit of regular distributions [also see
the calculation in Eq.~(\ref{GreRegDisqDis}) of Chap.~\ref{KapQDisEin}]
($0<q<1$):%
\begin{align}
\delta_{q^{1/2}}(z)  &  =\pm\frac{1+q^{1/2}}{2}\,D_{q}\hspace{0.01in}%
\theta_{q^{1/2}}^{\pm}(z)=\pm\frac{1+q^{1/2}}{2}\sum_{n\hspace{0.01in}%
=\hspace{0.01in}-\infty}^{\infty}D_{q}\hspace{0.01in}\phi_{n/2}^{\pm
}(z)\nonumber\\
&  =\pm\frac{1+q^{1/2}}{2}\sum_{n\hspace{0.01in}=\hspace{0.01in}-\infty
}^{\infty}\frac{\phi_{n/2}^{\pm}(z)-\phi_{n/2}^{\pm}(qz)}{(1-q)z}\nonumber\\
&  =\pm\sum_{n\hspace{0.01in}=\hspace{0.01in}-\infty}^{\infty}\frac{\phi
_{n/2}^{\pm}(z)-\phi_{(n\hspace{0.01in}-1)/2}^{\pm}(z)+\phi_{(n\hspace
{0.01in}-1)/2}^{\pm}(z)-\phi_{(n-2)/2}^{\pm}(z)}{2\hspace{0.01in}(1-q^{1/2}%
)z}\nonumber\\
&  =\pm(1-q^{1/2})^{-1}z^{-1}\lim_{n\hspace{0.01in}\rightarrow\hspace
{0.01in}\infty}\phi_{n/2}^{\pm}(z).
\end{align}
Note that the identities above refer to distributions. Hence, for any test
function $g$ on the lattice $\mathbb{G}_{q^{1/2},\,x_{0}}$, we have:%
\begin{equation}
\int\nolimits_{-x_{0}.\infty}^{\hspace{0.01in}x_{0}.\infty}\text{d}_{q^{1/2}%
}x\,\delta_{q^{1/2}}(x)\,g(x)=\pm\lim_{n\hspace{0.01in}\rightarrow
\hspace{0.01in}\infty}\int\nolimits_{-x_{0}.\infty}^{\hspace{0.01in}%
x_{0}.\infty}\text{d}_{q^{1/2}}x\,\frac{\phi_{n/2}^{\pm}(x)}{(1-q^{1/2}%
)\hspace{0.01in}x}\,g(x).
\end{equation}

We examine the \textbf{scaling properties of the}\textit{ }$q^{1/2}%
$\textbf{-step functions and the }$q^{1/2}$\textbf{-delta distribution}. From
the definition of the scaling operator, it follows ($k,n\in\mathbb{Z}$):%
\begin{equation}
\Lambda^{k/2}\phi_{n/2}^{\pm}(z)=\phi_{n/2}^{\pm}(q^{k/2}z)=\phi
_{(n-k)/2}^{\pm}(z).
\end{equation}
Moreover, we have:%
\begin{align}
\Lambda^{k/2}\hspace{0.01in}\theta_{q^{1/2}}^{\pm}(z)  &  =\sum_{n\hspace
{0.01in}=\hspace{0.01in}-\infty}^{\infty}\Lambda^{k}\hspace{0.01in}\phi
_{n/2}^{\pm}(z)=\sum_{n\hspace{0.01in}=\hspace{0.01in}-\infty}^{\infty}%
\phi_{(n-k)/2}^{\pm}(z)\nonumber\\
&  =\sum_{n\hspace{0.01in}=\hspace{0.01in}-\infty}^{\infty}\phi_{n/2}^{\pm
}(z)=\theta_{q^{1/2}}^{\pm}(z).
\end{align}
Using the results above, we can derive the scaling property of $\delta
_{q^{1/2}}$. With the help of Eq.~(\ref{ZusTheDelQEinDimqHal}), we do the
following calculation:%
\begin{align}
\delta_{q^{1/2}}(q^{k/2}z)  &  =\Lambda^{k/2}\delta_{q^{1/2}}(z)=\pm
\frac{1+q^{1/2}}{2}\,\Lambda^{k/2}D_{q}\hspace{0.01in}\theta_{q^{1/2}}^{\pm
}(z)\nonumber\\
&  =\pm\frac{1+q^{1/2}}{2}\,q^{-k/2}D_{q}\hspace{0.01in}\Lambda^{k/2}%
\theta_{q^{1/2}}^{\pm}(z)\nonumber\\
&  =\pm\frac{1+q^{1/2}}{2}\,q^{-k/2}D_{q}\hspace{0.01in}\theta_{q^{1/2}}^{\pm
}(z)=q^{-k/2}\delta_{q^{1/2}}(z).
\end{align}

We modify the identities of Eqs.~(\ref{DefVerDelFktQ1Dim}%
)-(\ref{ChaIdeDelQFkt1DimTra}) of Chap.~\ref{KapQDisEin} so that they apply to
the lattice $\mathbb{G}_{q^{1/2},\,x_{0}}$. To this end, we apply the
substitutions in Eq.~(\ref{UegRegQHal}) to these identities and write
$\delta_{q^{1/2}}$ instead of $\delta_{q}$. For example, we get%
\begin{align}
\int_{-x_{0}.\infty}^{\hspace{0.01in}x_{0}.\infty}\text{d}_{q^{1/2}}%
x\,\delta_{q^{1/2}}((\bar{\ominus}\,qa)\,\bar{\oplus}\,x)\,f(x)  &
=f(a),\nonumber\\[0.04in]
\int_{-x_{0}.\infty}^{\hspace{0.01in}x_{0}.\infty}\text{d}_{q^{1/2}}%
x\hspace{0.01in}f(x)\,\delta_{q^{1/2}}(x\,\bar{\oplus}\,(\bar{\ominus}\,qa))
&  =f(a) \label{ChaIdeQHalbDeltFkt}%
\end{align}
where%
\begin{gather}
\delta_{q^{1/2}}(x\,\bar{\oplus}\,(\bar{\ominus}\,qa))=\delta_{q^{1/2}%
}(x\,\bar{\oplus}\,(\bar{\ominus}\,qa))\nonumber\\
=\sum_{k\hspace{0.01in}=\hspace{0.01in}0}^{\infty}\frac{q^{k(k\hspace
{0.01in}+1)/2}(-a)^{k}}{[[k]]_{q}!}D_{q}^{k}\delta_{q^{1/2}}(x).
\end{gather}
Moreover, we have ($m\in\mathbb{Z}$, $0<q<1$):%
\begin{gather}
\left.  \delta_{q^{1/2}}((\bar{\ominus}\,qa)\,\bar{\oplus}\,x)\right\vert
_{a\hspace{0.01in}=\hspace{0.01in}\pm x_{0}q^{m/2}}=\left.  \delta_{q^{1/2}%
}(x\,\bar{\oplus}\,(\bar{\ominus}\,qa))\right\vert _{a\hspace{0.01in}%
=\hspace{0.01in}\pm x_{0}q^{m/2}}\nonumber\\
=\frac{\phi_{m/2}^{\pm}(x)}{(1-q^{1/2})\hspace{0.01in}x_{0}\hspace
{0.01in}q^{m/2}}. \label{VerQDelqHalBasFkt}%
\end{gather}

In the following, we give \textbf{expressions for }$q$\textbf{-de\-formed
Fourier transforms on the lattice }$\mathbb{G}_{q^{1/2},\,x_{0}}$. In complete
analogy to Eq.~(\ref{FouTrazInvEinQ1}) of the previous chapter, we can
calculate $\mathcal{F}_{q^{1/2}}(x^{-1})$ with the following result:%
\begin{align}
\mathcal{F}_{q^{1/2}}(x^{-1})(\hspace{0.01in}p)  &  =\operatorname*{vol}%
\nolimits_{q}^{-1/2}\int_{-x_{0}.\infty}^{\hspace{0.01in}x_{0}.\infty}%
\text{d}_{q^{1/2}}x\,x^{-1}\exp_{q}(x|\text{i}p)\nonumber\\
&  =\frac{2\text{i}\left[  \Theta_{q}(x_{0}\hspace{0.01in}p)+\Theta
_{q}(q^{1/2}x_{0}\hspace{0.01in}p)\right]  }{(1+q^{1/2})\operatorname*{vol}%
\nolimits_{q}^{1/2}}.
\end{align}
For the lattice $\mathbb{G}_{q^{1/2},\,x_{0}}$, the variable $p$ can only take
the values $\pm\hspace{0.01in}x_{0}\hspace{0.01in}p_{0}\hspace{0.01in}q^{m/2}$
with $m\in\mathbb{Z}$. Thus, it holds [cf. Eq.~(\ref{qDefVolEleEinDim})]:%
\begin{align}
&  \frac{2\text{i}\left[  \Theta_{q}(\pm\hspace{0.01in}x_{0}\hspace
{0.01in}p_{0}\hspace{0.01in}q^{m/2})+\Theta_{q}(\pm\hspace{0.01in}x_{0}%
\hspace{0.01in}p_{0}\hspace{0.01in}q^{(m\hspace{0.01in}+1)/2})\right]
}{(1+q^{1/2})\operatorname*{vol}\nolimits_{q}^{1/2}}=\pm\frac{2\text{i}\left[
\Theta_{q}(x_{0}\hspace{0.01in}p_{0})+\Theta_{q}(q^{1/2}x_{0}\hspace
{0.01in}p_{0})\right]  }{(1+q^{1/2})\operatorname*{vol}\nolimits_{q}^{1/2}%
}\nonumber\\
&  \qquad\qquad\qquad=\pm\frac{\text{i}}{\sqrt{2}}\big[\Theta_{q}(x_{0}%
\hspace{0.01in}p_{0})+\Theta_{q}(q^{1/2}x_{0}\hspace{0.01in}p_{0}%
)\big]^{1/2}\nonumber\\
&  \qquad\qquad\qquad=\pm\frac{\text{i}}{4}(1+q^{1/2})\operatorname*{vol}%
\nolimits_{q}^{1/2}.
\end{align}
Combining our results so far, we finally obtain:%
\begin{equation}
\mathcal{F}_{q^{1/2}}(x^{-1})(\hspace{0.01in}p)=\frac{\text{i}}{4}%
(1+q^{1/2})\operatorname*{vol}\nolimits_{q}^{1/2}\operatorname*{sgn}%
(\hspace{0.01in}p)\text{.} \label{QHalFouTraInvKooEinDim}%
\end{equation}

The expression for $\mathcal{F}_{q^{1/2}}(x^{-1})$ enables us to calculate an
expression for $\mathcal{F}_{q^{1/2}}(1)$. Analogous to Eq.~(\ref{FouTra1EinQ}%
) of the last chapter and considering Eq.~(\ref{ZusTheDelQEinDimqHal}), we
find:%
\begin{align}
\mathcal{F}_{q^{1/2}}(1)(\hspace{0.01in}p)  &  =\mathcal{F}_{q^{1/2}}%
(x\hspace{0.01in}x^{-1})(\hspace{0.01in}p)=-\text{i}D_{q,\hspace{0.01in}%
p}\hspace{0.01in}\mathcal{F}_{q^{1/2}}(x^{-1})(\hspace{0.01in}p)\nonumber\\
&  =\frac{1}{4}(1+q^{1/2})\,\operatorname*{vol}\nolimits_{q}^{1/2}%
D_{q,\hspace{0.01in}p}\big[\theta_{q^{1/2}}^{+}(\hspace{0.01in}p)-\theta
_{q^{1/2}}^{-}(\hspace{0.01in}p)\big]\nonumber\\
&  =\operatorname*{vol}\nolimits_{q}^{1/2}\delta_{q^{1/2}}(\hspace{0.01in}p).
\label{QHalFouTraIdeEinDim}%
\end{align}
This result implies the following identities for the shifted $q^{1/2}$-delta
function:%
\begin{align}
\delta_{q^{1/2}}(a\,\bar{\oplus}\,x)  &  =\frac{1}{\operatorname*{vol}%
\nolimits_{q}^{1/2}}\hspace{0.01in}\mathcal{F}_{q^{1/2}}(1)(a\,\bar{\oplus
}\,x)\nonumber\\
&  =\frac{1}{\operatorname*{vol}\nolimits_{q}}\int_{-p_{0}.\infty}%
^{\hspace{0.01in}p_{0}.\infty}\text{d}_{q^{1/2}}p\,\exp_{q}(\text{i}%
p|a\,\bar{\oplus}\,x)\nonumber\\
&  =\frac{1}{\operatorname*{vol}\nolimits_{q}}\int_{-p_{0}.\infty}%
^{\hspace{0.01in}p_{0}.\infty}\text{d}_{q^{1/2}}p\,\exp_{q}(\text{i}%
p|a)\exp_{q}(\text{i}p|\hspace{0.01in}x)\nonumber\\
&  =\frac{1}{\operatorname*{vol}\nolimits_{q}^{1/2}}\hspace{0.01in}%
\mathcal{F}_{q^{1/2}}(\exp_{q}(\text{i}p|a))(x).
\end{align}

Next, we come to the $q^{1/2}$-Fourier transform of $\theta_{q^{1/2}}%
^{+}-\hspace{0.01in}\theta_{q^{1/2}}^{-}$. Analogous to
Eq.~(\ref{FouTraThePluMinEinQ}) of the previous chapter, we do the following
calculation:%
\begin{align}
&  \mathcal{F}_{q^{1/2}}\big(\theta_{q^{1/2}}^{+}-\hspace{0.01in}%
\theta_{q^{1/2}}^{-}\big)(\hspace{0.01in}p)=\nonumber\\
&  \qquad=\operatorname*{vol}\nolimits_{q}^{-1/2}\left(  \int\limits_{0}%
^{x_{0}.\infty}\text{d}_{q^{1/2}}x\,\exp_{q}(x|\text{i}p)-\int\limits_{-x_{0}%
.\infty}^{0}\text{d}_{q^{1/2}}x\,\exp_{q}(-x|\text{i}p)\right) \nonumber\\
&  \qquad=\frac{(1-q^{1/2})\hspace{0.01in}x_{0}}{\operatorname*{vol}%
\nolimits_{q}^{1/2}}\sum_{m\hspace{0.01in}=\hspace{0.01in}-\infty}^{\infty
}q^{m/2}\left[  \exp_{q}(\text{i}x_{0}\hspace{0.01in}q^{m/2}p)-\exp
_{q}(-\text{i}x_{0}\hspace{0.01in}q^{m/2}p)\right] \nonumber\\
&  \qquad=\frac{2\text{i}(1-q^{1/2})\hspace{0.01in}x_{0}}{\operatorname*{vol}%
\nolimits_{q}^{1/2}}\sum_{m\hspace{0.01in}=\hspace{0.01in}-\infty}^{\infty
}q^{m/2}\sin_{q}(x_{0}\hspace{0.01in}q^{m/2}p). \label{FouTraQHalStufDif}%
\end{align}
We rewrite the series in the last expression in the following way [cf.
Eq.~(\ref{FouTraThePluMinEinQ}) of Chap.~\ref{KapEinDimQFouTraN}]:%
\begin{align}
&  \sum_{m\hspace{0.01in}=\hspace{0.01in}-\infty}^{\infty}q^{m/2}\sin
_{q}(x_{0}\hspace{0.01in}q^{m/2}p)=\sum_{m\hspace{0.01in}=\hspace
{0.01in}-\infty}^{\infty}q^{m}\left[  \sin_{q}(x_{0}\hspace{0.01in}%
q^{m}p)+q^{1/2}\sin_{q}(x_{0}\hspace{0.01in}q^{m\hspace{0.01in}+1/2}p)\right]
\nonumber\\
&  \qquad=\frac{(1-q)\hspace{0.01in}x_{0}\hspace{0.01in}p}{(1-q)\hspace
{0.01in}x_{0}\hspace{0.01in}p}\lim_{M\rightarrow\hspace{0.01in}\infty}%
\sum_{m\hspace{0.01in}=\hspace{0.01in}-M}^{\infty}\left[  q^{m}\sin_{q}%
(x_{0}\hspace{0.01in}q^{m}p)+q^{m+1/2}\sin_{q}(x_{0}\hspace{0.01in}%
q^{m\hspace{0.01in}+1/2}p)\right] \nonumber\\
&  \qquad=\frac{1}{(1-q)\hspace{0.01in}x_{0}\hspace{0.01in}p}\lim
_{M\rightarrow\hspace{0.01in}\infty}\left[  1-\cos_{q}(x_{0}\hspace
{0.01in}q^{-M}p)+1-\cos_{q}(x_{0}\hspace{0.01in}q^{-M+1/2}p)\right]
\nonumber\\
&  \qquad=\frac{2}{(1-q)\hspace{0.01in}x_{0}\hspace{0.01in}p}.
\end{align}
Plugging this result into Eq.~(\ref{FouTraQHalStufDif}), we finally obtain:%
\begin{equation}
\mathcal{F}_{q^{1/2}}\big(\theta_{q^{1/2}}^{+}-\hspace{0.01in}\theta_{q^{1/2}%
}^{-}\big)(\hspace{0.01in}p)=\frac{4\hspace{0.01in}\text{i}}{(1+q^{1/2}%
)\operatorname*{vol}\nolimits_{q}^{1/2}}\,p^{-1}. \label{FouQHalStuDifEnd}%
\end{equation}

Using the formula above and Eq.~(\ref{ZusTheDelQEinDimqHal}), we can
immediately calculate the $q^{1/2}$-Fourier transform of $\delta_{q^{1/2}}$:%
\begin{align}
\mathcal{F}_{q^{1/2}}(\delta_{q^{1/2}})(\hspace{0.01in}p)  &  =\frac
{1+q^{1/2}}{4}\,\mathcal{F}_{q^{1/2}}\big(D_{q}\big(\theta_{q^{1/2}}%
^{+}-\hspace{0.01in}\theta_{q^{1/2}}^{-}\big)\big)\nonumber\\
&  =-\frac{\text{i}\hspace{0.01in}(1+q^{1/2})}{4}\hspace{0.01in}\Lambda
\hspace{0.01in}p\hspace{0.01in}\mathcal{F}_{q^{1/2}}\big(\theta_{q^{1/2}}%
^{+}-\hspace{0.01in}\theta_{q^{1/2}}^{-}\big)=\operatorname*{vol}%
\nolimits_{q}^{-1/2}.
\end{align}
Furthermore, we can give the $q^{1/2}$-Fourier transform of $\theta_{q}^{+}$
and $\theta_{q}^{-}$. Taking into account%
\begin{equation}
2\hspace{0.01in}\theta_{q^{1/2}}^{\pm}=\pm\big(\theta_{q^{1/2}}^{+}%
-\hspace{0.01in}\theta_{q^{1/2}}^{-}\big)+1
\end{equation}
and Eqs.~(\ref{QHalFouTraIdeEinDim}) and (\ref{FouQHalStuDifEnd}), it follows:%
\begin{equation}
\mathcal{F}_{q^{1/2}}\big(\theta_{q^{1/2}}^{\pm}\big)=\frac{\pm\hspace
{0.01in}2\hspace{0.01in}\text{i}}{(1+q^{1/2})\operatorname*{vol}%
\nolimits_{q}^{1/2}}\,p^{-1}+\frac{1}{2}\operatorname*{vol}\nolimits_{q}%
^{1/2}\delta_{q^{1/2}}(\hspace{0.01in}p).
\end{equation}

We also determine an expression for $\mathcal{F}_{q^{1/2}}(x^{k})$ with
$k\in\mathbb{Z}$. We first consider the case of the positive exponents. Using
Eq.~(\ref{DarIdeJaclKoo}) from the previous chapter and
Eq.~(\ref{QHalFouTraIdeEinDim}), we perform the following calculation
($n\in\mathbb{N}_{0}$):%
\begin{align}
\delta_{q^{1/2}}(\hspace{0.01in}p)  &  =\operatorname*{vol}\nolimits_{q}%
^{-1/2}\mathcal{F}_{q^{1/2}}(1)(\hspace{0.01in}p)=\operatorname*{vol}%
\nolimits_{q}^{-1/2}\frac{\mathcal{F}_{q^{1/2}}(D_{q}^{n}x^{n})(\hspace
{0.01in}p)}{[[\hspace{0.01in}n]]_{q}!}\nonumber\\
&  =\operatorname*{vol}\nolimits_{q}^{-1/2}\frac{(-\text{i})^{n}}%
{[[\hspace{0.01in}n]]_{q}!}\hspace{0.01in}(\Lambda\hspace{0.01in}p)^{n}%
\hspace{0.01in}\mathcal{F}_{q^{1/2}}(x^{n})(\hspace{0.01in}p).
\end{align}
Solving for $\mathcal{F}_{q^{1/2}}(x^{n})$ gives:%
\begin{align}
\mathcal{F}_{q^{1/2}}(x^{n})(\hspace{0.01in}p)  &  =\text{i}^{n}%
\operatorname*{vol}\nolimits_{q}^{1/2}[[\hspace{0.01in}n]]_{q}!\,(\Lambda
\hspace{0.01in}p)^{-n}\hspace{0.01in}\delta_{q^{1/2}}(\hspace{0.01in}%
p)\nonumber\\
&  =\text{i}^{n}\operatorname*{vol}\nolimits_{q}^{1/2}q^{-n(n\hspace
{0.01in}+1)/2}\hspace{0.01in}[[\hspace{0.01in}n]]_{q}!\,p^{-n}\hspace
{0.01in}\delta_{q^{1/2}}(\hspace{0.01in}p). \label{FouTraqEinDimPotX}%
\end{align}
From this result, we obtain a formula to calculate the $q^{1/2}$-Fourier
transform of a function that has a power series about the origin:%
\begin{equation}
\mathcal{F}_{q^{1/2}}(f)(\hspace{0.01in}p)=\operatorname*{vol}\nolimits_{q}%
^{1/2}\delta_{q}(\hspace{0.01in}p)\sum_{n\hspace{0.01in}=\hspace{0.01in}%
0}^{\infty}\text{i}^{n}q^{-n(n\hspace{0.01in}+1)/2}\hspace{0.01in}%
p^{-n}\hspace{0.01in}(D_{q}^{n}f)(0).
\end{equation}
Finally, we treat the case of the negative exponents. Using
Eq.~(\ref{DarNegPotJacAblKooEinDim}) of the previous chapter and
Eq.~(\ref{QHalFouTraInvKooEinDim}), we obtain:%
\begin{align}
\mathcal{F}_{q^{1/2}}(x^{-n\hspace{0.01in}-1})(\hspace{0.01in}p)  &
=\frac{(-1)^{n}q^{n(n\hspace{0.01in}+1)/2}}{[[\hspace{0.01in}n]]_{q}%
!}\,\mathcal{F}_{q^{1/2}}(D_{q}^{n}x^{-1})(\hspace{0.01in}p)\nonumber\\
&  =\frac{\text{i}^{n}q^{n(n\hspace{0.01in}+1)/2}}{[[\hspace{0.01in}n]]_{q}%
!}\,(\Lambda\hspace{0.01in}p)^{n}\hspace{0.01in}\mathcal{F}_{q}(x^{-1}%
)(\hspace{0.01in}p)\nonumber\\
&  =\frac{\text{i}^{n\hspace{0.01in}+1}(1+q^{1/2})}{4\hspace{0.01in}%
[[\hspace{0.01in}n]]_{q}!}\hspace{0.01in}\operatorname*{vol}\nolimits_{q}%
^{1/2}p^{n}\operatorname*{sgn}(\hspace{0.01in}p).
\end{align}

Analogous to the previous chapter, we derive \textbf{expressions for the
inverse }$q^{1/2}$\textbf{-Fourier transforms}. This way, we have [cf.
Eq.~(\ref{InvFourSgn}) of Chap.~\ref{KapEinDimQFouTraN} and
Eq.~(\ref{VerJacInt1Dim})]%
\begin{align}
&  \mathcal{F}_{q^{1/2}}^{\hspace{0.01in}-1}(\theta_{q^{1/2}}^{+}%
-\hspace{0.01in}\theta_{q^{1/2}}^{-})(x)=\nonumber\\
&  \qquad=\operatorname*{vol}\nolimits_{q}^{-1/2}\,\underline{\bar{S}%
}_{\hspace{0.01in}x}\hspace{-0.02in}\left(  \int\limits_{0}^{p_{0}.\infty
}\text{d}_{q^{1/2}}p\,\exp_{q}(\text{i}p|qx)-\int\limits_{-p_{0}.\infty}%
^{0}\text{d}_{q^{1/2}}p\,\exp_{q}(\text{i}p|qx)\right) \nonumber\\
&  \qquad=-\frac{4\hspace{0.01in}\text{i}}{(1+q^{1/2})\operatorname*{vol}%
\nolimits_{q}^{1/2}x}%
\end{align}
and\footnote{Note that $\operatorname{sgn}=\theta_{q^{1/2}}^{+}-\hspace
{0.01in}\theta_{q^{1/2}}^{-}$.} [cf. Eq.~(\ref{FourInvPnSgn}) of the previous
chapter]%
\begin{align}
\mathcal{F}_{q^{1/2}}^{-1}(\hspace{0.01in}p^{n}\operatorname{sgn}%
(\hspace{0.01in}p))(x)  &  =\text{i}^{n}q^{-n(n\hspace{0.01in}+1)/2}%
\hspace{0.01in}\Lambda^{-n}\hspace{0.01in}D_{q,\hspace{0.01in}x}^{n}%
\hspace{0.01in}\mathcal{F}_{q^{1/2}}^{\hspace{0.01in}-1}(\operatorname{sgn}%
(\hspace{0.01in}p))(x)\nonumber\\
&  =(-\text{i})^{n\hspace{0.01in}+1}\hspace{0.01in}[[\hspace{0.01in}%
n]]_{q}!\hspace{0.01in}\frac{4}{(1+q^{1/2})\operatorname*{vol}\nolimits_{q}%
^{1/2}}\hspace{0.01in}x^{-n\hspace{0.01in}-1}.
\end{align}
Moreover, it holds [cf. Eq.~(\ref{FouInvPinv}) of the previous chapter]%
\begin{align}
\mathcal{F}_{q^{1/2}}^{\hspace{0.01in}-1}(\hspace{0.01in}p^{-1})(x)  &
=\operatorname*{vol}\nolimits_{q}^{-1/2}\,\underline{\bar{S}}_{\hspace
{0.01in}x}\hspace{-0.02in}\left(  \,\int\limits_{-p_{0}.\infty}^{p_{0}.\infty
}\text{d}_{q^{1/2}}\hspace{0.01in}p\,p^{-1}\exp_{q}(\text{i}p|qx)\right)
\nonumber\\
&  =-\frac{\text{i}}{4}\,(1+q^{1/2})\operatorname*{vol}\nolimits_{q}%
^{1/2}\operatorname*{sgn}(x) \label{InvFouQHalbpInv}%
\end{align}
and [cf. Eq.~(\ref{EinInvFourEin}) of the previous chapter]%
\begin{align}
\mathcal{F}_{q^{1/2}}^{\hspace{0.01in}-1}(1)(x)  &  =\mathcal{F}_{q^{1/2}%
}^{\hspace{0.01in}-1}(\hspace{0.01in}p\hspace{0.01in}p^{-1})(x)=\text{i}%
q^{-1}\Lambda^{-1}D_{q,\hspace{0.01in}x}\hspace{0.01in}\mathcal{F}_{q^{1/2}%
}^{\hspace{0.01in}-1}(\hspace{0.01in}p^{-1})(x)\nonumber\\
&  =\operatorname*{vol}\nolimits_{q}^{1/2}\delta_{q^{1/2}}(x).
\label{InvFouQHalbIde}%
\end{align}
We can use Eq.~(\ref{InvFouQHalbpInv}) to derive an expression for
$\mathcal{F}_{q^{1/2}}^{\hspace{0.01in}-1}(\hspace{0.01in}p^{-n\hspace
{0.01in}-1})$ with $n\in\mathbb{N}$ [cf. Eq.~(\ref{InvFouPMinEin}) of the
previous chapter]. This way, we get:%
\begin{equation}
\mathcal{F}_{q^{1/2}}^{-1}(\hspace{0.01in}p^{-n\hspace{0.01in}-1}%
)(x)=\frac{(1+q^{1/2})\operatorname*{vol}\nolimits_{q}^{1/2}}{4}%
\,\frac{(-\text{i})^{n\hspace{0.01in}+1}\hspace{0.01in}q^{n(n\hspace
{0.01in}+1)/2}}{[[\hspace{0.01in}n]]_{q}!}\,x^{n}\operatorname*{sgn}(x).
\end{equation}
Similarly, we obtain an expression for $\mathcal{F}_{q^{1/2}}^{\hspace
{0.01in}-1}(\hspace{0.01in}p^{n})$ with $n\in\mathbb{N}$ by using
Eq.~(\ref{InvFouQHalbIde}) [cf. Eq.~(\ref{InvFourEinPN}) of the previous
chapter]:%
\begin{equation}
\mathcal{F}_{q^{1/2}}^{\hspace{0.01in}-1}(\hspace{0.01in}p^{n})(x)=\text{i}%
^{-n}\operatorname*{vol}\nolimits_{q}^{1/2}\,[[\hspace{0.01in}n]]_{q}%
!\,x^{-n}\hspace{0.01in}\delta_{q^{1/2}}(x).
\end{equation}
Consequently, we have the following formula for a function $f$ being analytic
at the origin:%
\begin{equation}
\mathcal{F}_{q^{1/2}}^{\hspace{0.01in}-1}(f(\hspace{0.01in}%
p))(x)=\operatorname*{vol}\nolimits_{q}^{1/2}\delta_{q^{1/2}}(x)\hspace
{0.01in}\sum_{n\hspace{0.01in}=\hspace{0.01in}0}^{\infty}\hspace
{0.01in}\text{i}^{-n}(D_{q,\hspace{0.01in}p}^{n}f)(0)\,x^{-n}.
\end{equation}

In Eqs.~(\ref{NacInvFouDirEin}) and (\ref{NacInvFouDirEin2}) of the previous
chapter, we can replace $q$ by $q^{1/2}$ due to the identities%
\begin{align}
\mathcal{F}_{q^{1/2}}^{\hspace{0.01in}-1}(\delta_{q^{1/2}}(\hspace
{0.01in}p))(x)  &  =\operatorname*{vol}\nolimits_{q}^{-1/2},\nonumber\\
\mathcal{F}_{q^{1/2}}^{\hspace{0.01in}-1}(\hspace{0.01in}p^{-n}\hspace
{0.01in}\delta_{q^{1/2}}(\hspace{0.01in}p))(x)  &  =\operatorname*{vol}%
\nolimits_{q}^{-1/2}q^{n(n\hspace{0.01in}+1)/2}\,\frac{(-\text{i}x)^{n}%
}{[[\hspace{0.01in}n]]_{q}!},
\end{align}
and ($n\in\mathbb{N}$)%
\begin{align}
\mathcal{F}_{q^{1/2}}(x^{-n}\hspace{0.01in}\delta_{q^{1/2}}(x))(\hspace
{0.01in}p)  &  =\operatorname*{vol}\nolimits_{q}^{-1/2}\frac{(\text{i}p)^{n}%
}{[[\hspace{0.01in}n]]_{q}!},\nonumber\\
\mathcal{F}_{q^{1/2}}(x^{n}\operatorname*{sgn}(x))(\hspace{0.01in}p)  &
=\frac{\text{i}^{n+1}\hspace{0.01in}4\hspace{0.01in}[[\hspace{0.01in}n]]_{q}%
!}{(1+q^{1/2})\operatorname*{vol}\nolimits_{q}^{1/2}q^{n(n\hspace
{0.01in}+1)/2}}\hspace{0.01in}p^{-n\hspace{0.01in}-1}.
\end{align}

\bibliographystyle{unsrt}
\bibliography{acompat,habil}

\end{document}